\documentclass[letter,11pt]{amsart}
\usepackage{amsfonts,textcomp,amssymb,amsmath,amsthm}
\usepackage{color,ulem}
\usepackage{fullpage}
\usepackage{graphicx}
\usepackage[active]{srcltx}
\usepackage{mathabx}
\usepackage{dsfont,mathrsfs,stmaryrd,wasysym,amsbsy}
\usepackage{enumerate}
\usepackage{IEEEtrantools}

\usepackage{tikz}
\usetikzlibrary{positioning}

\newcommand{\esssupmath}[1]{{\underset{#1}{\textrm{\rm ess sup}} }}

\newcommand{\cF}{\mathcal{F}}

\newcommand{\cT}{\mathcal{T}}

\newcommand{\CC}{\mathbb{C}}

\newcommand{\EE}{\mathbb{E}}
\newcommand{\E}{\mathbb{E}}

\newcommand{\NN}{\mathbb{N}}

\newcommand{\PP}{\mathbb{P}}

\newcommand{\RR}{\mathbb{R}}
\newcommand{\rr}{\mathbb{R}}

\newcommand{\bone}{\mathbf{1}}

\theoremstyle{plain}
\newtheorem{theorem}{Theorem}[section]
\newtheorem{corollary}[theorem]{Corollary}
\newtheorem{lemma}[theorem]{Lemma}
\newtheorem{proposition}[theorem]{Proposition}

\newtheorem{ass}[theorem]{Assumption}

\theoremstyle{definition}
\newtheorem{remark}[theorem]{Remark}

\numberwithin{equation}{section}

\begin{document}
\title[Stefan problem with surface tension]{Stefan problem with surface tension: uniqueness of physical solutions under radial symmetry} 
\author{Yucheng Guo, Sergey Nadtochiy and Mykhaylo Shkolnikov}
\address{ORFE Department, Princeton University, Princeton, NJ 08544.} 
\email{yg7348@princeton.edu}
\address{Department of Applied Mathematics, Illinois Institute of Technology, Chicago, IL 60616.}
\email{snadtochiy@iit.edu}
\address{ORFE Department, Bendheim Center for Finance, and Program in Applied \& Computational Mathematics, Princeton University, Princeton, NJ 08544.}
\email{mshkolni@gmail.com}
\footnotetext[1]{S.~Nadtochiy is partially supported by the NSF CAREER grant DMS-1651294.}
\footnotetext[2]{M.~Shkolnikov is partially supported by the NSF grant DMS-2108680.}

\begin{abstract}
We study the Stefan problem with surface tension and radially symmetric initial data. In this context, the notion of a so-called physical solution, which exists globally despite the inherent  blow-ups of the melting rate, has been recently introduced in \cite{NaShsurface}.~The paper at hand is devoted to the proof that the physical solution is unique, the first such result when the free boundary is not flat, or when two phases are present.~The main argument relies on a detailed analysis of the hitting probabilities for a three-dimensional Brownian motion, as well as on a novel convexity property of the free boundary obtained by comparison techniques.~In the course of the proof, we establish a wide variety of regularity estimates for the free boundary and for the temperature function, of interest in their own right.
\end{abstract}

\maketitle


\section{Introduction}\label{se:intro}

The classical formulation of the Stefan problem with surface tension (also known as the Stefan-Gibbs-Thomson problem) reads as follows.~Given $\Gamma_0\subset\rr^3$ and $u(0,\cdot)\!:\rr^3\to\rr$, find $\{\Gamma_t\subset\rr^3\}_{t>0}$ and $\{u(t,\cdot)\!:\rr^3\to\rr\}_{t>0}$ such that
\begin{eqnarray}
&& \partial_t u(t,z) = \frac{1}{2}\Delta_z u(t,z),\quad z\in\rr^3\backslash\partial\Gamma_t,\;\;t>0,  \label{intro: heat eq} \\ 
&& u(t,z) = \gamma H_t(z),\quad z\in\partial\Gamma_t,\;\;t>0, \label{intro: GT} \\
&& V_t(z)=\frac{1}{2}(\nabla_z u(t,z))\cdot\overrightarrow{n_{+,t}}(z)
+\frac{1}{2}(\nabla_z u(t,z))\cdot\overrightarrow{n_{-,t}}(z),\quad z\in\partial\Gamma_t,\;\; t\ge0.
\label{intro: Stefan}
\end{eqnarray}
Here, for each $t\ge0$, the sets $\Gamma_t$ and $\rr^3\backslash\Gamma_t$ are to be interpreted as the regions occupied by a solid (e.g., ice) and the corresponding liquid (e.g., water) at the time $t$; the function $-u(t,\cdot)$ describes the temperature relative to the equilibrium freezing point in the two regions; $H_t$ stands for the mean curvature (``surface tension'') of the interface $\partial\Gamma_t$ (with the convention $H_t\ge0$ if $\Gamma_t$ is convex); $\gamma>0$ is the surface tension parameter; $\overrightarrow{n_{+,t}}$ and $\overrightarrow{n_{-,t}}$ are the unit normals along $\partial\Gamma_t$ pointing towards the liquid and the solid; and $V_t$ is the rate of freezing in the direction $\overrightarrow{n_{+,t}}$. The equations \eqref{intro: heat eq}, \eqref{intro: GT}, and \eqref{intro: Stefan} are referred to as the (standard) heat equation, the Gibbs-Thomson boundary condition, and the Stefan growth condition, respectively.

\medskip

The flat version of the problem \eqref{intro: heat eq}--\eqref{intro: Stefan} in which $\partial\Gamma_t=\{\Lambda_t\}\times\rr^2$ is a hyperplane for all $t$, so that \eqref{intro: GT} reduces to a zero boundary condition, has been introduced independently by \textsc{Lam\'{e}} and \textsc{Clapeyron} \cite{LC} and by \textsc{Stefan} \cite{Stefan1}, \cite{Stefan2}, \cite{Stefan3}, \cite{Stefan4}.~In non-flat settings, the correction \eqref{intro: GT} to the Dirichlet boundary condition postulates that the equilibrium freezing temperature at a positively curved interface is depressed  proportionally to the mean curvature (a so-called Gibbs-Thomson effect, see, e.g., \cite[Section 8.3]{Gli} for more details). One consequence of the presence of mean curvature in \eqref{intro: GT} is that the sign of the temperature cannot be determined by the phase, e.g., some parts of the liquid may be ``supercooled", with temperatures below the regular (flat) equilibrium value.~This lack of consistency between the temperature sign and the phase causes serious mathematical challenges (e.g., the lack of a comparison principle) and is the reason why a global well-posedness theory for such equations is still missing.~A concrete example of a challenge caused by supercooling is that \eqref{intro: heat eq}--\eqref{intro: Stefan} exhibits blow-ups of the melting rate~$-V_t$ (see \cite{meir} for a deliberation in the radially symmetric setting).\footnote{On the other hand, if the supercooling is imposed via the initial condition, the presence of surface tension in \eqref{intro: GT} is expected to improve the (space) regularity of the moving interface, thus giving hope for well-posedness of such Stefan problems. (This is in contrast to the supercooled Stefan problem with $\gamma=0$, which is believed to be ill-posed in spatial dimension above one, for generic initial data.)}~Due to such blow-ups, previous mathematical work focused on (i) \textit{classical} solutions to \eqref{intro: heat eq}--\eqref{intro: Stefan} under special circumstances: a small $\gamma$ (see \cite{FrRe}), a small initial data (see \cite{eps}), a nearly flat initial interface (see \cite{HaGu}), an initial interface close to a steady sphere (see \cite{ha}), and (ii) \textit{weak} solutions to \eqref{intro: heat eq}--\eqref{intro: Stefan} with generic initial conditions (see \cite{Lu}, \cite[Example 5]{RoSa} for constructions of a weak solution with a sharp interface, and \cite{vis} for a prior construction with a possibly non-sharp interface).~The notion of a weak solution is, however, too weak to ensure uniqueness (see \cite[Section 5]{Lu}).~Instead, we adapt (stronger) \textit{physical} solutions to \eqref{intro: heat eq}--\eqref{intro: Stefan} in the radially symmetric setting, proposed in \cite{NaShsurface}, and show their uniqueness.

\medskip

Suppose now that $\Gamma_0$ is a ball of radius $\Lambda_0>0$ around the origin and that $u(0,z)=v(0,|z|)$ for some $v(0,\cdot)\!:[0,\infty)\to\rr$.~Then, one is naturally led to look for solutions of \eqref{intro: heat eq}--\eqref{intro: Stefan} in which for all $t>0$, the set $\Gamma_t$ is a ball of radius $\Lambda_t$ around the origin and $u(t,z)=v(t,|z|)$ for some $v(t,\cdot)\!:[0,\infty)\to\rr$.~Since in this case $\gamma H_t(z)=\gamma/|z|$, it is convenient to reformulate \eqref{intro: heat eq}--\eqref{intro: Stefan} in terms of $\{\Lambda_t\}_{t\ge0}$ and $\{w(t,\cdot):=\gamma/\cdot-v(t,\cdot)\!:(0,\infty)\to\rr\}_{t\ge0}$, obtaining, respectively,
\begin{eqnarray}
&& \partial_t w(t,x) = \frac{1}{2}\partial_{xx}w(t,x)+\frac{1}{x}\partial_x w(t,x),\quad x\in(r_0,\infty)\backslash\{\Lambda_t\},\;\; t>0, \label{intro: heat rs} \\
&& w(t,\Lambda_t)=w(t,r_0)=0,\quad t>0, \label{intro: GT rs} \\
&& \Lambda_t'=-\frac{1}{2}\partial_xw(t,x+)|_{x=\Lambda_t}+\frac{1}{2}\partial_x w(t,x-)|_{x=\Lambda_t},\quad t\ge0. \label{intro: Stefan rs}
\end{eqnarray}
In \eqref{intro: heat rs}--\eqref{intro: Stefan rs}, we have let the solid have a ``cold core'':~a ball with a radius $r_0\!>\!0$ at whose boundary the Gibbs-Thomson temperature of $-\gamma/r_0$ is maintained.~This feature rules out the non-physical temperatures below absolute zero and, more broadly, an application of the Gibbs-Thomson relationship beyond its range of validity (cf.~\cite[footnote on p.~200]{Gli}) together with the mathematical anomalies that go with it.

\medskip

Due to the possible blow-ups of the melting rate $-\Lambda'$ we adapt the notion of a \textit{probabilistic} solution to \eqref{intro: heat rs}--\eqref{intro: Stefan rs} from \cite{NaShsurface}.~A right-continuous function with left limits $\Lambda\!:[0,\zeta)\!\to\!(r_0,\infty)$ is called a probabilistic solution to \eqref{intro: heat rs}--\eqref{intro: Stefan rs} with initial data $(\Lambda_{0-}\!>\! r_0,\,w(0-,\cdot)\!:(r_0,\infty)\!\to\!\rr)$ if
\begin{equation}
\frac{\Lambda_{0-}^3}{3}-\frac{\Lambda_t^3}{3}
=\int_{r_0}^\infty w(0-,x)\,\PP^x(\tau\le t\wedge\tau_{r_0})\,\nu(\mathrm{d}x),\quad t\in[0,\zeta), \label{intro: growth weak} 
\end{equation}
where under $\PP^x$ the canonical process $R$ is a three-dimensional Bessel process started from $x$, $\tau:=\inf\{s\ge0\!:(R_s-\Lambda_s)(R_0-\Lambda_{0-})\le0\}$, $\tau_{r_0}:=\inf\{s\ge0\!:R_s\le r_0\}$, and $\nu(\mathrm{d}x):=x^2\,\mathrm{d}x$.~Hereby, the terminal time $\zeta$ is to be understood as the time of instantaneous melting to the cold core.~Note that \eqref{intro: heat rs} can be rewritten as 
\begin{equation}\label{forwardPDE}
\partial_t\big(x^2 w(t,x)\big) 
=\partial_x\,\Big(\frac{1}{2}\partial_x\big(x^2 w(t,x)\big)
-\frac{1}{x}\big(x^2 w(t,x)\big)\!\Big)
=\Big(\frac{1}{2}\partial_{xx}+\frac{1}{x}\partial_x\Big)^*\big(x^2 w(t,x)\big),
\end{equation}
and the Fokker-Planck equation (formally) leads to the equality of measures
\begin{equation}\label{intro: def w}
x^2w(t,x)\,\mathrm{d}x=\int_{r_0}^\infty  y^2 w(0-,y)\,\PP^y\big(R_t\,\mathbf{1}_{\{t<\tau\wedge\tau_{r_0}\}}\in\cdot\big)\,\mathrm{d}y
\quad\text{on}\quad(r_0,\infty),
\end{equation}
which implicitly encodes \eqref{intro: heat rs}--\eqref{intro: GT rs}.~Moreover,~\eqref{intro: growth weak} follows from \eqref{intro: Stefan rs} via the (formal) computation
\begin{equation*}
\begin{split}
-\frac{\mathrm{d}}{\mathrm{d}t} \frac{\Lambda_t^3}{3}
&=\frac{x^2}{2}\partial_x w(t,x+)\Big|_{x=\Lambda_t}
-\frac{x^2}{2}\partial_x w(t,x-)\Big|_{x=\Lambda_t} \\
&=\Big(\frac{1}{2}\partial_x\big(x^2 w(t,x+)\big)
-\frac{1}{x}\big(x^2 w(t,x+)\big)\!\Big)\Big|_{x=\Lambda_t}
-\Big(\frac{1}{2}\partial_x\big(x^2 w(t,x-)\big)
-\frac{1}{x}\big(x^2 w(t,x-)\big)\!\Big)\Big|_{x=\Lambda_t} \\
&=\frac{\mathrm{d}}{\mathrm{d}t}\int_{r_0}^\infty x^2 w(0-,x)\,\PP^x(\tau\le t\wedge\tau_{r_0})\,\mathrm{d}x.
\end{split}
\end{equation*}
The latter equality can be obtained by solving the system of linear equations 
\begin{equation*}
\begin{split}
& \int_{r_0}^\infty x^2 w(0-,x)\,\frac{\mathrm{d}}{\mathrm{d}t}\PP^x(\tau\le t\!\wedge\!\tau_{r_0})\,\mathrm{d}x
+ \int_{r_0}^\infty x^2 w(0-,x)\,\frac{\mathrm{d}}{\mathrm{d}t}\PP^x(\tau_{r_0}\le t\!\wedge\!\tau)\,\mathrm{d}x = -\frac{\mathrm{d}}{\mathrm{d}t}\int_{r_0}^\infty x^2 w(t,x)\,\mathrm{d}x, \\
& \int_{r_0}^\infty x^2 w(0-,x)\frac{1}{\Lambda_t}
\frac{\mathrm{d}}{\mathrm{d}t}\PP^x(\tau\!\le\! t\!\wedge\!\tau_{r_0})\,\mathrm{d}x
\!+\!\int_{r_0}^\infty x^2 w(0-,x)\frac{1}{r_0}\frac{\mathrm{d}}{\mathrm{d}t}
\PP^x(\tau_{r_0}\!\le\! t\!\wedge\!\tau)\,\mathrm{d}x = -\frac{\mathrm{d}}{\mathrm{d}t}\int_{r_0}^\infty x w(t,x)\,\mathrm{d}x
\end{split}
\end{equation*}
(both make use of \eqref{intro: def w}, and the second one is a consequence of $\frac{\mathrm{d}}{\mathrm{d}t}\EE^x[1/R_{t\wedge\tau\wedge\tau_{r_0}}]=0$), in which we evaluate the right-hand sides by means of \eqref{intro: GT rs}, \eqref{forwardPDE} and \eqref{intro: GT rs}, $\partial_t(xw(t,x))=\frac{1}{2}\partial_{xx}(xw(t,x))$, respectively.

\medskip

We remark that $\Lambda\!:[0,\widetilde{\zeta})\to(r_0,\infty)$ is a probabilistic solution for any $0\le\widetilde{\zeta}<\zeta$ (instantaneous melting to the cold core at $\widetilde{\zeta}$ rather than at $\zeta$).~To rule out such arbitrary discontinuities we use the notion of a \textit{physical} solution to \eqref{intro: heat rs}--\eqref{intro: Stefan rs}:~A probabilistic solution to \eqref{intro: heat rs}--\eqref{intro: Stefan rs} is called physical if, for all $t\in[0,\zeta)$, 
\begin{equation}\label{intro: phys rig} 
\begin{split}
& (\Lambda_{t-}-\Lambda_t)_+ = \inf\bigg\{y\in(0,\Lambda_{t-}- r_0]:
\;\nu([\Lambda_{t-}-y,\Lambda_{t-})) \\
&\qquad\qquad\qquad\qquad\qquad\qquad\qquad\qquad\;\;
>\int_{r_0}^\infty x^2w(0-,x)\,\PP^x(R_t\in[\Lambda_{t-}-y,\Lambda_{t-}),
\,t\le\tau\wedge\tau_{r_0})\,\mathrm{d}x\bigg\}, \\ 
& (\Lambda_{t-}-\Lambda_t)_- = \inf\bigg\{y\in(0,\infty):
\;\nu((\Lambda_{t-},\Lambda_{t-}+y]) \\
&\qquad\qquad\qquad\qquad\qquad\qquad\quad\;\;\,
>\int_{r_0}^\infty -x^2w(0-,x)\,\PP^x(R_t\in(\Lambda_{t-},\Lambda_{t-}+y],\,t\le\tau\wedge\tau_{r_0})\,\mathrm{d}x
\bigg\}; \\ 
& \nu([\Lambda_{\zeta-}-y,\Lambda_{\zeta-}))
\le \int_{r_0}^\infty x^2w(0-,x)\,\PP^x(R_\zeta\in[\Lambda_{\zeta-}-y,\Lambda_{\zeta-}),\,\zeta\le\tau\wedge\tau_{r_0})\,\mathrm{d}x,\;\; y\in(0,\Lambda_{\zeta-}-r_0)\\
&\qquad\qquad\qquad\qquad\qquad\qquad\qquad\qquad\qquad\qquad\qquad\qquad\qquad\qquad\qquad\qquad\qquad\qquad\quad
\text{if}\;\;\zeta<\infty. 
\end{split}
\end{equation}
The physical interpretation of \eqref{intro: phys rig} becomes apparent once it is rewritten with the help of \eqref{intro: def w}:
\begin{equation}\label{intro: phys} 
\begin{split}
& (\Lambda_{t-}-\Lambda_t)_+ = \inf\bigg\{y\in(0,\Lambda_{t-}- r_0]:
\,\int_{\Lambda_{t-}-y}^{\Lambda_{t-}}\mathrm{d}\nu
>\int_{\Lambda_{t-}-y}^{\Lambda_{t-}} w(t-,\cdot)\,\mathrm{d}\nu\bigg\}, \\ 
& (\Lambda_{t-}-\Lambda_t)_- =  \inf\bigg\{y\in(0,\infty):
\,\int_{\Lambda_{t-}}^{\Lambda_{t-}+y} \mathrm{d}\nu
>\int_{\Lambda_{t-}}^{\Lambda_{t-}+y} -w(t-,\cdot)\,\mathrm{d}\nu\bigg\}; \\ 
& \int_{\Lambda_{\zeta-}-y}^{\Lambda_{\zeta-}} \mathrm{d}\nu
\le \int_{\Lambda_{\zeta-}-y}^{\Lambda_{\zeta-}} w(\zeta-,\cdot)\,\mathrm{d}\nu,
\;\; y\in(0,\Lambda_{\zeta-}-r_0)\quad\text{if}\;\;\zeta<\infty. 
\end{split}
\end{equation}
The inequality in the first line of \eqref{intro: phys} ensures that a unit volume of the solid can only melt if its temperature exceeds the associated Gibbs-Thomson value by one unit, and the inequalities in the second and third lines admit similar interpretations.~As mentioned above, this notion of a physical solution was proposed in \cite{NaShsurface}. It originates from the work of \textsc{Delarue}, \textsc{Inglis}, \textsc{Rubenthaler} and \textsc{Tanr\'{e}} \cite{DIRT2}, and was employed recently in \cite{dns} to establish the uniqueness in the one-phase supercooled Stefan problem in spatial dimension one. In addition, it is shown in \cite[Section 5]{nsz} that physical solutions lie between the classical and the weak ones.
Compared to the one-phase one-dimensional supercooled Stefan problem, the problem at hand has two new challenging features:~(i) the non-linear dependence of the left-hand side in \eqref{intro: growth weak} on $\Lambda$, and (ii) the presence of a second phase. As a consequence of (ii), the present problem does not enjoy the ``partial comparison principle" available in the one-phase supercooled Stefan problem, which forces us to develop a new method for deducing uniqueness from partial regularity. Section \ref{se:2} describes this method and constitutes the main mathematical contribution of this work.

\medskip

Our standing assumption on the initial data reads as follows. 

\begin{ass}\label{intro: main ass}
The initial data $(\Lambda_{0-}>r_0,\,w(0-,\cdot)\!:(r_0,\infty)\to\rr)$ satisfies:
\begin{enumerate}[(a)]
\item There exists a constant $C<\infty$ such that $0\le w(0-,x)\le C(x-\Lambda_{0-})x^{-2}$, $x\ge\Lambda_{0-}$.
\item On $(r_0,\Lambda_{0-}]$, the function $w(0-,\cdot)$ is non-negative, bounded, and changes monotonicity finitely often on compact intervals.
\end{enumerate}
\end{ass}

The subsequent theorem is the main result of this paper. 

\begin{theorem}\label{intro: main thm}
Under Assumption \ref{intro: main ass} the solution $(\Lambda,\zeta)$ to \eqref{intro: growth weak}, \eqref{intro: phys rig} is unique. 
\end{theorem}

It is worth commenting on the different parts of Assumption \ref{intro: main ass}.

\begin{remark}\label{intro:rmk}
\begin{enumerate}[(a)]
\item The requirement $w(0-,x)\ge 0$, $x>r_0$ and the growth condition \eqref{intro: growth weak} imply that $\Lambda$ is non-increasing, i.e., the ball of ice is melting.
\item The inequality $w(0-,x)\le C(x-\Lambda_{0-})x^{-2}$, $x\ge\Lambda_{0-}$ entails linear bounds $w(0-,x)\le\widetilde{C}(x-\Lambda_{0-})$ on compact subsets of $[\Lambda_{0-},\infty)$, as well as the asymptotics $w(0-,x)=O(x^{-1})$, $x\to\infty$.
\item The finite number of monotonicity changes in $w(0-,\cdot)$ on compact sub-intervals of $(r_0,\Lambda_{0-}]$ parallels the central assumption in \cite{dns} (see \cite[Theorem 1.1]{dns}).~We only impose it on the solid phase, as $w(0-,x)\ge 0$, $x\ge\Lambda_{0-}$ automatically ensures sufficient regularity of $w$ in the liquid phase. 
\end{enumerate}
\end{remark}

The proof of Theorem \ref{intro: main thm} builds on a regularity preservation result of independent interest.

\begin{theorem}\label{intro: reg thm}
Under Assumption \ref{intro: main ass}, the following are true for any solution $(\Lambda,\zeta)$ to \eqref{intro: growth weak}, \eqref{intro: phys rig} and any given $t_0\in[0,\zeta)$.
\begin{enumerate}[(a)]
\item Liquid phase:~There exists a constant $C_{t_0}<\infty$ such that $0\le w(t_0-,x)\le C_{t_0}(x-\Lambda_{t_0-})x^{-2}$, $x\ge\Lambda_{t_0-}$.
\item Solid phase:~On $(r_0,\Lambda_{t_0-}]$, the function $w(t_0-,\cdot)$ is non-negative, bounded, and changes monotonicity finitely often on compact intervals.
\item Boundary: there is an $\upsilon_0>0$ such that $\Lambda\in C([t_0,t_0+\upsilon_0))\cap C^1((t_0,t_0+\upsilon_0))$.
\end{enumerate}
\end{theorem}

\smallskip

The rest of the paper is devoted to the proofs of Theorems \ref{intro: main thm} and \ref{intro: reg thm}. Section \ref{se:1.5} establishes certain preliminary results on the ``forward'' and ``backward" representations of a physical solution, in particular proving some of the heuristic derivations presented above.~In Section \ref{se:2}, Theorem~\ref{intro: main thm} is deduced from Theorem \ref{intro: reg thm}.~Considering two solutions $\Lambda$, $\widetilde{\Lambda}$ on a right neighborhood of some given $t\in[0,\zeta\!\wedge\!\widetilde{\zeta})$ we split the argument into two main parts:~(i) a delicate estimate on the difference between the melting rate contributions of the solid phases of the two solutions, via a detailed analysis of appropriate hitting probabilities for a three-dimensional Brownian motion, and (ii) the proof that the difference between the melting rate contributions of the liquid phases of the two solutions is dominated by the one of the solid phases, thanks to a time-reversal of the growth condition and a novel convexity property of $\Lambda$, $\widetilde{\Lambda}$ (obtained by comparison techniques).~In Section~\ref{se:3}, we develop the regularity theory for physical solutions which includes all statements in Theorem~\ref{intro: reg thm} apart from the finite number of monotonicity changes on compact sub-intervals of $(r_0,\Lambda_{t_0-}]$.~The latter is shown in Appendix \ref{se:4}.~The general strategy in Section \ref{se:3} and Appendix \ref{se:4} parallels that in \cite[Sections 2--4]{dns}.~However, many proofs are carefully adapted due to the two new challenging features of the problem \eqref{intro: growth weak} mentioned above:~the non-linear dependence of the left-hand side in \eqref{intro: growth weak} on $\Lambda$, and the presence of a second phase.

\section{Preliminaries}\label{se:1.5}

We start our analysis by defining, for a physical solution $\Lambda$ and all $t\in[0,\zeta)$, $x>0$, 
\begin{equation} \label{eq.sec3.sergey.w.def}
w(t,x)=\EE^x[w(0-,R_t)\,\bone_{\{\tau_{\Lambda_{t-\cdot}}\wedge\tau_{r_0}>t\}}]\leq \EE^x[w(0-,R_t)\,\bone_{\{\tau_{\Lambda_{t-\cdot}}>t\}}],
\end{equation}
where we set $w(0-,x)=0$ for $x\in(0,r_0]$, 
\begin{align}
&\tau_{\Lambda_{t-\cdot}}= \inf\{s\in[0,t+1]\!:(R_s-\Lambda_{t-s})(R_0-\Lambda_t)\leq0\},\label{eq.sec3.sergey.tauLambda.tminus.def}
\end{align}
and $\Lambda_s=\Lambda_{0-}$ for $s\in[-1,0)$.
Note that for any $t\in[0,\zeta)$: $w(t,x)=0$, $x\in(0,r_0]$. We let
\begin{align}
&\tau^-_{\Lambda_{t-\cdot}}:= \inf\{s\in[0,t+1]\!:(R_s-\Lambda_{t-s})(R_0-\Lambda_{t-})\leq0\}, \\
&\tau_{\Lambda_{t+\cdot}}:= \inf\{s\geq 0\!:(R_s-\Lambda_{t+s})(R_0-\Lambda_t)\leq0\},\label{eq.sec3.sergey.tauLambda.tplus.def}\\
&\tau^-_{\Lambda_{t+\cdot}}:= \inf\{s\geq 0\!:(R_s-\Lambda_{t+s})(R_0-\Lambda_{t-})\leq0\},\label{eq.sec3.sergey.tauLambda.minus.tplus.def}
\end{align}
and remark in passing that the stopping time $\tau$ appearing in Section \ref{se:intro} is equal to $\tau^-_{\Lambda_{0+\cdot}}$.

\medskip

The following ``Markov property'' and ``forward representation'' of $w$ are used repeatedly. The latter, in particular, shows that $x^2w(t,x)$ is the density of $R_t$ on $\{\tau^-_{\Lambda_{0+\cdot}}\wedge\tau_{r_0}> t\}$.

\begin{proposition}\label{prop:MarkovSystem}
For any $0\leq t_0\leq t<\zeta$, $x>0$ and any Borel $A\subset(0,\infty)$, 
\begin{IEEEeqnarray*}{rCl}
&&w(t,x)=\EE^x[w(t_0,R_{t-t_0})\,\bone_{\{\tau_{\Lambda_{t-\cdot}}\wedge\tau_{r_0}>t-t_0\}}],\\
&&\int_{A} w(t,x)\,\nu(\mathrm{d}x)=\int_{r_0}^\infty w(0-,x)\,\PP^x\big(R_t\in A,\,\tau^-_{\Lambda_{0+\cdot}}\wedge\tau_{r_0}>t\big)\,\nu(\mathrm{d}x).
\end{IEEEeqnarray*}
\end{proposition}

\noindent\textbf{Proof.} For the first statement, we apply the Markov property of the Bessel process $R$: 
\begin{IEEEeqnarray*}{rCl}
\EE^x[w(0-,R_t)\,\bone_{\{\tau_{\Lambda_{t-\cdot}}\wedge\tau_{r_0}>t\}}\,|\,\cF_{t-t_0}]
&=&\EE^{R_{t-t_0}}[w(0-,R_{t_0})\,\bone_{\{\tau_{\Lambda_{t_0-\cdot}}\wedge\tau_{r_0}>t_0\}}]\,\bone_{\{\tau_{\Lambda_{t-\cdot}}\wedge\tau_{r_0}>t-t_0\}}\\
&=&w(t_0,R_{t-t_0})\,\bone_{\{\tau_{\Lambda_{t-\cdot}}\wedge\tau_{r_0}>t-t_0\}}.
\end{IEEEeqnarray*}
The tower property of conditional expectation then yields
\begin{IEEEeqnarray*}{rCl}
w(t,x)=\EE^x[w(0-,R_t)\,\bone_{\{\tau_{\Lambda_{t-\cdot}}\wedge\tau_{r_0}>t\}}]=\EE^x[w(t_0,R_{t-t_0})\,\bone_{\{\tau_{\Lambda_{t-\cdot}}\wedge\tau_{r_0}>t-t_0\}}].    
\end{IEEEeqnarray*}

\smallskip

For the second statement, consider the functional $F$ on $C([0,t])$ given by
\begin{IEEEeqnarray*}{rCl}
F(\omega)=f_1(\omega_0)\,f_2(\omega_t)\,\bone_{\{\inf_{0\le s\le t}\,(\omega_s-\Lambda_s)(\omega_0-\Lambda_{0-})>0,\,\inf_{0\le s\le t}\,(\omega_s-r_0)>0\}},    
\end{IEEEeqnarray*}
where $f_1(x):=w(0-,x)$, $f_2(x):=\bone_{A}(x)$. Let $\stackrel{\leftarrow}{F}$ be the time-reversed functional defined through
\begin{IEEEeqnarray*}{rCl}
\stackrel{\leftarrow}{F}(\omega)=F(\omega_{t-\cdot}),
\end{IEEEeqnarray*}
and let $\nu^*(\cdot):=\int_0^\infty \PP^x(R\in\cdot)\,\nu(\mathrm{d}x)$ on $C([0,t])$. Since $\nu$ is a reversible measure for $R$ (see, e.g., \cite[``time-reversal formula'' (11)]{Law}),
\begin{IEEEeqnarray*}{rCl}
\int_{C([0,t])} F(\omega)\,\nu^*(\mathrm{d}\omega)
=\int_{C([0,t])} \stackrel{\leftarrow}{F}(\omega)\,\nu^*(\mathrm{d}\omega).    
\end{IEEEeqnarray*}
Consequently,
\begin{equation*}
\begin{split}
\int_{r_0}^\infty w(0-,x)\,\PP^x\big(R_t\in A,\,\tau^-_{\Lambda_{0+\cdot}}\!\wedge\!\tau_{r_0}\!>\! t\big)\,\nu(\mathrm{d}x)
=\int_{r_0}^\infty \EE^x[F(R)]\,\nu(\mathrm{d}x) 
=\int_{C([0,t])} F(\omega)\mathop{\nu^*(\mathrm{d}\omega)} \\
=\int_{C([0,t])} \stackrel{\leftarrow}{F}(\omega)\,\nu^*(\mathrm{d}\omega)
=\int_{A} \EE^x[w(0-,R_{t})\,\bone_{\{\tau_{\Lambda_{t-\cdot}}\wedge\tau_{r_0}>t\}}]\,\nu(\mathrm{d}x)
=\int_{A} w(t,x)\,\nu(\mathrm{d}x).
\end{split}
\end{equation*}
This finishes the proof. \qed

\medskip

Recall the growth condition \eqref{intro: growth weak}. In view of Proposition \ref{prop:MarkovSystem} and the Markov property of the Bessel process killed upon hitting $r_0$ or $\Lambda$, we have for any $0\leq t_0<t<\zeta$:
\begin{equation}\label{eq.sec3.sergey.growth.cond.t0}
\frac{\Lambda_{t_0}^3}{3}-\frac{\Lambda_t^3}{3}=\int_{r_0}^{\infty}w(t_0,x)\,\PP^x(\tau_{\Lambda_{t_0+\cdot}}\leq (t-t_0)\wedge\tau_{r_0})\,\nu(\mathrm{d}x).    
\end{equation}
In addition, it is easy to see from \eqref{eq.sec3.sergey.w.def} (cf.~\cite[proof of Lemma 3.2]{NaShsurface}) that $w(t-,x):=\lim_{s\uparrow t} w(s,x)$ is well-defined for every 
$t\in[0,\zeta]$ and $x\in(0,\infty)\backslash\{\Lambda_t,\Lambda_{t-}\}$, that $w(t,x)=w(t-,x)$ for $x\notin[\Lambda_t,\Lambda_{t-}]$, and that $w(t,x)=0$ for $x\in(\Lambda_t,\Lambda_{t-})$. Passing to the limit in \eqref{eq.sec3.sergey.growth.cond.t0} we thus deduce for any $0\leq t_0\leq t<\zeta$:
\begin{equation}\label{eq.sec3.sergey.growth.cond.t0.minus}
\frac{\Lambda_{t_0-}^3}{3}-\frac{\Lambda_t^3}{3}=\int_{r_0}^{\infty}w(t_0-,x)\,\PP^x(\tau^-_{\Lambda_{t_0+\cdot}}\leq (t-t_0)\wedge\tau_{r_0})\,\nu(\mathrm{d}x).    
\end{equation}
Moreover, we find that for all $0\leq t_0\leq t<\zeta$ and $x>0$, $\PP^x$-almost surely, it holds $w(t_0-,R_{t-t_0})=w(t_0,R_{t-t_0})$ on $\{\tau_{\Lambda_{t-\cdot}}\wedge\tau_{r_0}>t-t_0\}$. Hence, with Proposition \ref{prop:MarkovSystem},
\begin{equation}\label{eq.section3.sergey.w.back.Markov.tminus}
w(t,x)=\EE^x[w(t_0-,R_{t-t_0})\,\bone_{\{\tau_{\Lambda_{t-\cdot}}\wedge\tau_{r_0}>t-t_0\}}].
\end{equation}
Lastly, taking limits in Proposition \ref{prop:MarkovSystem} we identify $x^2w(t-,x)$ as the density of $R_t$ on $\{\tau^-_{\Lambda_{0+\cdot}}\!\wedge\!\tau_{r_0}\geq t\}$:
\begin{align*}
&\int_A w(t-,x)\,\nu(\mathrm{d}x)=\int_{r_0}^\infty w(0-,x)\,
\PP^x(R_t\in A,\,\tau^-_{\Lambda_{0+\cdot}}\wedge\tau_{r_0}\geq t)\,\nu(\mathrm{d}x),\quad A\in\mathcal{B}_{(0,\infty)},\quad t\in[0,\zeta].
\end{align*}
In particular, this makes the equivalence of the physicality conditions \eqref{intro: phys rig}  and \eqref{intro: phys} rigorous.


\section{Proof of Theorem \ref{intro: main thm}} \label{se:2}

\subsection{Main line of the argument}


In this section, we prove Theorem \ref{intro: main thm} assuming Theorem \ref{intro: reg thm}. To this end, we let $\Lambda$, $\widetilde{\Lambda}$ be two physical solutions, with the associated $(\zeta,w)$ and $(\widetilde\zeta,\widetilde w)$, and set
\begin{equation}
t_0:=\inf\{t\in[0,\zeta\wedge\widetilde\zeta):\,\Lambda_t\neq\widetilde{\Lambda}_t\}. 
\end{equation}
Assuming $t_0<\zeta\wedge\widetilde\zeta$ (otherwise, $t_0=\infty$ and, hence, $\zeta=\widetilde\zeta$ and $\Lambda$ coincides with $\widetilde{\Lambda}$) we derive a contradiction. 
Since $\Lambda$ and $\widetilde\Lambda$ can be interchanged, it is no loss of generality under the latter assumption to conclude that, for any $T\in(t_0,\zeta\wedge\widetilde\zeta)$, there exists $t\in(t_0,T]$ at which
\begin{equation}\label{t_def}
\Lambda^3_t-\widetilde{\Lambda}^3_t = \max_{t_0\le s\le t}\,|\Lambda^3_s-\widetilde{\Lambda}^3_s|>0.
\end{equation}
Thus, assuming that, for any $T\in(t_0,\zeta\wedge\widetilde\zeta)$, a time $t\in(t_0,T]$ satisfying \eqref{t_def} exists, we aim to arrive at a contradiction. (The latter suffices to prove Theorem \ref{intro: main thm}.)

\smallskip

To this end, we note that $\Lambda_{t_0}=\widetilde{\Lambda}_{t_0}>r_0$ and that $\Lambda$, $\widetilde{\Lambda}$ are continuous and non-increasing on $[t_0,T_0]$ for some $T_0\in(t_0,\zeta\wedge\widetilde\zeta)$ (see Theorem \ref{intro: reg thm}(c) and Remark \ref{intro:rmk}(a)).
We also notice that, since $\Lambda$ and $\widetilde{\Lambda}$ coincide on $[0,t_0)$, and in view of the jump condition \eqref{intro: phys} (recall that only downward jumps are possible), we have $\Lambda_{t_0}=\widetilde{\Lambda}_{t_0}$. Then, we also have $w(t_0,\cdot)\equiv\widetilde{w}(t_0,\cdot)$, as follows, for example, from \eqref{eq.sec3.sergey.w.def}. We write $f$ for $w(t_0,\cdot)\equiv\widetilde{w}(t_0,\cdot)$ and use \eqref{eq.sec3.sergey.growth.cond.t0} to estimate, for any $t\in(t_0,\zeta\wedge\widetilde\zeta)$:
\begin{equation}\label{uniq preamble}
\begin{split}
\frac{\Lambda^3_t}{3}-\frac{\widetilde{\Lambda}_t^3}{3}
& = \int_{r_0}^\infty f(x)\,\PP^x\big(\widetilde{\tau}\le(t-t_0)\wedge\tau_{r_0}\big)\,\nu(\mathrm{d}x)
-\int_{r_0}^\infty f(x)\,\PP^x\big(\tau\le(t-t_0)\wedge\tau_{r_0}\big)\,\nu(\mathrm{d}x) \\
& \le \int_{r_0}^{\Lambda_{t_0}} f(x)\,\PP^x\big(\widetilde{\tau}\le (t-t_0)\wedge\tau_{r_0},\,\tau>(t-t_0)\wedge\tau_{r_0}\big)\,\nu(\mathrm{d}x)\\
&\;\;\;\, +\int_{\Lambda_{t_0}}^\infty f(x)\,\PP^x\big(\widetilde{\tau}\le (t-t_0)\wedge\tau_{r_0},\,\tau>(t-t_0)\wedge\tau_{r_0}\big)\,\nu(\mathrm{d}x),
\end{split}
\end{equation}
where we simplify the notation by denoting $\tau:=\tau_{\Lambda_{t_0+\cdot}}$ and $\widetilde\tau:=\tau_{\widetilde\Lambda_{t_0+\cdot}}$ (recall \eqref{eq.sec3.sergey.tauLambda.tplus.def}).

\smallskip

To analyze the first summand in the last expression, for any small enough $T\in(t_0,T_0]$ and $0<\varepsilon_2\le(\Lambda_T\!\wedge\!\widetilde{\Lambda}_T-r_0)/2$, we consider the function $\psi\!:(0,\Lambda_{t_0}-\Lambda_T\wedge\widetilde{\Lambda}_T+\varepsilon_2]\to(0,1]$ satisfying
\begin{equation}\label{psi bound}
f(\Lambda_{t_0}-x) =
 1-\psi(x),\quad x\in(0,\Lambda_{t_0}-\Lambda_T\wedge\widetilde{\Lambda}_T+\varepsilon_2]
\end{equation}
and such that either $\psi$ is non-decreasing and $\psi(0+)=0$, or $\psi(0+)>0$. Here, we have used the local monotonicity of $f$ (see Theorem \ref{intro: reg thm}(b)) and the jump condition \eqref{intro: phys}. 
The desired contradiction with \eqref{t_def} is, then, a direct consequence of \eqref{uniq preamble}, Theorem \ref{intro: reg thm} and the following two propositions.

\begin{proposition}\label{prop multi}
Let $\Lambda,\,\widetilde\Lambda\!:[t_0,T_0]\rightarrow(r_0,\infty)$ be continuous non-increasing functions and let $f\!:(r_0,\infty)\rightarrow\RR$ be such that Theorem \ref{intro: reg thm}(a),(b) hold with $f$ in place of $w(t_0-,\cdot)$ (which yields \eqref{psi bound}).
Then, for any fixed small enough $\varepsilon_2>0$, there is a $T\in(t_0,T_0]$ such that, for all $t\in(t_0,T]$,
\begin{align*}
&\int_{r_0}^{\Lambda_{t_0}} f(x)\,\PP^x\big(\widetilde{\tau}\le (t-t_0)\wedge\tau_{r_0},\,\tau>(t-t_0)\wedge\tau_{r_0}\big)\,\nu(\mathrm{d}x)\\
&\le \frac{1}{3}\max_{t_0\leq s\leq t}(\Lambda^3_s-\widetilde{\Lambda}_s^3)_+
-\frac{1}{36}\max_{t_0\leq s\leq t} (\Lambda^3_s-\widetilde{\Lambda}_s^3)_+\,
\EE\big[\psi\big(\max_{0\le s\le t-t_0} (B_s+\Lambda_{t_0}-\Lambda_{t_0+s})\wedge\varepsilon_2\big)\big],
\end{align*}
where $B$ is a standard Brownian motion in $\rr$.
\end{proposition}

\begin{proposition}\label{prop 2nd phase}
Let the assumptions of Proposition \ref{prop multi} hold and assume that either $\psi(0+)>0$ or $\Lambda,\,\widetilde\Lambda$ satisfy \eqref{eq.sec3.sergey.growth.cond.t0} for $t\in[t_0,T_0]$ with $f$ in place of $w(t_0,\cdot)$. Then, for any fixed small enough $\varepsilon_2>0$, there is a $T\in(t_0,T_0]$ such that, for all $t\in(t_0,T]$,
\begin{equation}
\begin{split}
&\,\int_{\Lambda_{t_0}}^\infty f(x)\,\PP^x\big(\widetilde{\tau}\le (t-t_0)\wedge\tau_{r_0},\,\tau>(t-t_0)\wedge\tau_{r_0}\big)\,\nu(\mathrm{d}x)\label{eq psi 2nd phase}\\
&\leq\frac{1}{72}\max_{t_0\leq s\leq t} |\Lambda^3_s-\widetilde{\Lambda}_s^3|\,
\EE\big[\psi\big(\max_{0\le s\le t-t_0} (B_s+\Lambda_{t_0}-\Lambda_{t_0+s})\wedge\varepsilon_2\big)\big],
\end{split}
\end{equation} 
where $B$ is a standard Brownian motion in $\rr$.
\end{proposition}

\begin{remark}\label{rmk dep}
The proofs of Propositions \ref{prop multi}, \ref{prop 2nd phase} shed light on the dependence of $T$ on $\Lambda,\,\widetilde\Lambda$. Namely, assume that $f:(r_0,\infty)\rightarrow\RR$ is such that Theorem \ref{intro: reg thm}(a),(b) hold with $f$ in place of $w(t_0-,\cdot)$, and that $\psi(0+)>0$.
Then, for any small enough $\varepsilon_2>0$, there exist $\varepsilon_3>0$ and $\upsilon>0$, depending only on $(r_0,\Lambda_{t_0},f)$, such that the estimates of Propositions \ref{prop multi}, \ref{prop 2nd phase} hold for $T=t_0+\upsilon$ and for all continuous non-increasing $\Lambda$, $\widetilde\Lambda$ satisfying $\Lambda_{t_0}-\Lambda_s,\Lambda_{t_0}-\widetilde\Lambda_s\leq \varepsilon_3$ for all $s\in[t_0,t_0+\upsilon]$. This observation is used in Section \ref{se:3}. (Note that neither the statements nor the proofs of Propositions \ref{prop multi}, \ref{prop 2nd phase} rely on the results of Section \ref{se:3}.)
\end{remark}

\subsection{Proof of Proposition \ref{prop multi}}


To ease the notation, we take $t_0=0$ throughout the proof. Denote by $B_r(z)$ the ball of radius $r>0$ centered at $z\in\RR^3$. Then, for any $t\in(t_0,T_0]$ and any $0<\varepsilon_1<\varepsilon_2\le(\Lambda_{T_0}\!\wedge\!\widetilde{\Lambda}_{T_0}-r_0)/2$, we have
\begin{align}
&\qquad\,\int_{r_0}^{\Lambda_{0}} f(x)\,\PP^x\big(\widetilde{\tau}\leq t\wedge\tau_{r_0},\,\tau>t\wedge\tau_{r_0}\big)\,\nu(\mathrm{d}x) \nonumber \\
& \qquad\leq C_1\EE\bigg[\bone_{A_{\varepsilon_2}}\,\int_{B_{\Lambda_{0}}(0)\setminus B_{r_{0}}(0)} f(|z|)\,\bone_{\{\max_{0\le s\le t} (|z+W_s|-\widetilde{\Lambda}_s)\geq0,
\,\max_{0\le s\le t} (|z+W_s|-\Lambda_s)<0\}}\,\mathrm{d}z\bigg] \nonumber \\
&\qquad\quad+C_1\EE\bigg[\bone_{A^c_{\varepsilon_2}}\,\int_{B_{\Lambda_{0}}(0)\setminus B_{r_{0}}(0)} f(|z|)\,\bone_{\{\max_{0\le s\le t\wedge\tau_{r_0}} (|z+W_s|-\widetilde{\Lambda}_s)\geq0,
\,\max_{0\le s\le t\wedge\tau_{r_0}} (|z+W_s|-\Lambda_s)<0\}}\,\mathrm{d}z\bigg] \nonumber \\
&\label{eq.first} \qquad \leq C_1 \EE \left[\bone_{A_{\varepsilon_1}}\,\int_{S^{2}} \int_{r_*}^{r_{**}} f(r)\,r^{2}\,\mathrm{d}r\,\sigma(\mathrm{d}v)\right]
+ C_1 \EE \left[\bone_{A_{\varepsilon_2}\setminus A_{\varepsilon_1}}\,\int_{S^{2}} \int_{r_*}^{r_{**}} f(r)\,r^{2}\,\mathrm{d}r\,\sigma(\mathrm{d}v)\right] \\
&\qquad\quad+C_1\EE\bigg[\bone_{A^c_{\varepsilon_2}}\,\int_{B_{\Lambda_{0}}(0)\setminus B_{r_{0}}(0)} f(|z|)\,\bone_{\{\max_{0\le s\le t\wedge\tau_{r_0}} (|z+W_s|-\widetilde{\Lambda}_s)\geq0,\,\max_{0\le s\le t\wedge\tau_{r_0}} (|z+W_s|-\Lambda_s)<0\}}\,\mathrm{d}z\bigg], \nonumber
\end{align}
where $W$ is a standard Brownian motion in $\rr^3$, both inequalities rely on $A_\varepsilon:=\{\max_{0\le s\le t}|W_s|\leq\varepsilon\}$, we have let $S^{2}$ be the unit sphere of dimension $2$, with area $1/C_1$ and surface measure $\sigma(\mathrm{d}v)$ (each element $v\in S^{2}$ is identified with a unit vector in $\RR^3$), and
\begin{align*}
&r_*:=r_*(W,v):=\inf\big\{r>r_0:\,\max_{0\le s\le t} (|rv+W_s|-\widetilde{\Lambda}_s)\geq0\big\}, \\
&r_{**}:=r_{**}(W,v):=\inf\big\{r\ge r_*:\,\max_{0\le s\le t} (|rv+W_s|-\Lambda_s)\geq0\big\}, \\
& s^*:=s^*(W,v):= \inf\big\{s>0:\,\max_{0\le s'\le s} (|r_*v+W_{s'}|-\widetilde{\Lambda}_{s'})\geq0\big\}, \\
&r^*:=r^*(W,v):=\inf\big\{r\geq r_*:\,|rv+W_{s^*}|-\Lambda_{s^*}\geq0\big\}\wedge \Lambda_{0}.
\end{align*}

\smallskip

\noindent{\bf Step 1}. In this step, we estimate the first term on the right-hand side of \eqref{eq.first}.

\begin{lemma}\label{le:rstar.main}
There exist $T,\varepsilon_0,\overline{C}_1,\overline{C}_2,\overline{C}_3>0$ such that, for all $t\in(0,T]$ and $\varepsilon\in(0,\varepsilon_0]$, almost surely on $A_{\varepsilon}$, for each $v\in S^{2}$ at least one of the following two statements holds:
\begin{enumerate}[{$\bullet$}]
\item $(r^*)^3 - (r_*)^3 \leq (\Lambda^3_{s^*} - \widetilde{\Lambda}^3_{s^*})_+$,

\item $(r^*)^3 - (r_*)^3 \leq (\Lambda^3_{s^*} - \widetilde{\Lambda}^3_{s^*})_+ - \frac{\widetilde\Lambda^2_{s^*}}{8}\,(\Lambda_{s^*} - \widetilde{\Lambda}_{s^*})_+\,\alpha^2 + 6 \widetilde{\Lambda}_{s^*}(\Lambda_{s^*} - \widetilde{\Lambda}_{s^*})_+\,\epsilon$ and
\begin{align*}
(\widehat{r}^*)^3 - (\widehat{r}_*)^3 
\leq (\Lambda^3_{s^*} - \widetilde{\Lambda}^3_{s^*})_+\vee	(\Lambda^3_{\widehat{s}^*} - \widetilde{\Lambda}^3_{\widehat{s}^*})_+  
- 6 \widetilde{\Lambda}_{s^*}(\Lambda_{s^*} - \widetilde{\Lambda}_{s^*})_+\,\epsilon + \frac{\widetilde{\Lambda}^2_{s^*}}{16}(\Lambda_{s^*} - \widetilde{\Lambda}_{s^*})_+\,\alpha^2,
\end{align*}
where $\widehat{r}_*,\widehat{s}^*,\widehat{r}^*$ correspond to $\widehat{v}:=-v$, and $\alpha,\widehat{\alpha},\epsilon$ are non-negative random variables depending on $v$ and satisfying
\begin{align*}
& \alpha \leq \overline{C}_1 \max_{q\in S^2:\,q\cdot v=0}\;\max_{0\le s\le t}\;(q\cdot W_s)\leq \overline {C}_1\varepsilon,
\quad\widehat{\alpha} \leq \overline{C}_2 \varepsilon,
\quad \epsilon \leq \overline{C}_3\, \alpha^2.
\end{align*}
\end{enumerate}
\end{lemma}

\smallskip

Before proving Lemma \ref{le:rstar.main} we record the following two corollaries of it. 

\begin{corollary}\label{cor:2.4}
For all $t,\varepsilon>0$ small enough, almost surely on $A_\varepsilon$, for all $v\in S^2$, we have
\begin{align*}
& (r^*)^3 - (r_*)^3 + (\widehat{r}^*)^3 - (\widehat{r}_*)^3
\leq 2\big((\Lambda^3_{s^*} - \widetilde\Lambda^3_{s^*})_+ \vee (\Lambda^3_{\widehat s^*} - \widetilde\Lambda^3_{\widehat s^*})_+\big).
\end{align*}
\end{corollary}

\smallskip

\noindent Next, we introduce the pairwise disjoint (random) sets
\begin{align*}
\Sigma_1&:=\big\{v\in S^2:\,(r^*)^3 - (r_*)^3\leq(\Lambda^3_{s^*} - \widetilde\Lambda^3_{s^*})_+ 
\vee (\Lambda^3_{\widehat s^*} - \widetilde\Lambda^3_{\widehat s^*})_+, \\ 
&\qquad\qquad\quad\;\;\; 
(\widehat{r}^*)^3 - (\widehat{r}_*)^3
\leq(\Lambda^3_{s^*} - \widetilde\Lambda^3_{s^*})_+ \vee (\Lambda^3_{\widehat s^*} - \widetilde\Lambda^3_{\widehat s^*})_+\big\}, \\
\Sigma_2&:=\big\{v\in S^2:\, (r^*)^3 - (r_*)^3>(\Lambda^3_{s^*} - \widetilde\Lambda^3_{s^*})_+ \vee (\Lambda^3_{\widehat s^*} - \widetilde\Lambda^3_{\widehat s^*})_+, \\
&\qquad\qquad\quad\;\;\; 
(\widehat{r}^*)^3 - (\widehat{r}_*)^3
\leq (\Lambda^3_{s^*} - \widetilde\Lambda^3_{s^*})_+ \vee (\Lambda^3_{\widehat s^*} - \widetilde\Lambda^3_{\widehat s^*})_+\big\}, \\
\Sigma_3&:=\big\{v\in S^2:\, (r^*)^3 - (r_*)^3\leq(\Lambda^3_{s^*} - \widetilde\Lambda^3_{s^*})_+ \vee (\Lambda^3_{\widehat s^*} - \widetilde\Lambda^3_{\widehat s^*})_+, \\
&\qquad\qquad\quad\;\;\; 
(\widehat{r}^*)^3 - (\widehat{r}_*)^3>(\Lambda^3_{s^*} - \widetilde\Lambda^3_{s^*})_+ \vee (\Lambda^3_{\widehat s^*} - \widetilde\Lambda^3_{\widehat s^*})_+\big\},
\end{align*}
and notice that
\begin{align*}
-\Sigma_1=\Sigma_1,\quad-\Sigma_2=\Sigma_3,\quad\Sigma_1\cup\Sigma_2\cup\Sigma_3=S^2,  
\end{align*}
where the last identity follows from Corollary \ref{cor:2.4}. We also define the corresponding weights
\begin{equation*}
\begin{split}
&\theta_1(v)=(\Lambda^3_{s^*} - \widetilde\Lambda^3_{s^*})_+ 
	\vee (\Lambda^3_{\widehat s^*} - \widetilde\Lambda^3_{\widehat s^*})_+, \\
&\theta_2(v)=(\Lambda^3_{s^*} - \widetilde\Lambda^3_{s^*})_+ - \frac{\widetilde\Lambda^2_{s^*}}{8}\,(\Lambda_{s^*} - \widetilde\Lambda_{s^*})_+\,\alpha^2 + 6 \widetilde\Lambda_{s^*}\,(\Lambda_{s^*} - \widetilde\Lambda_{s^*})_+\,\epsilon,\\
&\theta_3(v)=(\Lambda^3_{\widehat s^*} - \widetilde\Lambda^3_{\widehat s^*})_+ \vee (\Lambda^3_{s^*} - \widetilde\Lambda^3_{s^*})_+ 
- 6 \widetilde\Lambda_{\widehat s^*}\,(\Lambda_{\widehat s^*} - \widetilde\Lambda_{\widehat s^*})_+\,\widehat\epsilon + \frac{\widetilde\Lambda^2_{\widehat s^*}}{16}(\Lambda_{\widehat s^*} - \widetilde\Lambda_{\widehat s^*})_+\,\widehat\alpha^2.
\end{split}
\end{equation*}
Further, we define the (random) sets
\begin{align*}
& \Sigma_0=\bigg\{v\in S^2:\,
(\Lambda_{s^*}^3-\widetilde\Lambda_{s^*}^3)_+\leq\frac12\max_{0\le s\le t} (\Lambda_s^3-\widetilde\Lambda_s^3)_+\bigg\}, \\
& \Sigma'_0=\bigg\{v\in S^2:\,
(\Lambda_{s^*}^3-\widetilde\Lambda_{s^*}^3)_+\vee (\Lambda^3_{\widehat s^*} - \widetilde\Lambda^3_{\widehat s^*})_+\leq\frac12\max_{0\le s\le t} (\Lambda_s^3-\widetilde\Lambda_s^3)_+\bigg\}
\end{align*}
and the corresponding modified weights
\begin{align*}
& \theta_{0,i}(v)=\theta_i(v)\,\bone_{S^2\setminus\Sigma'_0}(v)
+\frac{2}{3}\max_{0\le s\le t} (\Lambda_s^3-\widetilde\Lambda_s^3)_+\,\bone_{\Sigma'_0}(v),\quad i=1,3, \\
& \theta_{0,2}(v)=\theta_2(v)\,\bone_{S^2\setminus\Sigma_0}(v)+\frac{2}{3}\max_{0\le s\le t} (\Lambda_s^3-\widetilde\Lambda_s^3)_+\,\bone_{\Sigma_0}(v).
\end{align*}

\begin{corollary}\label{cor2.5}
For all $t,\varepsilon>0$ small enough, almost surely on $A_\varepsilon$, for all $v\in S^2$, it holds: 
\begin{enumerate}[{$\bullet$}]
\item $(r^*)^3-(r_*)^3\leq \theta_{0,i}(v)$, $v\in\Sigma_i$, $i=1,2,3$;
\item $\theta_{0,1}(v)\leq\max_{0\le s\le t} (\Lambda_s^3-\widetilde\Lambda_s^3)_+$, $v\in\Sigma_1$;
\item $\theta_{0,2}(v)+\theta_{0,3}(-v)\leq 2\max_{0\le s\le t} (\Lambda_s^3-\widetilde\Lambda_s^3)_+$, $v\in\Sigma_2$;
\item $\theta_{0,i}(v)\ge\frac13\max_{0\le s\le t} (\Lambda_s^3-\widetilde\Lambda_s^3)_+$, $v\in\Sigma_i$, $i=1,2,3$.
\end{enumerate}
\end{corollary}


\smallskip

\noindent\textbf{Proof of Lemma \ref{le:rstar.main}.} Note that we always have $r_*<\Lambda_0$, and on $A_\varepsilon$, we have $r_*>0$. Fix an arbitrary $v\in S^2$, write $z^*$ for $r_*v + W_{s^*}$, which satisfies $|z^*|=\widetilde\Lambda_{s^*}$, and set $\underline{r}=z^*\cdot v$. Let us also observe that 
\begin{equation*}
r^* \leq r_* + h =: \overline{r},
\end{equation*}
where $h:=\sup\{r\geq 0\!:|z^*+rv| \leq \widetilde\Lambda_{s^*}\vee\Lambda_{s^*}\}$. We distinguish between three cases.


\medskip

\noindent\textbf{Case 1: $r_*=\underline{r}$.} Define $\alpha\in[0,\pi/2]$ as the angle between $v$ and $z^*$, and consider 
\begin{equation*}
R:= \widetilde\Lambda_{s^*}=|z^*|,\quad \Delta:=\widetilde\Lambda_{s^*}\vee\Lambda_{s^*} - R=(\Lambda_{s^*}-\widetilde\Lambda_{s^*})_+.
\end{equation*}
Then, by simple trigonometry we find that
\begin{equation*}
\underline{r} = R \cos \alpha,\quad \overline r = \left( (R+\Delta)^2 - R^2 \sin^2\alpha \right)^{1/2},
\quad h = \left( (R+\Delta)^2 - R^2 \sin^2\alpha \right)^{1/2} - R \cos \alpha.
\end{equation*}
In particular, 
\begin{equation}\label{size of h}
h = \Delta+O(\Delta\alpha^2),\quad\text{where}\;\;\alpha\leq \overline C_1\max_{q\in S^2:\,q\cdot v=0}\;\max_{0\le s\le t}\;(q\cdot W_s).
\end{equation}
The next lemma completes the proof of Lemma \ref{le:rstar.main} in the first case.

\begin{lemma}\label{lem2.6}
For all $\delta\in(0,1)$, there exist $\alpha_0,\Delta_0,C>0$ such that
\begin{equation*}
\begin{split}
\left( (R+\Delta)^2 - R^2 \sin^2\alpha \right)^{3/2} - R^3\cos^3\alpha
\leq (R+\Delta)^3 - R^3 - \frac{R^2}{4}\Delta\alpha^2 + C\Delta\alpha^4,\\
\alpha\in[0,\alpha_0),\;\Delta\in[0,\Delta_0),\;R\in[\delta,1/\delta].
\end{split}
\end{equation*}
\end{lemma}

\noindent\textbf{Proof.} The desired inequality is equivalent to
\begin{align*}
& (1+\Delta/R)^3 \left( 1 - \sin^2\alpha/(1+\Delta/R)^2 \right)^{3/2} - \cos^3\alpha
\leq (1+\Delta/R)^3 - 1 - \frac{1}{4R}\Delta\alpha^2 + O(\Delta\alpha^4).
\end{align*}
Taylor expanding in $\Delta/R$ on the left-hand side we obtain
\begin{align*}
&\,(1+\Delta/R)^3 \left( 1 - \sin^2\alpha/(1+\Delta/R)^2 \right)^{3/2} - \cos^3\alpha\\
&= (1+\Delta/R)^3 \left( 1 - \sin^2\alpha + 2\sin^2\alpha\,(\Delta/R) + \sin^2\alpha\,O(\Delta^2/R^2) \right)^{3/2} - \cos^3\alpha\\
& = (1+\Delta/R)^3 \big(\cos^2\alpha + \sin^2\alpha\,(\Delta/R)\big(2 + O(\Delta/R)\big)\big)^{3/2} - \cos^3\alpha\\
& = \cos^3\alpha\,\Big((1+\Delta/R)^3 \left(1 +\tan^2\alpha\,(\Delta/R) 
\big(2 + O(\Delta/R)\big)\right)^{3/2} - 1\Big).
\end{align*}
By another Taylor expansion we express the latter as
\begin{align*}
& \,\cos^3\alpha\,\bigg(\!(1+\Delta/R)^3 \bigg(1 + \frac{3}{2}\tan^2\alpha\,(\Delta/R)
\big(2 + O(\Delta/R)\big)
+ \tan^4\alpha\,O(\Delta^2/R^2)\!\bigg) - 1\bigg) \\
& = \cos^3\alpha\,\Big(\!(1+\Delta/R)^3 -1 + (1+\Delta/R)^3(\Delta/R)\tan^2\alpha\,
\big(3 + O(\Delta/R) + \tan^2\alpha\,O(\Delta/R)\big)\!\Big).
\end{align*}
Dividing the above by $\Delta$ and Taylor expanding the result in $\alpha$ we get
\begin{align*}
&\,\cos^3\alpha\,\bigg(\frac{(1+\Delta/R)^3 -1}{\Delta} + \frac{(1+\Delta/R)^3}{R} \tan^2\alpha 
\,\big(3+O(\Delta/R)\big)\!\bigg) \\
& = \bigg(1-\frac{3}{2}\alpha^2 + O(\alpha^4)\!\bigg)
\bigg(\frac{(1+\Delta/R)^3 -1}{\Delta} + \frac{(1+\Delta/R)^3}{R} (\alpha^2 + O(\alpha^4))
(3+O(\Delta/R))\!\bigg) \\
& = \frac{(1+\Delta/R)^3 -1}{\Delta} + \frac{(1+\Delta/R)^3}{R} \alpha^2 \left( 3+ O\left( \Delta/R\right)\right)
-\frac{3}{2} \alpha^2 \frac{(1+\Delta/R)^3 -1}{\Delta}
+ O(\alpha^4)\\
&\leq \frac{(1+\Delta/R)^3 - 1}{\Delta}
+ \frac{1}{R} \alpha^2 \frac{7}{2}
-\frac{3}{2R} \alpha^2 \frac{5}{2}
+ O(\alpha^4)
= \frac{(1+\Delta/R)^3 - 1}{\Delta} - \frac{1}{4R} \alpha^2 + O(\alpha^4). 
\end{align*}
The lemma readily follows. \qed


\smallskip

\noindent\textbf{Case 2: $r_*<\underline{r}$.} In this case, we only need to again apply Lemma \ref{lem2.6}, which yields
\begin{equation*}
(r^*)^3-(r_*)^3\le (r_*+h)^3 - (r_*)^3 
\le (\underline{r}+h)^3 - \underline{r}^3 \leq (\Lambda^3_{s^*} - \widetilde\Lambda^3_{s^*})_+
- \frac{R^2}{4}\Delta\alpha^2 + C\Delta\alpha^4.
\end{equation*}

\smallskip

\noindent\textbf{Case 3: $r_*>\underline{r}$.} By estimating the angles of the right triangles with the vertices $z^*$, $\underline{r}v$, $r_*v$, and $0$, $\underline{r}v$, $z^*$, we infer that
\begin{equation}\label{size of epsilon}
0<\epsilon:=r_* - \underline{r} \leq \overline C_3 \alpha^2.
\end{equation}
Now, consider $\widehat v\!=\!-v$ and the associated $\widehat r_*,\widehat s^*,\widehat z^*,\widehat h$. Then, $|(\underline{r}-\epsilon)\widehat{v}+W_{s^*}|\ge |z^*|=\widetilde{\Lambda}_{s^*}$ implies
\begin{equation*}
0<\widehat r_* \leq \underline{r} - \epsilon \leq r_*,
\end{equation*}
and thus, 
\begin{equation*}
(\widehat r_*+ \widehat h)^3 - (\widehat r_*)^3
= 3(\widehat{r}_*)^2 \widehat h + 3 \widehat{r}_* \widehat h^2 + \widehat h^3
\leq 3(\underline{r}^2 - 2 \underline{r} \epsilon + \epsilon^2) \widehat h + 3 (\underline{r} - \epsilon) \widehat h^2 + \widehat h^3.
\end{equation*}
Therefore, if $\widehat h\le h$, then the above, along with Lemma \ref{lem2.6}, \eqref{size of h}, and \eqref{size of epsilon}, lead to
\begin{align*}
(\widehat r_*+ \widehat h)^3 - (\widehat r_*)^3 
&\le 3(\underline{r}^2-2\underline{r}\epsilon+\epsilon^2)h+3(\underline{r}-\epsilon) h^2 + h^3,\\
&\le (\underline{r}+h)^3-\underline{r}^3-6\underline{r}\epsilon h+3\epsilon^2h \\
&\le (\Lambda^3_{s^*} - \widetilde\Lambda^3_{s^*})_+ - \frac{R^2}{8}\Delta\alpha^2-6\epsilon R\Delta(1+O(\alpha^2))+O(\Delta\alpha^4) \\
&\le (\Lambda^3_{s^*} - \widetilde\Lambda^3_{s^*})_+ - \frac{R^2}{8}\Delta\alpha^2-6\epsilon R\Delta+\frac{R^2}{16}\Delta\alpha^2.
\end{align*}

\smallskip

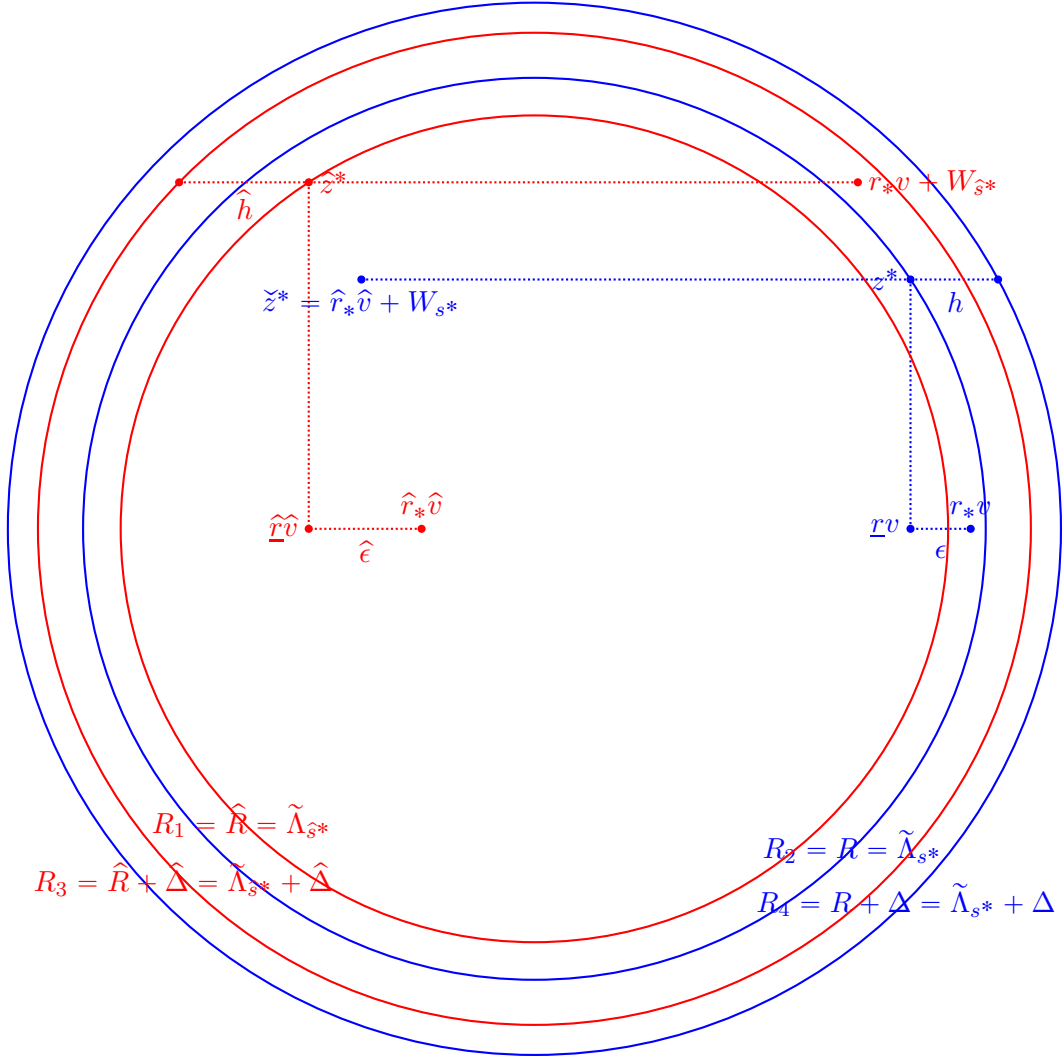
\begin{figure}[h]
	\begin{tikzpicture}
		\draw[blue, thick] (0,0) circle (7) node at (4.94, -4.94) {$R_4=R+\Delta=\widetilde\Lambda_{s^*}+\Delta$};
		\draw[blue, thick] (0,0) circle (6) node at (4.23,-4.23) {$R_2=R=\widetilde\Lambda_{s^*}$};
		\filldraw[blue, thick] (5,3.317) circle (0.04) node [anchor = east] {$z^*$};
		\filldraw[blue, thick] (5.8,0) circle (0.04) node[anchor = south]{$r_*v$};
		\filldraw[blue, thick] (5,0) circle (0.04) node [anchor = east]{$\underline{r}v$};
		\filldraw[blue, thick] (6.164,3.317) circle (0.04);
		\draw[densely dotted][blue, thick] (5,3.317) -- (5,0);
		\draw[densely dotted][blue, thick] (5.8,0) -- (5,0) node at (5.4,-0.3) {$\epsilon$};
		\draw[densely dotted][blue, thick] (5,3.317) -- (6.164,3.317) node at (5.6,3.017) {$h$};
		\draw[red, thick] (0,0) circle (6.6) node at (-4.667,-4.667) {$R_3=\widehat R+\widehat\Delta=\widetilde\Lambda_{\widehat s^*}+\widehat\Delta$};
		\draw[red, thick] (0,0) circle (5.5) node at (-3.889, -3.889) {$R_1=\widehat R=\widetilde\Lambda_{\widehat s^*}$};
		\filldraw[red, thick] (-3,4.610) circle (0.04) node [anchor = west] {$\widehat z^*$};
		\filldraw[red, thick] (-1.5,0) circle (0.04) node[anchor = south]{$\widehat r_*\widehat v$};
		\filldraw[red, thick] (-3,0) circle (0.04) node [anchor = east]{$\underline{\widehat r}\widehat v$};
		\filldraw[red, thick] (-4.723,4.610) circle (0.04);
		\draw[densely dotted][red, thick] (-3,4.610) -- (-3,0);
		\draw[densely dotted][red, thick] (-1.5,0) -- (-3,0) node at (-2.25,-0.3) {$\widehat\epsilon$};
		\draw[densely dotted][red, thick] (-4.723,4.610) -- (-3,4.610) node at (-3.85,4.310) {$\widehat h$};
		\filldraw[blue, thick] (-2.3,3.317) circle (0.04) node [anchor = north] {$\widecheck z^*=\widehat r_*\widehat v+W_{s^*}$};
		\draw[densely dotted][blue, thick] (-2.3,3.317) -- (5,3.317);
		\filldraw[red, thick] (4.3,4.610) circle (0.04) node [anchor = west] {$r_*v+W_{\widehat s^*}$};
		\draw[densely dotted][red, thick]
		(-3,4.610) -- (4.3,4.610);
	\end{tikzpicture}
\end{figure}
If $\widehat{h}>h$, we consider $\widehat{\underline r}$ and $\widehat{\epsilon}:= \widehat{\underline{r}}-\widehat{r}_*$. We claim that $\widehat\epsilon\geq\epsilon$. Indeed, assume without loss of generality that $v$ is the first canonical basis vector in $\rr^3$. Then, the first coordinate of $z^*$ is positive, while the first coordinate of $\widehat z^*$ is negative. Denote by $\widecheck{z}^*=\widehat r_* \widehat v + W_{s^*}$ the point obtained by shifting $z^*$ parallelly to the first axis by $-(r_*+\widehat r_*)$, and write $x$ for the first coordinate of $\widecheck z^*$. Then, $x\geq -\widehat{\underline{r}}$, since if we assume the opposite, $x< -\widehat{\underline{r}}$, then $\underline{\widehat{r}}<|x|$ and $\underline{\widehat{r}}<\underline{r}$. To see the latter, notice that $\widecheck z^*\in \overline{B_{\widetilde\Lambda_{s^*}}(0)}$ (otherwise we would encounter a contradiction to the definition of $\widehat r_*$), which implies $|x|\leq \underline{r}$ that can be combined with the assumed $\widehat{\underline{r}}<|x|$. Together $\underline{\widehat{r}}<|x|$, $\underline{\widehat{r}}<\underline{r}$ give
\begin{align*}
& 2\underline{\widehat{r}} < |x| + \underline{r} = r_*+\widehat{r}_*.
\end{align*}
At this point, we recall $|\widehat z^*|=\widetilde\Lambda_{\widehat s^*}$ and infer from $2\underline{\widehat{r}}<r_*+\widehat{r}_*$ that $\widehat z^*$ shifted parallelly to the first axis by $r_*+\widehat r_*$, yielding $r_*v+W_{\widehat{s}^*}$, lies outside of $\overline{B_{\widetilde\Lambda_{\widehat{s}^*}}(0)}$, a contradiction to the definition of~$r_*$.
Thus, $\underline{r} - r_* - \widehat r_*=x\ge -\widehat{\underline{r}}$, which gives the desired $\widehat{\epsilon}=\widehat{\underline{r}} - \widehat r_*\ge r_* - \underline{r}=\epsilon$.

\medskip


Furthermore, it is not hard to see that, on $A_\varepsilon$,
\begin{align*}
(\widetilde\Lambda_t^2-\varepsilon^2)^{1/2}
\leq\underline{r},\widehat{\underline{r}}\leq\widetilde\Lambda_0,
\quad\widetilde\Lambda_t-\varepsilon\leq r_*,\widehat r_*\leq\widetilde\Lambda_0.
\end{align*}
As a result, relying on the local continuity of $\Lambda$, $\widetilde{\Lambda}$, we conclude that, for any given $\varepsilon'>0$, the following holds on $A_\varepsilon$ for all small enough $t>0$ and $\varepsilon>0$: 
\begin{align*}
|\widehat{\underline{r}} - \underline{r}| \leq \varepsilon',
\quad \widehat{\underline{r}}/\underline{r},\,\underline{r}/\widehat{\underline{r}} \geq 1-\varepsilon',
\quad \widehat\epsilon \leq (1-\varepsilon')\,\underline{r}.
\end{align*}
Then, in view of $\epsilon\leq\widehat\epsilon$ and $h<\widehat{h}$, we obtain
\begin{align*}
&R(1-O(\alpha^2))\leq \underline{r}\leq R,
\quad \widehat R(1-O(\widehat\alpha^2))\leq \widehat{\underline{r}} \leq \widehat R,
\quad \widehat\epsilon\leq \widehat{\underline{r}},
\quad \Delta\leq h\leq \Delta+O(\Delta\alpha^2), \\
& (- 2\widehat{\underline{r}}\widehat\epsilon + \widehat\epsilon^2) \widehat h
\leq (- 2 \widehat{\underline{r}} \widehat\epsilon + \widehat\epsilon^2) h
\leq (- 2(1-\varepsilon') \underline{r} \widehat\epsilon + \widehat\epsilon^2) h
\leq (- 2(1-\varepsilon') \underline{r} \epsilon + \epsilon^2) h.
\end{align*}

\smallskip

Finally, we use the latter estimates, \eqref{size of epsilon}, and Lemma \ref{lem2.6} to deduce
\begin{align*}
(\widehat r_*+ \widehat h)^3 - (\widehat r_*)^3
= 3(\widehat{r}_*)^2 \widehat h + 3 \widehat{r}_* \widehat h^2 + \widehat h^3
= 3(\widehat{\underline{r}}^2 - 2 \widehat{\underline{r}} \widehat\epsilon + \widehat\epsilon^2) \widehat h + 3 (\widehat{\underline{r}} - \widehat\epsilon) \widehat h^2 + \widehat h^3 \\
\leq 3\widehat{\underline{r}}^2 \widehat h + 3 \widehat{\underline{r}} \widehat h^2 + \widehat h^3
- 6(1-\varepsilon')\, \underline{r}\epsilon h + O(\Delta \alpha^4) \\
\leq (\widehat{\underline{r}}+\widehat h)^3 - \widehat{\underline{r}}^3 
- 6(1-\varepsilon') R(1-O(\alpha^2))\epsilon\Delta + O(\Delta \alpha^4) \\
\leq (\Lambda^3_{\widehat s^*} - \widetilde\Lambda^3_{\widehat s^*})_+ 
- \frac{\widehat R^2}{8}\widehat\Delta\widehat\alpha^2 
- 6 R \Delta \epsilon + O(\varepsilon' \Delta \alpha^2) + O(\Delta \alpha^4), 
\end{align*}
\begin{align*}
(r_*+ h)^3 - (r_*)^3
= 3(r_*)^2 h + 3 r_* h^2 + h^3
= 3(\underline{r}^2 + 2 \underline{r} \epsilon + \epsilon^2) h + 3 (\underline{r} + \epsilon) h^2 + h^3 \\
\leq (\underline{r}+ h)^3 - \underline{r}^3 + 6 \underline{r}h\epsilon 
+ O(\Delta \alpha^4) + O(\Delta^2 \alpha^2) \\
\leq (\underline{r}+h)^3 - \underline{r}^3 + 6R\epsilon\Delta(1+O(\alpha^2))
 + O(\Delta \alpha^4) + O(\Delta^2 \alpha^2) \\
\leq (\Lambda^3_{s^*} - \widetilde\Lambda^3_{s^*})_+ 
- \frac{R^2}{4}\Delta\alpha^2 + 6 R\Delta\epsilon + O(\Delta \alpha^4) + O(\Delta^2 \alpha^2) \\
\leq (\Lambda^3_{s^*} - \widetilde\Lambda^3_{s^*})_+ - \frac{R^2}{8}\Delta\alpha^2 
+ 6 R\Delta\epsilon,
\end{align*}
for all small enough $t,\varepsilon>0$.
Choosing $\varepsilon'>0$ as a sufficiently small multiple of $R^2$ we get
\begin{align*}
& (\widehat{r}_*+ \widehat h)^3 - (\widehat r_*)^3 
\leq (\Lambda^3_{\widehat s^*} - \widetilde\Lambda^3_{\widehat s^*})_+ 
- \frac{\widehat R^2}{8}\widehat\Delta\widehat\alpha^2 
- 6 R \Delta \epsilon + \frac{R^2}{16} \Delta \alpha^2
\end{align*}
for all small enough $t,\varepsilon>0$. This finishes the proof of Lemma \ref{le:rstar.main}. \qed

\medskip

Returning to the proof of Proposition \ref{prop multi}, we use that $r_*\geq \Lambda_t\wedge\widetilde\Lambda_t-\varepsilon_2$ on $A_{\varepsilon_1}$, as well as $r+\max_{0\le s\le t} (v\cdot W_s-\Lambda_s)\leq0$, $r\in[r_*,r_{**}]$ (which implies $\Lambda_{0}-r\ge\max_{0\le s\le t} (v\cdot W_s-\Lambda_s+\Lambda_{0})$) and Corollary \ref{cor2.5}, to conclude that, for any fixed small enough $T\in(0,T_0]$ and $0<\varepsilon_2\le(\Lambda_T\!\wedge\!\widetilde{\Lambda}_T-r_0)/2$ (such that \eqref{psi bound} holds), the following holds for all small enough $t,\varepsilon_1>0$:
\begin{align}
& \; C_1\,\bone_{A_{\varepsilon_1}} \int_{S^{2}}\int_{r_*}^{r_{**}} f(r)\, r^{2}\,\mathrm{d}r\,\sigma(\mathrm{d}v)\nonumber \\
&=C_1\,\bone_{A_{\varepsilon_1}} \int_{S^2}\int_{r_*}^{r_{**}} (1-\psi(\Lambda_{0}-r))\,r^{2}\,\mathrm{d}r\,\sigma(\mathrm{d}v)\nonumber \\
&\le C_1\,\bone_{A_{\varepsilon_1}} \int_{S^2}\int_{r_*}^{r^*} \big(1-\psi\big(\max_{0\le s\le t} (v\cdot W_s-\Lambda_s+\Lambda_{0})\big)/2\big)\,r^{2}\,\mathrm{d}r\,\sigma(\mathrm{d}v)\nonumber \\
&= \sum_{i=1}^3 C_1\,\bone_{A_{\varepsilon_1}} \int_{\Sigma_i}\frac13((r^*)^3-(r_*)^3)\big(1-\psi\big(\max_{0\le s\le t} (v\cdot W_s-\Lambda_s+\Lambda_{0})\big)/2\big)\,\sigma(\mathrm{d}v)\nonumber \\
&\le \frac13\sum_{i=1}^3 C_1\,\bone_{A_{\varepsilon_1}} \int_{\Sigma_i} \theta_{0,i}(v)\big(1-\psi\big(\max_{0\le s\le t} (v\cdot W_s-\Lambda_s+\Lambda_{0})\big)/2\big)\,\sigma(\mathrm{d}v)\nonumber \\
&= \frac13\sum_{i=1}^3 C_1\,\bone_{A_{\varepsilon_1}} \int_{\Sigma_i} \theta_{0,i}(v)\,\sigma(\mathrm{d}v)
-\frac16\sum_{i=1}^3C_1\,\bone_{A_{\varepsilon_1}} \int_{\Sigma_i} \theta_{0,i}(v)\,\psi\big(\max_{0\le s\le t}(v\cdot W_s-\Lambda_s+\Lambda_{0})\big)\,\sigma(\mathrm{d}v)\nonumber \\
&\le \frac13C_1\bigg(\bone_{A_{\varepsilon_1}}\bigg(\sigma(\Sigma_1)\max_{0\le s\le t} (\Lambda_s^3-\widetilde\Lambda_s^3)_+ + \int_{\Sigma_2} \theta_{0,2}(v)\,\sigma(\mathrm{d}v)+\int_{\Sigma_2} \theta_{0,3}(-v)\,\sigma(\mathrm{d}v)\bigg)\!\bigg)\nonumber \\
&\quad\,-C_1\,\bone_{A_{\varepsilon_1}} \int_{S^2}\psi\big(\max_{0\le s\le t} (v\cdot W_s-\Lambda_s+\Lambda_{0})\big)\,\sigma(\mathrm{d}v)\,\frac{1}{18}\max_{0\le s\le t} (\Lambda_s^3-\widetilde\Lambda_s^3)_+ \label{eq.sec2.Prop21.step1.finalEq}\\
&\le \frac13C_1\,\bone_{A_{\varepsilon_1}}(\sigma(\Sigma_1)+2\sigma(\Sigma_2))\max_{0\le s\le t} (\Lambda_s^3-\widetilde\Lambda_s^3)_+\nonumber \\
&\quad\,-C_1\,\bone_{A_{\varepsilon_1}} \int_{S^2}\psi\big(\max_{0\le s\le t} (v\cdot W_s-\Lambda_s+\Lambda_{0})\big)\,\sigma(\mathrm{d}v)\,\frac{1}{18}\max_{0\le s\le t} (\Lambda_s^3-\widetilde\Lambda_s^3)_+\nonumber \\
&\leq \bone_{A_{\varepsilon_1}}\,\frac13\max_{0\le s\le t} (\Lambda_s^3-\widetilde\Lambda_s^3)_+
-C_1\,\bone_{A_{\varepsilon_1}}\int_{S^2}\psi\big(\max_{0\le s\le t} (v\cdot W_s-\Lambda_s+\Lambda_{0})\big)\,\sigma(\mathrm{d}v)\,\frac{1}{18}\max_{0\le s\le t} (\Lambda_s^3-\widetilde\Lambda_s^3)_+,\nonumber
\end{align}
where we recall that $1/C_1$ is the area of $S^{2}$. The above, together with the fact that the law of $(\bone_{A_{\varepsilon_1}}, v\cdot W)$ does not depend on $v$, yields the following bound for the first term on the right-hand side of \eqref{eq.first}: For any small enough $T\in(0,T_0]$ and $0<\varepsilon_2\le(\Lambda_T\!\wedge\!\widetilde{\Lambda}_T-r_0)/2$,
\begin{equation}\label{eq.sec2.prop2.1.step1.result}
\begin{split}
&\; C_1 \EE \left[\bone_{A_{\varepsilon_1}}\,\int_{S^{2}} \int_{r_*}^{r_{**}} f(r)\,r^{2}\,\mathrm{d}r\,\sigma(\mathrm{d}v)\right] \\
&\leq \PP(A_{\varepsilon_1})\,\frac13\max_{0\le s\le t} (\Lambda_s^3-\widetilde\Lambda_s^3)_+
-\EE\big[\bone_{A_{\varepsilon_1}}\,\psi\big(\max_{0\le s\le t} (W_s^{(1)}-\Lambda_s+\Lambda_{0})\big)\big]\,\frac{1}{18}
\max_{0\le s\le t} (\Lambda_s^3-\widetilde\Lambda_s^3)_+
\end{split}
\end{equation}
for all small enough $t,\varepsilon_1>0$, where $W^{(1)}$ is the first coordinate of $W$.

\medskip

\noindent\textbf{Step 2}. In this step, we estimate the second term on the right-hand side of \eqref{eq.first}.

\begin{lemma}\label{le:SupOrthogonalPerturbation}
For any $v\in S^2$, and for any $v'\in S^2$ forming an angle $\beta\in[0,\frac\pi2]$ with $v$, we have
\begin{align*}
& \max_{q\in S^2:\,q\cdot v'=0} (q\cdot x) 
\leq \max_{q\in S^{2}:\,q\cdot v=0} (q\cdot x + \sqrt{2}\beta|x|),\quad x\in\RR^3.
\end{align*}
\end{lemma}

\noindent\textbf{Proof.} It suffices to consider $\beta\in(0,\frac\pi2]$. Let $O$ be the orthogonal matrix which acts as a rotation by $\beta$ on the span of $v$, $v'$ and as the identity map on its orthogonal complement.~Then, $v'=Ov$ and $|O-I|=2\sqrt{2}\sin(\beta/2)$. Moreover,
\begin{align*}
\max_{q\in S^{2}:\,q\cdot v'=0} (q\cdot x)
=\max_{q\in S^{2}:\,q\cdot(Ov)=0} (q\cdot x)
=\max_{q\in S^{2}:\,(O^\intercal q)\cdot v=0} ((O^\intercal q)\cdot(O^\intercal x))
=\max_{q\in S^{2}:\,q\cdot v=0} (q\cdot(O^\intercal x)).
\end{align*}
The desired result now follows from $|O^\intercal x-x|\leq 2\sqrt{2}\sin(\beta/2)|x|\leq\sqrt{2}\beta|x|$. \qed

\medskip

\noindent Next, for any $t,\varepsilon>0$, we choose (via measurable selection) a random vector $v^*$ satisfying
\begin{align} \label{eq.sec2.prop2.1.step2.vstar}
v^*\in\text{argmax}_{v\in S^2} \max_{0\le s\le t}\;(v\cdot W_s)
\end{align}
and define, for any $\beta_0>0$, the (random) set $\Sigma\subset S^2$ as the set of all $v'\in S^2$ that form angles of at most $\beta_0$ with $v^*$.
Then, it is easy to see that there exist universal constants $\beta_0,c>0$ such that for any $t,\varepsilon>0$, the following holds almost surely on $A_\varepsilon^c$:
\begin{align*}
& \max_{0\le s\le t}\;(v\cdot W_s) \geq c\varepsilon,\quad v\in\Sigma.
\end{align*}
Notice that $c_{\beta_0}:=\sigma(\Sigma)$ is deterministic and that 
\begin{align}\label{cepsilon}
& r_{**} \leq \Lambda_{0} - c\varepsilon,\quad v\in\Sigma.
\end{align}

\smallskip

Using \eqref{cepsilon}, Lemma \ref{le:rstar.main}, Lemma \ref{le:SupOrthogonalPerturbation}, and Corollary \ref{cor:2.4} we infer that, for any small enough $T,\varepsilon_2>0$, the following holds for all small enough $\varepsilon_1,t>0$, almost surely on $A_{\varepsilon_2}\setminus A_{\varepsilon_1}$:
\begin{equation}\label{eq.secondTerm.bound}
\begin{split}
&\; C_1 \int_{S^{2}}\!\int_{r_*}^{r_{**}}\! f(r) r^{2}\,\mathrm{d}r\,\sigma(\mathrm{d}v)
\!=\! C_1 \int_{\Sigma\cup(-\Sigma)} \int_{r_*}^{r_{**}} \! f(r) r^{2}\,\mathrm{d}r\,\sigma(\mathrm{d}v)
\!+\! C_1\int_{S^{2}\setminus(\Sigma\cup(-\Sigma))} \int_{r_*}^{r_{**}}\! f(r) r^{2}\,\mathrm{d}r\,\sigma(\mathrm{d}v) \\
&\leq \! C_1 \sigma(\Sigma)\bigg(1\!-\!\frac{\psi(c \varepsilon_1)}{2}\bigg)
\!\bigg(\!\max_{0\le s\le t} \frac{(\Lambda^3_s\!-\!\widetilde\Lambda^3_s)_+}{3} 
\!+\!\max_{0\le s\le t}\,(\Lambda_s\!-\!\widetilde\Lambda_s)_+\,
O\big(\!\big(\max_{q\in S^2:\,q\cdot v^*=0}\max_{0\le s\le t}\,(q\cdot W_s \!+\! \sqrt{2}\beta_0)\big)^2\big)\!\!\bigg) \\
&\quad\! +C_1 \sigma(-\Sigma)\,\bigg(\max_{0\le s\le t} \frac{(\Lambda^3_s-\widetilde\Lambda^3_s)_+}{3} + \max_{0\le s\le t}\, (\Lambda_s-\widetilde\Lambda_s)_+\,O\big(\!\big(\max_{q\in S^2:\,q\cdot v^*=0}\,\max_{0\le s\le t}\;(q\cdot W_s + \sqrt{2}\beta_0)\big)^2\big)\!\bigg) \\
&\quad\! +C_1 \sigma(S^2\setminus(\Sigma\cup(-\Sigma)))\,\max_{0\le s\le t} \frac{
(\Lambda^3_s-\widetilde\Lambda^3_s)_+}{3} \\
& \leq \max_{0\le s\le t} \frac{(\Lambda^3_s-\widetilde\Lambda^3_s)_+}{3}
- C c_{\beta_0}\, \max_{0\le s\le t} \frac{(\Lambda^3_s-\widetilde\Lambda^3_s)_+}{3}\,
\big(c'\psi(c\varepsilon_1)-\max_{q\in S^2:\,q\cdot v^*=0}\,\max_{0\le s\le t}\;(q\cdot W_s)^2 - \beta^2_0\big),
\end{split}
\end{equation}
where $C,c'>0$ are universal constants.
Thus, the second term on the right-hand side of \eqref{eq.first} admits the bound
\begin{equation}\label{eq.sec2.prop2.1.step2.result}
\begin{split}
& \max_{0\le s\le t} \frac{(\Lambda^3_s-\widetilde\Lambda^3_s)_+}{3}\, \PP(A_{\varepsilon_2}\setminus A_{\varepsilon_1}) \\
&- C c_{\beta_0}\,\max_{0\le s\le t} \frac{(\Lambda^3_s-\widetilde\Lambda^3_s)_+}{3}\,\EE\big[ \bone_{A_{\varepsilon_2}\setminus A_{\varepsilon_1}} \big(c'\psi(c\varepsilon_1)-\max_{q\in S^2:\,q\cdot v^*=0}\,\max_{0\le s\le t} \;(q\cdot W_s)^2 - \beta^2_0\big)\big].
\end{split}
\end{equation}

\smallskip

\noindent\textbf{Step 3.} In this step, we estimate the third term on the right-hand side of \eqref{eq.first}.

\begin{lemma}\label{le:GaussianTail}
There exists a constant $\widecheck C>0$ such that, for all $t,\varepsilon_2>0$ small enough,
\begin{align*}
& \,C_1\EE\bigg[\int_{B_{\Lambda_{0}}(0)\setminus B_{r_{0}}(0)} f(|z|)\,\bone_{\{\max_{0\le s\le t\wedge\tau_{r_0}} (|z+W_s|-\widetilde{\Lambda}_s)\geq0,\,\max_{0\le s\le t\wedge\tau_{r_0}} (|z+W_s|-\Lambda_s)<0\}}\,\mathrm{d}z\,\bigg|\,A^c_{\varepsilon_2}\bigg] \\
& \leq \widecheck C\max_{0\le s\le t} \frac{(\Lambda^3_s-\widetilde\Lambda^3_s)_+}{3}.
\end{align*}
\end{lemma}

\noindent\textbf{Proof.} For this proof only, we introduce the notations
\begin{align*}
\tau=\inf\{s\geq 0:\,|W_s| \geq \varepsilon_{2}\}, 
\quad \sigma^z=\inf\{s\geq 0:\,|z+W_s| \geq \widetilde{\Lambda}_{s}\}, 
\quad R_s^z=|z+W_s|.
\end{align*}
With these and thanks to $\tau_{r_0}>s$ on  $\{\max_{[0,s]}\,(R^z-\widetilde{\Lambda})\ge0\}\cap\{\tau=s\}$, it holds
\begin{align}
&\;\mathbb{E}\bigg[\mathbf{1}_{A_{\varepsilon_{2}}^{c}} \int_{B_{\Lambda_{0}}(0)\setminus B_{r_{0}}(0)} f(|z|)\,\bone_{\{\max_{0\le s\le t\wedge\tau_{r_0}} (|z+W_s|-\widetilde{\Lambda}_s)\geq0,\,\max_{0\le s\le t\wedge\tau_{r_0}} (|z+W_s|-\Lambda_s)<0\}}\,\mathrm{d}z\bigg]\nonumber\\
& =\int_{0}^{t} \int_{B_{\Lambda_{0}}(0)\setminus B_{r_{0}}(0)} f(|z|)\,\mathbb{P}\big(\max_{[0,t]}\,(R_{\cdot\wedge\tau_{r_0}}^{z}-\widetilde{\Lambda}) \geq 0,\,\max_{[0,t]}\,(R_{\cdot\wedge\tau_{r_0}}^{z}-\Lambda)<0\,\big|\,\tau=s\big)\,\mathrm{d}z\,
\mathbb{P}(\tau \in \mathrm{d}s)\nonumber\\
& \leq \int_0^t \int_{B_{\Lambda_{0}}(0)\setminus B_{r_{0}}(0)} f(|z|)\,\mathbb{P}\big(\max_{[0, t]}\,(R_{\cdot\wedge\tau_{r_0}}^z\!-\!\widetilde{\Lambda}) \!\geq\! 0,\,\max _{[0,t]}\, (R_{\cdot\wedge\tau_{r_0}}^z\!-\!\Lambda)\!<\!0,\,\sigma^z \!\leq\! \tau\,\big|\, \tau\!=\!s\big)
\,\mathrm{d}z\,\mathbb{P}(\tau \in \mathrm{d}s)\nonumber\\
&\quad+\int_{0}^{t} \int_{B_{\Lambda_{0}}(0)\setminus B_{r_{0}}(0)} f(|z|)\,\mathbb{P}\big(\max_{[s, t]}\,(R_{\cdot\wedge\tau_{r_0}}^z\!-\!\widetilde{\Lambda})\!\geq\! 0,\,\max _{[s,t]}\,(R_{\cdot\wedge\tau_{r_0}}^z\!-\!\Lambda)\!<\!0,\,\sigma^z\!>\!\tau\,\big|\, \tau\!=\! s\big)\,\mathrm{d}z\, \mathbb{P}(\tau \in\mathrm{d}s) \nonumber\\
& \leq \int_0^t \int_{B_{\Lambda_{0}}(0)\setminus B_{r_{0}}(0)} f(|z|)\,\mathbb{P}\big(\max_{[0,s]}\,(R^z-\widetilde{\Lambda}) \geq 0,\, \max_{[0, s]}\, (R^z-\Lambda)<0\,\big|\,\tau=s\big)\,\mathrm{d}z\,
\mathbb{P}(\tau\in\mathrm{d}s) \nonumber \\
&\quad +\int_0^t \int_{B_{\Lambda_{0}}(0)\setminus B_{r_{0}}(0)} f(|z|)\,\mathbb{P}\big(\max_{[s,t]}\, (R_{\cdot\wedge\tau_{r_0}}^z-\widetilde{\Lambda}) \geq 0,\, \max _{[s,t]}\, (R_{\cdot\wedge\tau_{r_0}}^{z}-\Lambda)<0\,\big|\,\tau=s\big)\,\mathrm{d}z\, \mathbb{P}(\tau\in\mathrm{d}s) \nonumber \\
& =\int_{0}^{t} \mathbb{E}\bigg[\int_{B_{\Lambda_{0}}(0)\setminus B_{r_{0}}(0)} f(|z|)\, \mathbf{1}_{\{\max_{[0,s]} (|z+W|-\widetilde{\Lambda}) \geq 0,\, \max_{[0,s]} (|z+W|-\Lambda)<0\}}\,\mathrm{d}z\,\bigg|\,\tau=s\bigg]\,\mathbb{P}(\tau\in\mathrm{d}s)
\label{eq.sec2.prop2.1.step3.eq.1} \\
&\quad +\int_0^t \int_{B_{\Lambda_{0}}(0)\setminus B_{r_{0}}(0)} f(|z|)\,\mathbb{P}\big(\max_{[s, t]}\,(R_{\cdot\wedge\tau_{r_0}}^{z}-\widetilde{\Lambda}) \geq 0,\, \max_{[s, t]}\,(R_{\cdot\wedge\tau_{r_0}}^z-\Lambda)<0\,\big|\,\tau=s\big)\,\mathrm{d}z
\,\mathbb{P}(\tau\in\mathrm{d}s). \nonumber
\end{align}
On $\{\tau=s\}\subset\{\max_{[0,s]} |W|=\varepsilon_2\}$, we have 
\begin{equation}\label{eq.sec2.prop2.1.step3.eq.2}
\begin{split}
&\,\int_{B_{\Lambda_{0}}(0)\setminus B_{r_{0}}(0)} f(|z|)\,\mathbf{1}_{\{\max _{[0,s]} (|z+W|-\widetilde{\Lambda}) \geq 0, \, \max_{[0,s]} (|z+W|-\Lambda)<0\}}\,\mathrm{d}z \\
&\leq \int_{S^2} \int_{r_*}^{r_{**}} f(r)\,r^2\,\mathrm{d}r\,\sigma(\mathrm{d}v)
\leq \frac{1}{3C_1} \max _{0\le s'\le s} (\Lambda_{s'}^3-\widetilde{\Lambda}_{s'}^3)_+, 
\end{split}
\end{equation}
where the last inequality is shown by repeating \eqref{eq.sec2.Prop21.step1.finalEq}. This takes care of the first term on the right-hand side of \eqref{eq.sec2.prop2.1.step3.eq.1}.

\smallskip

To estimate the second term on the right-hand side of \eqref{eq.sec2.prop2.1.step3.eq.1}, we notice that
\begin{align*}
&\,\mathbb{P}\big(\max_{[s,t]}\,(R_{\cdot\wedge\tau_{r_0}}^z-\widetilde{\Lambda}) \geq 0, \,\max_{[s,t]}\,(R_{\cdot\wedge\tau_{r_0}}^z-\Lambda)<0\,\big|\,\tau=s\big) \\
&\le\mathbb{P}\big(\max_{[s,t]}\, (R_{\cdot\wedge\tau^s_{r_0}}^z-\widetilde{\Lambda}) \geq 0,\, \max_{[s,t]}\, (R_{\cdot\wedge\tau^s_{r_0}}^z-\Lambda)<0\,\big|\,\tau=s\big) \\
&=\mathbb{P}\big(\max_{[0,t-s]}\, (\overline{R}_{\cdot\wedge\overline\tau_{r_0}}^{|z|}-\widetilde{\Lambda}_{s+\cdot}) \geq 0, 
\max_{[0,t-s]} (\overline{R}_{\cdot\wedge\overline\tau_{r_0}}^{|z|}-\Lambda_{s+\cdot})<0\big),
\end{align*}
where 
\begin{align*}
& \tau^s_{r_0}:=\inf\{s'\geq s:\,R^z_{s'}\leq r_0\},
\quad \overline\tau_{r_0}:=\inf\{s'\geq 0:\,\overline{R}^{|z|}_{s'} \leq r_0\},
\end{align*}
and $\overline{R}^{|z|}$ is a three-dimensional Bessel process with the initial probability density
$\phi(\cdot;|z|)$, given by the density of $R_\tau^z$. Therefore, 
\begin{align*}
& \,\int_{B_{\Lambda_{0}}(0)\setminus B_{r_{0}}(0)} f(|z|)\,\mathbb{P}\big(\max_{[s, t]}\,(R_{\cdot\wedge\tau_{r_0}}^{z}-\widetilde{\Lambda}) \geq 0,\, \max_{[s, t]}\,(R_{\cdot\wedge\tau_{r_0}}^z-\Lambda)<0\,\big|\,\tau=s\big)\,\mathrm{d}z \\
& \leq C_1^{-1}\int_{r_0}^{\Lambda_{0}} f(r)\,r^2\,\mathbb{P}\big(\max_{[0,t-s]}\, (\overline{R}_{\cdot\wedge\overline\tau_{r_0}}^r-\widetilde{\Lambda}_{s+\cdot}) \geq 0, 
\max_{[0,t-s]} (\overline{R}_{\cdot\wedge\overline\tau_{r_0}}^r-\Lambda_{s+\cdot})<0\big)\,\mathrm{d}r \\
& = \mathbb{P}\big(\max_{[0,t-s]}\, (\widehat{R}_{\cdot\wedge\widehat\tau_{r_0}}-\widetilde{\Lambda}_{s+\cdot}) \geq 0, \max_{[0, t-s]}\, (\widehat{R}_{\cdot\wedge\widehat\tau_{r_0}}-\Lambda_{s+\cdot})<0\big),
\end{align*}
with $\widehat\tau_{r_0}:=\inf\{s'\geq 0:\, \widehat{R}_{s'}\leq r_0\}$ and with a three-dimensional Bessel process $\widehat{R}$ of initial density
\begin{align*}
\chi(x):= C^{-1}_1 \int_{r_0}^{\Lambda_{0}} \phi(x;r)\,f(r)\,r^2\,\mathrm{d}r
= C^{-1}_1 \int_{(x-\varepsilon_2)\vee r_0}^{(x+\varepsilon_{2})\wedge\Lambda_0} \phi(x;r)\,f(r)\,r^2\,\mathrm{d}r,
\quad x\in[r_0-\varepsilon_2,\Lambda_0+\varepsilon_2].
\end{align*}

\smallskip

Next, we apply the representation of the three-dimensional Bessel process as a Doob $h$-transform of the standard Brownian motion stopped at $0$ and Fubini's theorem, to obtain from the above
\begin{equation}\label{eq.sec2.prop2.1.step3.eq.3}
\begin{split}
& \,\int_{B_{\Lambda_{0}}(0)\setminus B_{r_{0}}(0)} f(|z|)\,\mathbb{P}\big(\max_{[s, t]}\,(R_{\cdot\wedge\tau_{r_0}}^{z}-\widetilde{\Lambda}) \geq 0,\, \max_{[s, t]}\,(R_{\cdot\wedge\tau_{r_0}}^z-\Lambda)<0\,\big|\,\tau=s\big)\,\mathrm{d}z \\
& \le\mathbb{E}\bigg[\int_{r_0-\varepsilon_2}^{\Lambda_0+\varepsilon_2} \chi(x)\,\frac{x+B_{(t-s)\wedge\tau^B_0}}{x} \,\mathbf{1}_{\{\max_{[0, t-s]} (x+B_{\cdot\wedge\tau^B_{r_0}}-\widetilde{\Lambda}_{s+\cdot}) \geq 0,\,\max_{[0,t-s]} (x+B_{\cdot\wedge\tau^B_{r_0}}-\Lambda_{s+\cdot})<0\}}
\,\mathrm{d}x\bigg],
\end{split}
\end{equation}
where $\tau^B_0:=\inf\{s'\geq 0\!:x+B_{s'}\le 0\}$ and $\tau^B_{r_0}:=\inf\{s'\geq 0\!:x+B_{s'} \le r_0\}$. We now estimate $\chi$ by bounding $\phi$. When $\phi(\cdot;r)>0$, it is proportional to the radius of the circle given by the intersection of the spheres $\partial B_x((0,0,0))$ and $\partial B_{\varepsilon_2}((r,0,0))$, which computes (e.g., by Heron's formula) to 
\begin{align*}
&\frac{(r+x+\varepsilon_2)^{1/2}(r+x-\varepsilon_2)^{1/2}(r+\varepsilon_2-x)^{1/2}(x+\varepsilon_2-r)^{1/2}}{2r}.
\end{align*}
For $r\in[r_0,\Lambda_0]$, the latter lies between two positive multiples of $(r+\varepsilon_2-x)^{1/2}(x+\varepsilon_2-r)^{1/2}$, thus,
\begin{align*}
\phi(x;r)=O\big(\varepsilon^{-2}_2(r+\varepsilon_2-x)^{1/2}(x+\varepsilon_2-r)^{1/2}\big).
\end{align*}
Consequently, $\chi$ can be estimated by a uniform constant.
~Finally, we consider the random variables
\begin{align*}
&\overline X:=\sup\big\{y\in[r_0-\varepsilon_2,\Lambda_0+\varepsilon_2]:\,\max_{[0,t-s]}\,(x+B_{\cdot\wedge\tau^B_{r_0}}-\Lambda_{s+\cdot})<0\big\}\vee (r_0-\varepsilon_2), \\
&\underline X:=\inf\big\{y\in[r_0-\varepsilon_2,\Lambda_0+\varepsilon_2]:\,\max_{[0,t-s]}\,(x+B_{\cdot\wedge\tau^B_{r_0}}-\widetilde\Lambda_{s+\cdot})\ge0\big\}. 
\end{align*}
In view of $\overline X-\underline X\le\max_{0\le s'\le t}\,(\Lambda_{s'}-\widetilde\Lambda_{s'})_+$, we arrive at
\begin{equation}\label{eq.sec2.prop2.1.step3.eq.4}
\begin{split}
&\,\mathbb{E}\bigg[\int_{r_0-\varepsilon_2}^{\Lambda_0+\varepsilon_2} \chi(x)\,\frac{x+B_{(t-s)\wedge\tau^B_0}}{x} \,\mathbf{1}_{\{\max_{[0, t-s]} (x+B_{\cdot\wedge\tau^B_{r_0}}-\widetilde{\Lambda}_{s+\cdot}) \geq 0,\,\max_{[0,t-s]} (x+B_{\cdot\wedge\tau^B_{r_0}}-\Lambda_{s+\cdot})<0\}}
\,\mathrm{d}x\bigg] \\
&=\mathbb{E}\bigg[\int_{\underline X}^{\overline X} \chi(x)\,\frac{x+B_{(t-s)\wedge\tau^B_0}}{x} \,\mathrm{d}x\bigg]
=O\big(\max_{s\le s'\le t}\,(\Lambda_{s'}-\widetilde\Lambda_{s'})_+\big).
\end{split}
\end{equation}

\smallskip

Collecting \eqref{eq.sec2.prop2.1.step3.eq.1}--\eqref{eq.sec2.prop2.1.step3.eq.4} we end up with the desired  
\begin{align*}
& \,\mathbb{E}\bigg[\mathbf{1}_{A_{\varepsilon_{2}}^{c}} \int_{B_{\Lambda_{0}}(0)\setminus B_{r_{0}}(0)} f(|z|)\,\bone_{\{\max_{0\le s\le t\wedge\tau_{r_0}} (|z+W_s|-\widetilde{\Lambda}_s)\geq0,\,\max_{0\le s\le t\wedge\tau_{r_0}} (|z+W_s|-\Lambda_s)<0\}}\,\mathrm{d}z\bigg] \\
&\le\PP(A_{\varepsilon_{2}}^{c})\, \widecheck C\max_{0\le s\le t} \frac{(\Lambda^3_s-\widetilde\Lambda^3_s)_+}{3}. 
\qquad\qquad\qquad\qquad\qquad\qquad\qquad\qquad\qquad\qquad\qquad\qquad\qquad\quad\;\;\; \qed
\end{align*}

\smallskip

Putting together \eqref{eq.first}, \eqref{eq.sec2.prop2.1.step1.result}, \eqref{eq.sec2.prop2.1.step2.result} and Lemma \ref{le:GaussianTail}  we find that, for any small enough $T,\varepsilon_2>0$, the following holds for all small enough $\varepsilon_1,t>0$:
\begin{align*}
& \,\int_{r_0}^{\Lambda_{0}} f(x)\,\PP^x\big(\widetilde{\tau}\leq t\wedge\tau_{r_0},\,\tau>t\wedge\tau_{r_0}\big)\,\nu(\mathrm{d}x) \\
& \le \PP(A_{\varepsilon_2})\,\frac13\max_{0\le s\le t} (\Lambda_s^3-\widetilde\Lambda_s^3)_+
-\EE\big[\bone_{A_{\varepsilon_1}}\,\psi\big(\max_{0\le s\le t} (W_s^{(1)}-\Lambda_s+\Lambda_{0})\big)\big]\,\frac{1}{18}
\max_{0\le s\le t} (\Lambda_s^3-\widetilde\Lambda_s^3)_+ \\
&\quad - \max_{0\le s\le t} \frac{(\Lambda^3_s-\widetilde\Lambda^3_s)_+}{3} \,
\Big(C c_{\beta_0}\,\EE\big[ \bone_{A_{\varepsilon_2}\setminus A_{\varepsilon_1}} \big(c'\psi(c\varepsilon_1)-\max_{q\in S^2:\,q\cdot v^*=0}\,\max_{0\le s\le t} \;(q\cdot W_s)^2 - \beta^2_0\big)\big]-\widecheck{C}\PP(A^c_{\varepsilon_2})\Big),
\end{align*}
whereby we recall the definition of $v^*$ from \eqref{eq.sec2.prop2.1.step2.vstar}.

\medskip

\noindent\textbf{Step 4.} At this point, to complete the proof of the proposition it is enough to show that, for any small enough $T,\varepsilon_2>0$ and any $\varepsilon_1\in(0,\varepsilon_2)$, 
\begin{align}
&\EE\big[\bone_{A_{\varepsilon_1}}\,\psi\big(\max_{0\le s\le t} (W_s^{(1)}-\Lambda_s+\Lambda_{0})\big)\big]
\ge\frac12\EE\big[\psi\big(\max_{0\le s\le t} (B_s-\Lambda_s+\Lambda_0)\wedge\varepsilon_2\big)\big], \label{Gaussian tail} \\
& \widecheck C\,\frac{\PP(A^c_{\varepsilon_2})}{\PP(A_{\varepsilon_2}\setminus A_{\varepsilon_1})} \leq  Cc_{\beta_0}\,(c'\psi(c\varepsilon_1) - \beta^2_0)
-Cc_{\beta_0} \,\EE\big[\max_{q\in S^2:\,q\cdot v^*=0}\,\max_{0\le s\le t} \;(q\cdot W_s)^2
\,\big|\,A_{\varepsilon_2}\setminus A_{\varepsilon_1}\big]\nonumber
\end{align}
for all small enough $t>0$. We show the former in Lemma \ref{lem2.9}, while the latter follows upon decreasing $\beta_0>0$ (to ensure $c'\psi(c\varepsilon_1) - \beta^2_0>0$) and using the fast tail decay exhibited by the law of $\max_{0\le s\le t} |W_s|$, as well as Lemma \ref{lem2.10}.

\begin{lemma}\label{lem2.9}
For any small enough $T,\varepsilon_2>0$ and any $\varepsilon_1\in(0,\varepsilon_2)$, the inequality \eqref{Gaussian tail} holds for all sufficiently small $t>0$.
\end{lemma}

\begin{lemma}\label{lem2.10}
For any $0<\varepsilon_1<\varepsilon_2$,
\begin{align}\label{eq.last.issue}
& \lim_{t\downarrow0}\,\EE\big[\max_{q\in S^2:\,q\cdot v^*=0}\,\max_{0\le s\le t} \;(q\cdot W_s)^2
\,\big|\,A_{\varepsilon_2}\setminus A_{\varepsilon_1}\big] = 0.
\end{align}
\end{lemma}

\noindent\textbf{Proof of Lemma \ref{lem2.9}.} Consider any small enough $T,\varepsilon_2>0$, such that \eqref{psi bound} holds. We start by estimating the left-hand side of \eqref{Gaussian tail}:
\begin{align*}
&\;\EE\big[\bone_{A_{\varepsilon_1}}\,\psi\big(\max_{0\le s\le t} (W_s^{(1)}-\Lambda_s+\Lambda_{0})\big)\big] \\
&=\EE\big[\bone_{A_{\varepsilon_1}}\,\psi\big(\max_{0\le s\le t} (W_s^{(1)}-\Lambda_s+\Lambda_{0})\wedge\varepsilon_2\big)\big] \\
&=\EE\big[\psi\big(\max_{0\le s\le t} (W_s^{(1)}-\Lambda_s+\Lambda_{0})\wedge\varepsilon_2\big)\big]
-\EE\big[\bone_{A_{\varepsilon_1}^c}\,\psi\big(\max_{0\le s\le t} (W_s^{(1)}-\Lambda_s+\Lambda_{0})\wedge\varepsilon_2\big)\big] \\
&\ge\EE\big[\psi\big(\max_{0\le s\le t} (W_s^{(1)}-\Lambda_s+\Lambda_{0})\wedge\varepsilon_2\big)\big]
-\sup_{0\le x\le\varepsilon_2} \psi(x)\,\PP(A_{\varepsilon_1}^c).
\end{align*}
On the other hand,
\begin{align*}
\EE\big[\psi\big(\max_{0\le s\le t} (B_s-\Lambda_s+\Lambda_0)\wedge\varepsilon_2\big)\big]
\ge\inf_{\varepsilon_1/4\le x\le \varepsilon_2}\psi(x)\,\PP\big(\max_{0\le s\le t} B_s\ge\varepsilon_1/4).    
\end{align*}
The desired result follows from the fact that, for all sufficiently small $t>0$, 
\begin{align*}
\sup_{0\le x\le\varepsilon_2} \psi(x)\,\PP(A_{\varepsilon_1}^c)
\leq\frac12\,\inf_{\varepsilon_1/4\le x\le \varepsilon_2}\psi(x)\,\PP\big(\max_{0\le s\le t} B_s\ge\varepsilon_1/4),
\end{align*}
which is due to $\PP(A_{\varepsilon_1}^c)\le 6\PP(\max_{0\le s\le t} B_s\ge\varepsilon_1/2)$ (easily obtained by a union bound). \qed

\medskip

\noindent\textbf{Proof of Lemma \ref{lem2.10}.} For this proof only, we introduce the notations
\begin{align*}
& \tau=\inf\{s\geq0:\,|W_s|\geq \varepsilon_1\},
\quad \widecheck v^*= W_\tau/|W_\tau|= W_\tau/\varepsilon_1
\in\text{argmax}_{v\in S^2} \max_{0\le s\le\tau} (v\cdot W_s).
\end{align*}
We also consider a random time $\widecheck\tau\in[0,t]$ satisfying $|W_{\widecheck\tau}|=\max_{0\le s\le t} |W_s|$, as well as $v^*:=\frac{W_{\widecheck\tau}}{|W_{\widecheck\tau}|}$. Let $\beta\in[0,\frac\pi2]$ be the angle between the lines spanned by $\widecheck v^*$ and by $v^*$. Then, on $A_{\varepsilon_1}^c$,
\begin{align*}
\max_{\tau\le s\le t} |W_s\!-\!W_\tau|^2
\!\ge\!|W_{\widecheck\tau}\!-\!W_\tau|^2
\!=\!|W_{\widecheck\tau}|^2\!+\!|W_\tau|^2\!\pm\!2|W_{\widecheck\tau}||W_\tau|\cos\beta
&\!=\!|W_\tau|^2\sin^2\beta\!+\!(|W_\tau|\cos\beta\!\pm\!|W_{\widecheck\tau}|)^2 \\
&\!\ge\!|W_\tau|^2\sin^2\beta\!\ge\!\frac{4}{\pi^2}\max_{0\le s\le\tau} |W_s|^2\,\beta^2.    
\end{align*}
Next, we note that $\max_{0\le s\le t} |q\cdot W_s|\le \max_{0\le s\le\tau} |q\cdot W_s|+\max_{\tau\le s\le t} |q\cdot(W_s-W_\tau)|$, $q\in S^2$ on~$A_{\varepsilon_1}^c$, 
by distinguishing whether the maximum on the left-hand side is attained on $[0,\tau]$ or on $[\tau,t]$. Together with Lemma \ref{le:SupOrthogonalPerturbation} this yields, on $A_{\varepsilon_1}^c$,
\begin{align*}
\max_{q\in S^2:\,q\cdot v^*=0}\,\max_{0\le s\le t}\;(q\cdot W_s)^2
&\le 2\max_{q\in S^2:\,q\cdot v^*=0}\,\max_{0\le s\le\tau}\;(q\cdot W_s)^2
+2\max_{\tau\le s\le t} |W_s-W_\tau|^2 \\
&\le 4\max_{q\in S^2:\,q\cdot\widecheck v^*=0}\,\max_{0\le s\le\tau}\;(q\cdot W_s)^2
+8\beta^2\max_{0\le s\le \tau} |W_s|^2+2\max_{\tau\le s\le t} |W_s-W_\tau|^2 \\
&\le 4\max_{q\in S^2:\,q\cdot\widecheck v^*=0}\,\max_{0\le s\le\tau}\;(q\cdot W_s)^2
+(2\pi^2+2)\max_{\tau\le s\le t} |W_s-W_\tau|^2.
\end{align*}
As a result, we obtain
\begin{align*}
&\;\EE\big[\max_{q\in S^2:\,q\cdot v^*=0}\,\max_{0\le s\le t}\;(q\cdot W_s)^2\,\bone_{A_{\varepsilon_2}\setminus A_{\varepsilon_1}}\big] \\
&\leq 4\EE\big[\max_{q\in S^2:\,q\cdot\widecheck v^*=0}\,\max_{0\le s\le\tau}\;(q\cdot W_s)^2\,\bone_{A_{\varepsilon_2}\setminus A_{\varepsilon_1}}\big]
+(2\pi^2+2)\EE\big[\max_{\tau\le s\le t} |W_s-W_\tau|^2\,\,\bone_{A_{\varepsilon_2}\setminus A_{\varepsilon_1}}\big].
\end{align*}

\smallskip

We proceed by introducing 
\begin{align*}
&\widecheck A^c_{\varepsilon_1}= \{|W_t|\geq\varepsilon_1\} \subset A^c_{\varepsilon_1}
\end{align*}
and using the strong Markov property of Brownian motion to deduce 
\begin{align*}
&\;\EE\big[\max_{q\in S^2:\,q\cdot\widecheck v^*=0}\,\max_{0\le s\le\tau}\;(q\cdot W_s)^2\,\bone_{A^c_{\varepsilon_1}\setminus{\widecheck A^c_{\varepsilon_1}}}\big]
=\EE\big[\max_{q\in S^2:\,q\cdot\widecheck v^*=0}\,\max_{0\le s\le\tau}\;(q\cdot W_s)^2\,\EE[\bone_{A^c_{\varepsilon_1}\setminus{\widecheck A^c_{\varepsilon_1}}}\,\vert\,\mathcal{F}^W_\tau]\big] \\
&\qquad\qquad\qquad\qquad\qquad\qquad\qquad\qquad\quad
\leq \frac{1}{2} \EE\big[\max_{q\in S^2:\,q\cdot\widecheck v^*=0}\,\max_{0\le s\le\tau}\;(q\cdot W_s)^2\,\bone_{A^c_{\varepsilon_1}}\big] \\
& \Longrightarrow\;\EE\big[\max_{q\in S^2:\,q\cdot\widecheck v^*=0}\,\max_{0\le s\le\tau}\;(q\cdot W_s)^2\,\bone_{\widecheck A^c_{\varepsilon_1}}\big] \\
&\qquad\, = \EE\big[\max_{q\in S^2:\,q\cdot\widecheck v^*=0}\,\max_{0\le s\le\tau}\;(q\cdot W_s)^2\,\bone_{A^c_{\varepsilon_1}}\big] 
-\EE\big[\max_{q\in S^2:\,q\cdot\widecheck v^*=0}\,\max_{0\le s\le\tau}\;(q\cdot W_s)^2\,\bone_{A^c_{\varepsilon_1}\setminus{\widecheck A^c_{\varepsilon_1}}}\big] \\
&\qquad\, \geq \frac{1}{2}  \EE\big[\max_{q\in S^2:\,q\cdot\widecheck v^*=0}\,\max_{0\le s\le\tau}\;(q\cdot W_s)^2\,\bone_{A^c_{\varepsilon_1}}\big].
\end{align*}
In view of $\lim_{t\downarrow0}\,\frac{\PP(A_{\varepsilon_2}\setminus A_{\varepsilon_1})}{\PP(A_{\varepsilon_1}^c)}=1$, the latter inequality, and
\begin{align*}
	\EE\big[\max_{\tau\le s\le t} |W_s-W_\tau|^2\,\,\bone_{A_{\varepsilon_2}\setminus A_{\varepsilon_1}}\big]
	&=\EE\big[\EE\big[\max_{\tau\le s\le t} |W_s-W_\tau|^2\,\,\bone_{A_{\varepsilon_2}\setminus A_{\varepsilon_1}}\,\big|\,\mathcal{F}^W_{\tau}\big]\big] \\
	&\le\EE\big[\bone_{A^c_{\varepsilon_1}}\,
	\EE\big[\max_{\tau\le s\le t} |W_s-W_\tau|^2\,\big|\,\mathcal{F}^W_{\tau}\big]\big]
	\leq 24t\,\PP(A^c_{\varepsilon_1}),
\end{align*}
one can reduce \eqref{eq.last.issue} to
\begin{align}\label{reduced}
&\lim_{t\downarrow0}\,\frac{\EE\big[\max_{q\in S^2:\,q\cdot\widecheck v^*=0}\,\max_{0\le s\le\tau}\;(q\cdot W_s)^2\,\bone_{\widecheck A^c_{\varepsilon_1}}\big]}{\PP(\widecheck A^c_{\varepsilon_1})}=0.
\end{align}

\smallskip

To see \eqref{reduced}, consider $\overline v^* = W_t/|W_t|$ and let $\overline{\beta}\in[0,\frac\pi2]$ be the angle between the lines spanned by $\overline v^*$ and by $\widecheck v^*$. Then, as above, we infer that
\begin{align*}
|W_t\!-\!W_\tau|^2=|W_t|^2+|W_\tau|^2\pm2|W_t||W_\tau|\cos\overline\beta
=\varepsilon_1^2\sin^2\overline\beta+(|W_t|\!\pm\!\varepsilon_1\cos\overline\beta)^2
\ge\varepsilon_1^2\sin^2\overline\beta\ge\frac{4}{\pi^2}\varepsilon_1^2\overline\beta^2  
\end{align*}
on $\widecheck A^c_{\varepsilon_1}$, and in conjunction with Lemma \ref{le:SupOrthogonalPerturbation},
\begin{align*}
\max_{q\in S^2:\,q\cdot\widecheck v^*=0}\,\max_{0\le s\le\tau}\;(q\cdot W_s)^2
&\le 2\max_{q\in S^2:\,q\cdot\overline v^*=0}\,\max_{0\le s\le\tau}\;(q\cdot W_s)^2
+4\varepsilon_1^2\overline\beta^2 \\ 
&\leq 2\max_{q\in S^2:\,q\cdot\overline v^*=0}\,\max_{0\le s\le\tau}\;(q\cdot W_s)^2
+\pi^2|W_t-W_\tau|^2
\end{align*}
on $\widecheck A^c_{\varepsilon_1}$. Consequently,
\begin{align*}
\EE\big[\max_{q\in S^2:\,q\cdot\widecheck v^*=0}\,\max_{0\le s\le\tau}\;(q\cdot W_s)^2\,\bone_{\widecheck A^c_{\varepsilon_1}}\big]
&\le2\EE\big[\max_{q\in S^2:\,q\cdot\overline v^*=0}\,\max_{0\le s\le\tau}\;(q\cdot W_s)^2\,\bone_{\widecheck A^c_{\varepsilon_1}}\big]
+ \pi^2\EE\big[|W_t-W_\tau|^2\,\bone_{A^c_{\varepsilon_1}}\big] \\
& \leq 2\EE\big[\max_{q\in S^2:\,q\cdot\overline v^*=0}\,\max_{0\le s\le t}\;(q\cdot W_s)^2\,\bone_{\widecheck A^c_{\varepsilon_1}}\big]
+24\pi^2t\,\PP(\widecheck A^c_{\varepsilon_1}).
\end{align*}
Further, we have 
\begin{align*}
\EE\big[\max_{q\in S^2:\,q\cdot\overline v^*=0}\,\max_{0\le s\le t}\;(q\cdot W_s)^2
\,\bone_{\widecheck A^c_{\varepsilon_1}}\big]
&=\int_{\RR^3\backslash B_{\varepsilon_1}(0)} \EE\big[\max_{q\in S^2:\,q\cdot\overline v^*=0}\,\max_{0\le s\le t}\;(q\cdot W_s)^2\,\big|\,W_t=z\big]\,
\PP(W_t\in\mathrm{d}z) \\
&=\int_{\RR^3\backslash B_{\varepsilon_1}(0)} \EE\big[\max_{q\in S^2:\,q\cdot z=0}\,\max_{0\le s\le t}\;(q\cdot W_s)^2\,\big|\,W_t=z\big]\,
\PP(W_t\in\mathrm{d}z).
\end{align*}
Due to the rotational symmetry of Brownian motion, it suffices to consider $z=(r,0,0)$:
\begin{align*}
\EE\big[\max_{q\in S^2:\,q\cdot z=0}\,\max_{0\le s\le t}\;(q\cdot W_s)^2\,\big|\,W_t\!=\!z\big]
&=\EE\big[ \max_{0\le s\le t} \big((W^{(2)}_s)^2\!+\!(W^{(3)}_s)^2\big)\,\big|\,W^{(1)}_t\!=\!r,\,W^{(2)}_t\!=\!0,\,W^{(3)}_t\!=\!0
\big] \\
& = \EE\big[ \max_{0\le s\le t} \big((W^{(2)}_s)^2+(W^{(3)}_s)^2\big)\,\big|\,W^{(2)}_t=0,\,W^{(3)}_t=0
\big]
\leq 16t.
\end{align*}
In particular, 
\begin{align*}
&\EE\big[\max_{q\in S^2:\,q\cdot\overline v^*=0}\,\max_{0\le s\le t}\;(q\cdot W_s)^2
\,\bone_{\widecheck A^c_{\varepsilon_1}}\big]
\le 16t\,\PP(\widecheck A^c_{\varepsilon_1}),
\end{align*}
which completes the proof of \eqref{reduced}. \qed

\subsection{Proof of Proposition \ref{prop 2nd phase}}

To ease the notation, we take $t_0=0$ throughout the proof.

\medskip

\noindent\textbf{Step 1.} To bound the integral in \eqref{eq psi 2nd phase} we couple the underlying Bessel processes according to
\begin{equation}
\mathrm{d}R^x_s = \frac{1}{R^x_s}\,\mathrm{d}s + \mathrm{d}B_s,\quad R^x_0=x.
\end{equation}
Notice that $R^{\underline{x}}_s\le R^{\overline{x}}_s$, $s\in[0,t]$, $\Lambda_{0}\le \underline{x}\le \overline{x}$, and also that 
\begin{equation}\label{gap}
	\mathrm{d}R^{\overline{x}}_s-\mathrm{d}R^{\underline{x}}_s = -\frac{R^{\overline{x}}_s-R^{\underline{x}}_s}{R^{\overline{x}}_s R^{\underline{x}}_s}\,\mathrm{d}s,\quad R^{\overline{x}}_0-R^{\underline{x}}_0=\overline{x}-\underline{x},\quad\text{i.e.:}\quad
	R^{\overline{x}}_s-R^{\underline{x}}_s = (\overline{x}-\underline{x})e^{-\int_0^s (R^{\overline{x}}_{s'} R^{\underline{x}}_{s'})^{-1}\,\mathrm{d}s'}.
\end{equation}
Together with $f(\Lambda_{0}+y)(\Lambda_{0}+y)^2\le \overline{C}y$, $y\ge0$,
this yields
\begin{equation*}
\begin{split}
\int_{\Lambda_{0}}^\infty f(x)\,\PP^x(\widetilde{\tau}\le t,\tau>t)\,\nu(\mathrm{d}x) 
& \le \overline{C}\,\EE\bigg[\int_{\Lambda_{0}}^\infty (x-\Lambda_{0})\,\mathbf{1}_{\{\min_{0\le s\le t} (R^x_s-\widetilde{\Lambda}_{s})\le 0,\,\min_{0\le s\le t} (R^x_s-\Lambda_{s})>0\}}\,\mathrm{d}x\bigg] \\
& \le \overline{C}\,\EE\big[(\overline{X}-\Lambda_{0})(\overline{X}-\underline{X})_+\big],
\end{split}
\end{equation*}
where $\overline{X}$ and $\underline{X}$ are defined via
\begin{align*}
&0=\min_{0\le s\le t}\,(R^{\overline{X}}_s-\widetilde{\Lambda}_{s})
=\overline{X}
+\min_{0\le s\le t} \bigg(\int_0^s \frac{1}{R^{\overline{X}}_r}\,\mathrm{d}r
+B_s-\widetilde{\Lambda}_{s}\bigg)\;\Longrightarrow\;
\overline{X}\le \max_{0\le s\le t}\,(-B_s+\widetilde{\Lambda}_{s}), \\
&0=\min_{0\le s\le t}\,(R^{\underline{X}}_s-\Lambda_{s}).
\end{align*}
In view of \eqref{gap}, it holds 
\begin{equation*}
(\overline{X}-\underline{X})_+\le (\widetilde{\Lambda}_{\widetilde{\tau}^{\overline{X}}}
-\Lambda_{\widetilde{\tau}^{\overline{X}}})_+
\,e^{t/(\Lambda_t\widetilde{\Lambda}_t)}
\le \max_{0\le s\le t} |\Lambda_s-\widetilde{\Lambda}_s|\,
e^{t/(\Lambda_t\widetilde{\Lambda}_t)}.
\end{equation*}
Altogether, for any $t>0$, we see that the integral in \eqref{eq psi 2nd phase} cannot exceed
\begin{equation*}
\begin{split}
\overline{C}\,\EE\big[\max_{0\le s\le t}
\,(B_s+\widetilde{\Lambda}_{s})-\Lambda_{0}\big]
\,\max_{0\le s\le t} |\Lambda_s-\widetilde{\Lambda}_s|\,
e^{t/(\Lambda_t\widetilde{\Lambda}_t)}.
\end{split}
\end{equation*}
This proves Proposition \ref{prop 2nd phase} if $\psi(0+)>0$. Hence, for the rest of the proof we assume that $\psi(0+)=0$ and that $\Lambda,\,\widetilde\Lambda$ satisfy \eqref{eq.sec3.sergey.growth.cond.t0} for $t\in[0,T_0]$, with $t_0=0$ and with $f$ in place of $w(t_0,\cdot)$.

\smallskip

\noindent\textbf{Step 2.} Let us fix arbitrary $T\in(0,T_0]$, $0<\varepsilon_2\le(\Lambda_{T}\!\wedge\!\widetilde{\Lambda}_{T}-r_0)/2$ and consider the probabilistic formulation of the \textit{one-phase} radially symmetric Stefan problem with surface tension and with zero Dirichlet boundary condition at the radius $r_0$, equipped with the initial data $\widehat{\Lambda}_{0-}:=\Lambda_{0}$, $\widehat{w}(0-,x):=f(x)\,\mathbf{1}_{[\Lambda_T\wedge\widetilde{\Lambda}_T-\varepsilon_2,\,\Lambda_{0}]}(x)$:
\begin{equation}\label{one phase}
\begin{split}
&\frac{\widehat\Lambda_{0-}^3}{3}-\frac{\widehat{\Lambda}_s^3}{3} = \int_{\Lambda_T\wedge\widetilde{\Lambda}_T-\varepsilon_2}^{\widehat\Lambda_{0-}} f(x)\,\PP^x\big(\tau^-_{\widehat\Lambda_{0+\cdot}}\le s\wedge\tau_{r_0}\big)\,\nu(\mathrm{d}x),\quad s\in[0,T], \\
& \tau_{r_0}=\inf\{s\geq0:\,R_s\leq r_0\},\quad 
\tau^-_{\widehat\Lambda_{0+\cdot}}=\inf\{s\geq0:\,R_s\geq\widehat\Lambda_s\},\quad \PP^x\text{-a.s.},\,\,\,x\in(0,\widehat\Lambda_{0-}],
\end{split}
\end{equation}
where $R$ is a three-dimensional Bessel process.~A solution $(\widehat\Lambda,\widehat w)$ to \eqref{one phase} on $[0,T]$ can be constructed via the following iterative scheme: $\widehat\Lambda^0\equiv\widehat{\Lambda}_{0-}$, and for $n\geq1$,
\begin{align}
& \widehat\Lambda^n_{0-}=\widehat{\Lambda}_{0-},
\quad \frac{(\widehat\Lambda_{0-})^3}{3}-\frac{(\widehat{\Lambda}^n_s)^3}{3} = \int_{\Lambda_T\wedge\widetilde{\Lambda}_T-\varepsilon_2}^{\widehat\Lambda_{0-}} f(x)\,\PP^x\big(\tau^-_{\widehat\Lambda^{n-1}_{0+\cdot}}\le s\wedge\tau_{r_0}\big)\,\nu(\mathrm{d}x),\;\;s\in[0,T].\label{eq.sec2.prop2.2.step2.iterScheme}
\end{align}
The mapping given by the latter integral is non-increasing: For any $\Lambda'\leq \Lambda''$,
\begin{align*}
&\int_{\Lambda_T\wedge\widetilde{\Lambda}_T-\varepsilon_2}^{\widehat\Lambda_{0-}} f(x)\,\PP^x\big(\tau^-_{\widehat\Lambda_{0+\cdot}'}\le s\wedge\tau_{r_0}\big)\,\nu(\mathrm{d}x)
\geq \int_{\Lambda_T\wedge\widetilde{\Lambda}_T-\varepsilon_2}^{\widehat\Lambda_{0-}} f(x)\,\PP^x\big(\tau^-_{\widehat\Lambda_{0+\cdot}''} 
\le s\wedge\tau_{r_0}\big)\,\nu(\mathrm{d}x),\quad s\in[0,T].
\end{align*}
Using this monotonicity together with $\widehat\Lambda^1\leq\widehat\Lambda^0$ we deduce inductively that the sequence $(\widehat\Lambda^n)_{n\ge0}$ is pointwise non-increasing in $n$, and thus, has a pointwise limit.~The latter admits a non-increasing right-continuous modification, denoted by $\widehat{\Lambda}$ and extended constantly to $[-1,0)$ with $\widehat{\Lambda}_{0-}\!=\!\Lambda_{0}\geq\widehat{\Lambda}_{0}$. The verification that $\widehat{\Lambda}$ solves \eqref{one phase} on $[0,T]$ is standard:~cf., e.g., \cite[proof of Proposition 6.1]{BaSh} and observe that the crossing property of Brownian motion (\cite[Chapter 2, Theorem 9.23(ii)]{KaSh}) is inherited by $R_{\cdot\wedge\tau_{r_0}}$.

\medskip

In comparison, \eqref{eq.sec3.sergey.growth.cond.t0} and the fact that $f(x)=0$ for $x\in(\Lambda_0,\Lambda_{0-})$ imply
\begin{align*}
&\frac{(\widehat\Lambda_{0-})^3}{3}-\frac{(\Lambda_s)^3}{3} \geq \int_{\Lambda_T\wedge\widetilde{\Lambda}_T-\varepsilon_2}^{\widehat\Lambda_{0-}} f(x)\,\PP^x\big(\tau^-_{\Lambda_{0+\cdot}}\le s\wedge\tau_{r_0}\big)\,\nu(\mathrm{d}x),\quad s\in[0,T].
\quad
\end{align*}
In particular, $\Lambda\leq \widehat\Lambda^0$ on $[0,T]$. Moreover, if $\Lambda\le \widehat\Lambda^{n-1}$ on $[0,T]$ for some $n\ge1$, then 
\begin{align*}
\frac{(\widehat\Lambda_{0-})^3}{3}-\frac{(\Lambda_s)^3}{3} &\ge \int_{\Lambda_T\wedge\widetilde{\Lambda}_T-\varepsilon_2}^{\widehat\Lambda_{0-}} f(x)\,\PP^x\big(\tau^-_{\Lambda_{0+\cdot}}\le s\wedge\tau_{r_0}\big)\,\nu(\mathrm{d}x) \\
& \ge \int_{\Lambda_T\wedge\widetilde{\Lambda}_T-\varepsilon_2}^{\widehat\Lambda_{0-}} f(x)\,\PP^x\big(\tau^-_{\widehat\Lambda^{n-1}_{0+\cdot}}\le s\wedge\tau_{r_0}\big)\,\nu(\mathrm{d}x)
=\frac{(\widehat\Lambda_{0-})^3}{3}-\frac{(\widehat{\Lambda}^n_s)^3}{3},\quad s\in[0,T],
\end{align*}
and, thus, $\Lambda\leq \widehat\Lambda^n$ on $[0,T]$. Arguing inductively, and taking $n\to\infty$ we deduce $\Lambda\leq \widehat\Lambda$ on $[0,T]$. Similarly, we see that $\widetilde\Lambda\leq \widehat\Lambda$ on $[0,T]$. In addition, $\widehat{\Lambda}_0\ge\Lambda_0=\widehat\Lambda_{0-}$ implies $\widehat\Lambda_0=\widehat\Lambda_{0-}$.

\medskip

\noindent\textbf{Step 3.} Thanks to the results of Steps 1 and 2 it remains to show that, given any small enough $\varepsilon>0$, the inequality
\begin{equation}\label{start of step 3 claim}
\begin{split}
&\EE\big[\sup_{0\le s\le t}
(B_s+\widehat{\Lambda}_s-\widehat{\Lambda}_{0})\big]
\leq \varepsilon
\EE\big[\psi\big(\sup_{0\le s\le t} (B_s+\widehat{\Lambda}_{0}-\widehat{\Lambda}_s)\wedge\varepsilon_2\big)\big]
\end{split}
\end{equation}
holds for all sufficiently small $\varepsilon_2>0$ and for all $t\in(0,t_1(\varepsilon_2))$, with some $t_1(\varepsilon_2)>0$. Indeed, the comparisons established in Step 2 yield
\begin{align*}
& \EE\big[\max_{0\le s\le t}\,
(B_s+\widetilde{\Lambda}_{s} - \Lambda_{0})\big]
\leq \EE\big[\sup_{0\le s\le t}
(B_s+\widehat{\Lambda}_s-\widehat{\Lambda}_{0})\big],\\
& \EE\big[\psi\big(\sup_{0\le s\le t} (B_s+\widehat{\Lambda}_{0}-\widehat{\Lambda}_s)\wedge\varepsilon_2\big)\big]
\leq \EE\big[\psi\big(\max_{0\le s\le t}\, (B_s+\Lambda_{0}-\Lambda_{s})\wedge\varepsilon_2\big)\big],
\end{align*}
so that the result of Step 1 and \eqref{start of step 3 claim} then complete the proof of the proposition. We show \eqref{start of step 3 claim} in Step 4, for which in this step we prepare the following lemma.

\begin{lemma}\label{le:ComparisonPsiandGrowthWidehat}
Given any small enough $\varepsilon_2,\eta>0$, the following holds for all sufficiently small $t>0$:
\begin{align*}
\EE[\Psi(\sup_{0\le s\le t} (B_s+\widehat\Lambda_{0}-\widehat\Lambda_s)\wedge\varepsilon_2)]
\ge\eta\EE\big[\sup_{0\le s\le t} (B_s+\widehat\Lambda_t-\widehat\Lambda_{t-s})\big],
\end{align*}
where $\Psi(x):=\int_0^x \psi(r)\,\mathrm{d}r$, $x\ge0$, and $\psi$ is extended by $1$ to $[0,\infty)$.
\end{lemma}

\noindent\textbf{Proof.} The inequalities $\Psi(x)\le\Psi(x\wedge\varepsilon_2)+\int_0^{(x-\varepsilon_2)_+} \psi(y+\varepsilon_2)\,\mathrm{d}y\le \Psi(x\wedge\varepsilon_2)+(x-\varepsilon_2)_+$ yield
\begin{equation}\label{Psi term appearance}
\begin{split}
&\;\EE[\Psi(\sup_{0\le s\le t} (B_s+\widehat{\Lambda}_{0}-\widehat{\Lambda}_s)\wedge\varepsilon_2)] \\
&\ge \EE[\Psi(\sup_{0\le s\le t} (B_s+\widehat{\Lambda}_{0}-\widehat{\Lambda}_s))]
-\E[(\sup_{0\le s\le t} (B_s+\widehat{\Lambda}_{0}-\widehat{\Lambda}_s)-\varepsilon_2)_+].
\end{split}
\end{equation}
Moreover, given any $\eta'>0$ we can bound the last expectation according to
\begin{equation}\label{gymnastics2}
\begin{split}
&\;\E\big[\big(\sup_{0\le s\le t}\, (B_s+\widehat{\Lambda}_{0}-\widehat{\Lambda}_s)-\varepsilon_2\big)_+\big] \\
&\le \E\big[\big(\max_{0\le s\le t} B_s\big)
\,\mathbf{1}_{\{\max_{0\le s\le t} B_s>\varepsilon_2/2\}}\big] \\
&\le \eta'\,\frac{\varepsilon_2}{8}\,\PP\big(\max_{0\le s\le t} B_s\in[\varepsilon_2/4,\varepsilon_2/2]\big)\,
\PP\big(\max_{0\le s\le t} B_s 
	\ge \max_{0\le s\le t} (-B_s)\,\big|\,\max_{0\le s\le t} B_s\in[\varepsilon_2/4,\varepsilon_2/2]\big) \\
&= \eta'\,\frac{\varepsilon_2}{8}\,\PP\big(\max_{0\le s\le t} |B_s|\in[\varepsilon_2/4,\varepsilon_2/2],\,\max_{0\le s\le t} B_s 
\ge \max_{0\le s\le t} (-B_s)\big) \\
&\le \eta'\,\E\big[\mathbf{1}_{\{\max_{0\le s\le t} |B_s|\in[\varepsilon_2/4,\varepsilon_2/2],\,\max_{0\le s\le t} B_s 
	\ge \max_{0\le s\le t} (-B_s)\}}\,
\sup_{0\le s\le t} (B_s+\widehat{\Lambda}_{t}-\widehat{\Lambda}_{t-s})\big] \\
&\le\eta'\,\EE\big[\mathbf{1}_{\{\max_{0\le s \le t} |B_s|\le\widehat{\Lambda}_{t}/4\}}\,\sup_{0\le s\le t} (B_s+\widehat{\Lambda}_{t}-\widehat{\Lambda}_{t-s})\big],
\end{split}
\end{equation}
for all small enough $T>0$ and all $t\in(0,T]$.

\medskip

To estimate the first expectation on the right-hand side of \eqref{Psi term appearance} we assume that $T,\varepsilon_2>0$ are small enough, such that \eqref{psi bound} holds, and observe that the growth condition in \eqref{one phase} implies
\begin{equation}\label{eq.sec2.prop2.2.step3.main1}
\begin{split}
\frac{\widehat{\Lambda}_{0}^3}{3}-\frac{\widehat{\Lambda}_{t}^3}{3} 
& = \int_{r_0}^{\widehat{\Lambda}_{0}} (1-\psi(\widehat{\Lambda}_{0}-x))\,\PP^x\big(\sup_{0\le s\leq t\wedge\tau_{r_0}} (R_s-\widehat{\Lambda}_s)\ge0\big)\,\nu(\mathrm{d}x).
\end{split}
\end{equation}
Now, consider the processes $\{(\overline{R}^x_s)_{0\le s\le t}\}_{x>0}$ defined as the (unique) strong solutions of the SDEs 
\begin{equation*}
\overline{R}^x_s = x - \int_0^{s\wedge\overline{\tau}^x_0} \frac{1}{\overline{R}^x_{s'}}\,\mathrm{d}s' 
+ B_{t-s\wedge\overline{\tau}^x_0}-B_{t},\quad s\in[0,t],
\end{equation*}
where $\overline{\tau}^x_0:=\inf\{s\geq0\!:\overline{R}^x_s\leq0\}$. Using strong uniqueness it is easy to see that
\begin{align}
& \overline{R}^{R^x_{t}}_s = R^x_{t-s},\;\; s\in[0,t], \label{eq.sec2.prop2.2.Rbar.prop.0} \\
& R^{\overline{R}^x_{t}}_s = \overline{R}^x_{t-s},\;\; s\in[0,t]\quad\text{if}\quad \overline{R}^x_{t}>0.\label{eq.sec2.prop2.2.Rbar.prop.1}
\end{align}
Moreover, It\^o's formula and the integrability of $(\overline{R}^x)^4$ reveal  
that $(\overline{R}^x)^3$ is a martingale in its natural filtration. In particular, $\EE[(\overline{R}^x_t)^3]=x^3$.

\medskip

In view of the continuity of $x\mapsto\overline{R}^x_{t}$ on $\{x>0\!:\overline{R}^x_t>0\}$ and \eqref{eq.sec2.prop2.2.Rbar.prop.1}, we have that
\begin{align*}
&\;\big\{x\in[r_0,\widehat\Lambda_{0}]:\,\min_{0\le s\leq t} R^x_s>r_0,\,
\sup_{0\le s\leq t} (R^x_s-\widehat\Lambda_s)\ge0\big\} \\
&= \big\{\overline{R}^x_t\in[r_0,\widehat\Lambda_{0}]:\,
\min_{0\le s\leq t} R^{\overline{R}^x_{t}}_s>r_0,
\,\sup_{0\le s\leq t} (R^{\overline{R}^x_{t}}_s-\widehat\Lambda_s)\geq0\big\} \\
& = \big\{\overline{R}^x_t\in[r_0,\widehat\Lambda_{0}]:\,
\min_{0\le s\le t} \overline{R}^x_{s}>r_0,
\,\sup_{0\le s\leq t} (\overline{R}^x_{s}-\widehat\Lambda_{t-s})\geq0\big\}.
\end{align*}
Let $\underline{X}\in[r_0,\widehat{\Lambda}_0]$ be the infimum of the above (random) interval.~It holds $\overline{R}^x_s\leq \overline{R}^{y}_s + x - y$, $s\in[0,\overline\tau^x_0]$, $0<x\leq y$ due to the monotonicity of the drift of $\overline{R}$. We also have $\overline{R}^x_s\geq x+B_{t-s}-B_{t} - s/r_0$, $0\le s\leq\overline{\tau}^x_{r_0}:=\inf\{s\geq0\!:\overline{R}^x_s\leq r_0\}$, $x>0$.
Using these observations and introducing 
\begin{align*}
&\widecheck{B}_s=B_{t-s}-B_{t},
\quad \xi= \sup_{0\le s\le t}\, (\widehat\Lambda_t+\widecheck{B}_{s}-s/r_0-\widehat\Lambda_{t-s})
\end{align*}
we deduce that, on the event $\{\inf_{0\le s\le t} \,(\widehat\Lambda_{t}+\widecheck{B}_{s}-s/r_0-r_0)>\xi\}$, we have
\begin{align*}
& \underline{X} \leq \overline{R}^{\widehat\Lambda_{t}}_{t} - \xi.
\end{align*}

\smallskip

Returning to the main line of the argument we lower bound
\begin{equation}\label{eq.sec2.prop2.2.step3.main2}
\begin{split}
\int_{r_0}^{\widehat\Lambda_0} \PP^x\big(\sup_{0\le s\leq t\wedge\tau_{r_0}} (R_s-\widehat{\Lambda}_s)\ge0\big)\,\nu(\mathrm{d}x)\ge
\EE\bigg[\int_{r_0}^{\widehat{\Lambda}_{0}} \mathbf{1}_{\{\min_{0\le s\le t} (R^x_s-r_0)>0,\,\sup_{0\le s\le t} (R^x_s-\widehat{\Lambda}_s)\ge0\}}\,x^2\,\mathrm{d}x\bigg] \\
= \frac{1}{3}\EE\big[(\widehat{\Lambda}_0)^3 - \underline{X}^3\big] \\
\ge \frac{1}{3}\EE\big[\big((\widehat{\Lambda}_0)^3 - (\overline{R}^{\widehat\Lambda_{t}}_{t}-\xi)^3\big)\,\bone_{\{\inf_{0\le s\le t} (\widehat\Lambda_t+\widecheck{B}_{s}-s/r_0-r_0)>\xi\}}\big] \\
= \frac{1}{3} \EE\big[(\widehat{\Lambda}_0)^3 - (\overline{R}^{\widehat\Lambda_{t}}_{t}-\xi)^3\big]
- \frac{1}{3}\EE\big[\big((\widehat{\Lambda}_0)^3 - (\overline{R}^{\widehat\Lambda_{t}}_{t}-\xi)^3\big)\,\bone_{\{\inf_{0\le s\le t} (\widehat\Lambda_t+\widecheck{B}_{s}-s/r_0-r_0)\le\xi\}}\big] \\
\geq \frac{1}{3}\big((\widehat{\Lambda}_0)^3 - (\widehat{\Lambda}_t)^3 \big) + \EE\big[\xi\big((\overline{R}^{\widehat\Lambda_{t}}_{t})^2-\xi\overline{R}^{\widehat\Lambda_{t}}_{t}+\xi^2/3\big)\big]
-\frac{1}{3}\EE\big[\big((\widehat{\Lambda}_0)^3+\xi^3\big)\,\bone_{\{\inf_{0\le s\le t} (\widehat\Lambda_t+\widecheck{B}_{s}-s/r_0-r_0)\le\xi\}}\big].
\end{split}
\end{equation}
Next, we recall that the Bessel process $(R^x_s)_{0\le s\le t}$ can be constructed by taking a Brownian motion $(x+B_{s\wedge\tau^x_0})_{0\le s\le t}$ stopped at the first hitting time of zero, denoted by $\tau^x_0$, and changing the measure via the Radon-Nikodym density $x^{-1}(x+B_{t\wedge\tau^x_0})$. Thus,
\begin{equation*}
\begin{split}
 &\;\int_{r_0}^{\widehat{\Lambda}_{0}} \psi(\widehat{\Lambda}_{0}-x)\,\PP^x\big(\sup_{0\le s\leq t\wedge\tau_{r_0}} (R_s-\widehat{\Lambda}_s)\ge0\big)\,\nu(\mathrm{d}x) \\
 & = \int_{r_0}^{\widehat{\Lambda}_{0}} \psi(\widehat{\Lambda}_{0}-x)\,x\,
 \EE\big[(x+B_{t\wedge\tau^x_{0}\wedge\tau^x_{\widehat{\Lambda}}})\,\mathbf{1}_{\{x+\sup_{0\le s\le t} (B_{s\wedge\tau^x_{r_0}}-\widehat{\Lambda}_s)\ge0\}}\big]\,\mathrm{d}x \\
 & \le(\widehat{\Lambda}_0)^2\,\EE\bigg[\int_{r_0}^{\widehat{\Lambda}_{0}} \psi(\widehat{\Lambda}_{0}-x)\,\mathbf{1}_{\{x+\sup_{0\le s\le t} (B_s-\widehat{\Lambda}_s)\ge0\}}\,\mathrm{d}x\bigg] 
 \le (\widehat{\Lambda}_0)^2\,\EE\big[\Psi\big(\widehat{\Lambda}_{0}+\sup_{0\le s\le t} (B_s-\widehat{\Lambda}_s)\big)\big],
\end{split}
\end{equation*}
where $\tau^x_{\widehat\Lambda}$ is the first crossing time of $\widehat\Lambda$ by $x+B$. Combining this with \eqref{eq.sec2.prop2.2.step3.main1}, \eqref{eq.sec2.prop2.2.step3.main2} we arrive at 
\begin{align}
&\;\EE\big[\Psi\big(\sup_{0\le s\le t} (B_s+\widehat{\Lambda}_{0}-\widehat{\Lambda}_s)\big)\big]
\nonumber \\
&\ge\frac{1}{(\widehat{\Lambda}_0)^2}\,\EE\big[\xi\big((\overline{R}^{\widehat\Lambda_{t}}_{t})^2-\xi\overline{R}^{\widehat\Lambda_{t}}_{t}+\xi^2/3\big)\big]
-\frac{1}{3(\widehat\Lambda_0)^2}\EE\big[\big((\widehat{\Lambda}_0)^3+\xi^3\big)\,\bone_{\{\inf_{0\le s\le t} (\widehat\Lambda_t+\widecheck{B}_{s}-s/r_0-r_0)\le\xi\}}\big]
\nonumber \\
&\ge \frac{1}{4(\widehat{\Lambda}_0)^2}\,\EE\big[\xi(\overline{R}^{\widehat\Lambda_t}_t)^2\big]
-\frac{1}{3(\widehat\Lambda_0)^2}\EE\big[\big((\widehat{\Lambda}_0)^3+\xi^3\big)\,\bone_{\{\inf_{0\le s\le t} (\widehat\Lambda_t+\widecheck{B}_{s}-s/r_0-r_0)\le\xi\}}\big]. \label{gymnastics3}
\end{align}

\smallskip

It remains to estimate the terms in \eqref{gymnastics3}. To bound the first term from below, we define the events $\overline{A}_{\overline{\varepsilon}}=\{\max_{0\le s\le t} |\widecheck{B}_{s}|\le\overline{\varepsilon}\}$, for $\overline{\varepsilon}>0$. Setting $\overline{\varepsilon}=\widehat{\Lambda}_t/2$ we note that on $\overline{A}_{\overline{\varepsilon}}$ it holds 
\begin{equation}\label{R LBD}	
\overline{R}^{\widehat\Lambda_{t}}_s\ge\widehat{\Lambda}_{t}-\frac{4t}{\widehat{\Lambda}_{t}}
-\frac{\widehat{\Lambda}_{t}}{2}
\end{equation}
for all $s\in[0,t]$ preceding the first time $\overline{R}^{\widehat\Lambda_{t}}$ hits $\widehat{\Lambda}_{t}/4$. But for all small enough $t>0$, the right-hand side of \eqref{R LBD} is strictly greater than $\widehat{\Lambda}_{t}/4$, and so $\overline{R}^{\widehat\Lambda_{t}}_s>\widehat{\Lambda}_{t}/4$, $s\in[0,t]$ on $\overline{A}_{\overline{\varepsilon}}$. Consequently, 
\begin{equation}\label{gymnastics4}
\begin{split}
\frac{1}{4(\widehat{\Lambda}_0)^2}\,\EE\big[\xi(\overline{R}^{\widehat\Lambda_{t}}_{t})^2\big]
&\ge \frac{(\widehat{\Lambda}_t)^2}{64(\widehat{\Lambda}_0)^2}
\,\EE\big[\mathbf{1}_{\overline{A}_{\overline{\varepsilon}}}\,\sup_{0\le s\le t}\, (\widehat\Lambda_t+\widecheck{B}_{s}-s/r_0-\widehat\Lambda_{t-s})\big].
\end{split}
\end{equation}
Next, we rewrite the latter expression by means of Girsanov's Theorem: 
\begin{equation}\label{gymnastics5}
\begin{split}
&\;\frac{(\widehat{\Lambda}_t)^2}{64(\widehat{\Lambda}_0)^2}
\,\EE\Big[e^{-\frac{\overline{B}_{t}}{r_0}
-\frac{t}{2r_0^2}}
\,\mathbf{1}_{\overline{A}_{\overline{\varepsilon}}}\,\sup_{0\le s\le t}\, (\widehat\Lambda_t+\overline{B}_{s}-\widehat\Lambda_{t-s})\Big] \\
&\ge \frac{(\widehat{\Lambda}_t)^2}{64(\widehat{\Lambda}_0)^2}\,
e^{-\frac{\overline{\varepsilon}}{r_0}}	\,
\EE\big[\mathbf{1}_{\{\max_{0\le s \le t} |\overline{B}_s|\le \overline{\varepsilon}-t/r_0\}}
\,\sup_{0\le s\le t}\,(\widehat{\Lambda}_{t}+\overline{B}_s-\widehat{\Lambda}_{t-s})\big] \\
&\ge \frac{(\widehat{\Lambda}_t)^2
\,e^{-\frac{\widehat\Lambda_t}{2r_0}}}{64(\widehat{\Lambda}_0)^2}\,
\EE\big[\mathbf{1}_{\{\max_{0\le s \le t} |\overline{B}_s|\le\widehat{\Lambda}_t/4\}}
\,\sup_{0\le s\le t}\,(\widehat{\Lambda}_{t} 
+\overline{B}_s-\widehat{\Lambda}_{t-s}) \big],
\end{split}
\end{equation}
where $\overline{B}_s:=\widecheck{B}_s-s/r_0$, $s\in[0,t]$. Further, by repeating the estimates in \eqref{gymnastics2}, we infer that
\begin{equation}\label{indicator removal}
\EE\big[\mathbf{1}_{\{\max_{0\le s \le t} |\overline{B}_s|\le\widehat{\Lambda}_t/4\}}
\,\sup_{0\le s\le t}\,(\widehat{\Lambda}_{t} 
+\overline{B}_s-\widehat{\Lambda}_{t-s}) \big]
\ge\frac{1}{2}\,\EE\big[\sup_{0\le s\le t}\,(\widehat{\Lambda}_{t}+B_s
-\widehat{\Lambda}_{t-s})\big]
\end{equation}
for all small enough $t>0$.
Combining \eqref{gymnastics4}--\eqref{indicator removal} we obtain
\begin{align}
& \frac{1}{4(\widehat{\Lambda}_0)^2}\,\EE[\xi(\overline{R}^{\widehat\Lambda_{t}}_{t})^2]
\geq \frac{(\widehat{\Lambda}_t)^2
\,e^{-\frac{\widehat\Lambda_t}{2r_0}}}{128(\widehat{\Lambda}_0)^2}\,\EE\big[\sup_{0\le s\le t}\,(B_s+\widehat{\Lambda}_{t} -\widehat{\Lambda}_{t-s})\big].\label{eq.sec2.prop2.2.step3.xiR2}
\end{align}

To bound the second term in \eqref{gymnastics3}, we argue as in \eqref{gymnastics2} that for any given $r\in(0,\widehat\Lambda_t-r_0)$,
\begin{align*}
&\,\frac{1}{3(\widehat\Lambda_0)^2}\EE\big[\big((\widehat{\Lambda}_0)^3+\xi^3\big)\,\bone_{\{\inf_{0\le s\le t} (\widehat\Lambda_t+\widecheck{B}_{s}-s/r_0-r_0)\le\xi\}}\big] \\
&\leq\frac{1}{3(\widehat\Lambda_0)^2}\EE\big[\big((\widehat{\Lambda}_0)^3+\xi^3
\big)\,\bone_{\{\xi\ge r/2\}}\big]  +\frac{1}{3(\widehat\Lambda_0)^2}\EE\big[\big((\widehat{\Lambda}_0)^3+(r/2)^3\big)
\,\bone_{\{\inf_{0\le s\le t} (\widecheck{B}_{s}-s/r_0)\le-r/2\}}\big] \\
& \leq \PP(\xi\ge r/4)
\leq \frac{4}{r} \EE[\xi\,\bone_{\{\xi\ge r/4\}}],
\end{align*}
provided $t>0$ is small enough.~Arguing again as in \eqref{gymnastics2}, we control the latter by
\begin{align*}
&\frac{(\widehat{\Lambda}_t)^2
	\,e^{-\frac{\widehat\Lambda_t}{2r_0}}}{256(\widehat{\Lambda}_0)^2}\,\EE\big[\sup_{0\le s\le t}\,(B_s+\widehat{\Lambda}_{t} -\widehat{\Lambda}_{t-s})\big]
\end{align*}
for all small enough $t>0$. Putting this together with \eqref{eq.sec2.prop2.2.step3.xiR2} and \eqref{gymnastics3}, we conclude that
\begin{align*}
&\EE\big[\Psi\big(\sup_{0\le s\le t} (B_s+\widehat{\Lambda}_{0}-\widehat{\Lambda}_s)\big)\big]
\ge \frac{(\widehat{\Lambda}_t)^2
\,e^{-\frac{\widehat\Lambda_t}{2r_0}}}{256(\widehat{\Lambda}_0)^2}\,\EE\big[\sup_{0\le s\le t}\,(B_s+\widehat{\Lambda}_{t} -\widehat{\Lambda}_{t-s})\big].
\end{align*}
In view of \eqref{Psi term appearance} and \eqref{gymnastics2}, the lemma readily follows. \qed

\medskip

\noindent\textbf{Step 4.} Lemma \ref{le:ComparisonPsiandGrowthWidehat} and the computation
\begin{equation*}
\EE\big[\sup_{0\le s\le t} (B_s+\widehat\Lambda_t-\widehat\Lambda_{t-s})\big]
=\EE\big[\sup_{0\le s\le t} (B_{t-s}-B_t+\widehat\Lambda_t-\widehat\Lambda_s)\big]
=\EE\big[\sup_{0\le s\le t} (B_s+\widehat\Lambda_t-\widehat\Lambda_s)\big]
\end{equation*}
result in the next corollary, which allows us to compare $\widehat\Lambda_0-\widehat\Lambda_t$ with the solution to a one-dimensional one-phase supercooled Stefan problem.


\begin{corollary}
Given any small enough $\varepsilon_2,\eta>0$, one has for all sufficiently small $t>0$: 
\begin{align*}
 \widehat\Lambda_{0}-\widehat\Lambda_t&\ge\EE\big[\sup_{0\le s\le t} (B_s+\widehat\Lambda_{0}-\widehat\Lambda_s)\big]-\eta^{-1}\EE\big[\Psi\big(\sup_{0\le s\le t} (B_s+\widehat\Lambda_{0}-\widehat\Lambda_s)\wedge\varepsilon_2\big)\big] \\
&\ge\EE\bigg[\int_0^{\sup_{0\le s\le t}\,(B_s+\widehat\Lambda_{0}-\widehat\Lambda_s)} \underline{f}(x)\,\mathrm{d}x\bigg],   
\end{align*}
where 
\begin{align*}
\underline{f}(x):=
\begin{cases}
1-\eta^{-1}\psi(x),&\quad x\in[0,\varepsilon_2], \\
0,&\quad x\in(\varepsilon_2,\infty)
\end{cases}
\end{align*}
is a sub-probability density.
\end{corollary}

Let us fix arbitrary $\varepsilon_2,\eta,t>0$ as in the above corollary, and let $\lambda\!:[0,\infty)\to[0,\infty)$ be a probabilistic solution to the one-dimensional one-phase supercooled Stefan problem with initial density $\underline{f}$: i.e., $\lambda_0=0$, $\lambda$ is non-decreasing and right-continuous, and 
\begin{align*}
\lambda_s &=\int_0^\infty \underline{f}(x)\,\EE[\bone_{\{\inf_{0\le s'\le s}\, (x-B_{s'}-\lambda_{s'})\le0\}}]\,\mathrm{d}x \\
&=\EE\bigg[\int_0^{\sup_{0\le s'\le s}\,(B_{s'}+\lambda_{s'})}
\underline{f}(x)\,\mathrm{d}x\bigg] \\
&=\EE\big[\sup_{0\le s'\le s}\,(B_{s'}+\lambda_{s'})\big]
-\EE\big[\eta^{-1}\Psi\big(\sup_{0\le s'\le s}\, (B_{s'}+\lambda_{s'})\wedge\varepsilon_2\big)
+\big(\sup_{0\le s'\le s}\,(B_{s'}+\lambda_{s'})-\varepsilon_2\big)_+\big].
\end{align*}
Using an iteration argument as in Step 2 we construct a function $\lambda$ that satisfies the above and deduce easily (by the comparison between $\widehat\Lambda_0-\widehat\Lambda$ and the approximations of $\lambda$) that $\widehat\Lambda_0-\widehat\Lambda_s\ge\lambda_s$ for $s\in[0,t]$. In addition, the following sub-additivity property holds. 
\begin{proposition}\label{prop sub}
For all $s\in[0,t]$, it holds $\lambda_s+\lambda_{t-s}\ge\lambda_t$.    
\end{proposition}

\noindent\textbf{Proof of Proposition \ref{prop sub}.} We will show a stronger statement.~For any $q\in(0,1]$, the scaling property of Brownian motion reveals the rescaled function $\lambda_s^{(q)}:=\frac{\lambda_{qs}}{\sqrt{q}}$ as a solution of the problem
\begin{align*}
\lambda_s^{(q)}=\int_0^\infty \underline{f}(\sqrt{q}x)\,\EE[\bone_{\{\inf_{0\le s'\le s}\, (x-B_{s'}-\lambda^{(q)}_{s'})\le0\}}]\,\mathrm{d}x
=\EE\bigg[\int_0^{\sup_{0\le s'\le s}\,(B_s+\lambda_s^{(q)})} \underline{f}(\sqrt{q}x)\,\mathrm{d}x\bigg].
\end{align*}
Since $\underline{f}$ is non-increasing, we see, by considering the iterates defining $\lambda$, as well as the analogous iterates for $\lambda^{(q)}$, that $\lambda^{(q)}_s\ge\lambda_s$, $s\in[0,t]$, i.e., $\lambda_{qs}\ge\sqrt{q}\lambda_s\ge q\lambda_s$, $s\in[0,t]$. Thus, for all $s\in[0,t]$,
\begin{align*}
\lambda_s\ge\frac{s}{t}\lambda_t\quad\text{and}\quad\lambda_{t-s}\ge\frac{t-s}{t}\lambda_t.     
\end{align*}
Adding these two inequalities yields Proposition \ref{prop sub}. \qed

\medskip

Lastly, by repeating \eqref{gymnastics2} we conclude that, for any $\varepsilon_2,\eta'>0$, the following holds for all small enough $t>0$:
\begin{align}
\EE\big[\big(\sup_{0\le s\le t}\,(B_s+\lambda_s)-\varepsilon_2\big)_+\big]
\le\eta'\EE\big[\sup_{0\le s\le t}\,(B_s+\lambda_{t-s}-\lambda_t)\big]. \label{eq.le:Comparisonlambda1}   
\end{align}
On the other hand, for any small enough $\varepsilon_2,\eta>0$, the following holds for all small enough $t>0$: 
\begin{align*}
\EE\big[\sup_{0\le s\le t}\,(B_s+\widehat\Lambda_s-\widehat\Lambda_{0})\big]
\le\EE\big[\sup_{0\le s\le t}\,(B_s-\lambda_s)\big]
\le\EE\big[\sup_{0\le s\le t}\,(B_s+\lambda_{t-s}-\lambda_t)\big],    
\end{align*}
which can be upper bounded further as follows, using the definition of $\lambda$ and \eqref{eq.le:Comparisonlambda1}:
\begin{align*}
\EE\big[\sup_{0\le s\le t}\,(B_s+\lambda_{t-s}-\lambda_t)\big]
&= \EE\big[\sup_{0\le s\le t}\,(B_s-B_t+\lambda_{t-s}-\lambda_t)\big] \\
&=\EE\big[\sup_{0\le s\le t}\,(B_s+\lambda_s)\big]-\lambda_t \\
&=\EE\big[\eta^{-1}\Psi\big(\sup_{0\le s\le t}\,(B_s+\lambda_s)\wedge\varepsilon_2\big)
+\big(\sup_{0\le s\le t}\,(B_s+\lambda_s)-\varepsilon_2\big)_+\big] \\
&\le\eta^{-1}\varepsilon_2\EE\big[\psi\big(\sup_{0\le s\le t}\,(B_s+\lambda_s)\wedge\varepsilon_2\big)\big]
+\frac12\EE\big[\sup_{0\le s\le t}\,(B_s+\lambda_{t-s}-\lambda_t)\big] \\
&=\eta^{-1}\varepsilon_2\EE\big[\psi\big(\sup_{0\le s\le t}\,(B_s+\lambda_s)\wedge\varepsilon_2\big)\big]
+\frac12\big(\EE\big[\sup_{0\le s\le t}\,(B_s+\lambda_s)\big]-\lambda_t\big).
\end{align*}
Thus, for any $\varepsilon_2,\eta>0$, the following holds for all small enough $t>0$:
\begin{align*}
\EE\big[\sup_{0\le s\le t}\,(B_s+\widehat\Lambda_s-\widehat\Lambda_{0})\big]
\le\EE\big[\sup_{0\le s\le t}\,(B_s+\lambda_s)\big]-\lambda_t
&\le 2\eta^{-1}\varepsilon_2\EE\big[\psi\big(\sup_{0\le s\le t}\,(B_s+\lambda_s)\wedge\varepsilon_2\big)\big] \\
&\le 2\eta^{-1}\varepsilon_2\EE\big[\psi\big(\sup_{0\le s\le t}\,(B_s+\widehat\Lambda_{0}-\widehat\Lambda_s)\wedge\varepsilon_2\big)\big].
\end{align*}
This verifies \eqref{start of step 3 claim}, with $\varepsilon=2\eta^{-1}\varepsilon_2$, and completes the proof of Proposition \ref{prop 2nd phase}. \qed

\section{Boundary and Interior Regularity} \label{se:3}

This section is devoted to the regularity theory for \eqref{intro: growth weak}, \eqref{intro: phys rig}, akin to the one in \cite[Sections 2, 3]{dns} for the one phase supercooled Stefan problem. Hereby, the proofs which are easy adaptations of those in \cite{dns} are relegated to the appendix. 

\subsection{Hölder Continuity} \label{sec:Holder}


\begin{proposition}
Fix an arbitrary $t_0\in[0,\zeta)$ and assume that $w(t_0-,\cdot)$ satisfies at least one of the following two conditions: 
\begin{enumerate}[(a)]
\item $\lim_{y_0\downarrow0}\,\esssupmath{0<y<y_0}\; w(t_0-,\Lambda_{t_0-}-y)<1$, or
\item $w(t_0-,\Lambda_{t_0-}-\cdot)$ is monotone in a right neighborhood of any point in $[0,\Lambda_{t_0-}-r_0)$.
\end{enumerate}
Then, there exist $\upsilon_0,y_0>0$ and $c\!:(0,\upsilon_0)\to[0,1)$ such that, for any $\upsilon\in(0,\upsilon_0)$,
\begin{equation*}
\int_0^y w(s-,\Lambda_{s-}-y')(\Lambda_{s-}-y')^2\,\mathrm{d}y'
\le c(\upsilon)\int_0^y (\Lambda_{s-}-y')^2\,\mathrm{d}y',\quad y\in[0,y_0],\quad s\in[t_0+\upsilon,t_0+\upsilon_0].
\end{equation*}
In case (a), there is an extension $c\!:[0,\upsilon_0)\to[0,1)$ and the above bound is also true for $\upsilon=0$.
\end{proposition}

\noindent\textbf{Proof.} In case (a), $\Lambda_{t_0}=\Lambda_{t_0-}$ and $w(t_0,x)=w(t_0-,x)$, $x\in(0,\infty)\backslash\{\Lambda_{t_0}\}$. Moreover, there exist $y_0\in(0,\Lambda_{t_0}/2)$ and $c\in[0,1)$ such that $w(t_0,\Lambda_{t_0}-y)\leq c$ for all $y\in[0,y_0]$. In particular, we can find a non-decreasing function $\psi\!:[0,y_0]\to[0,\infty)$ with $\psi(0)=\frac{1-c}{2}$ such that 
\begin{equation}\label{1-psi}
w(t_0,\Lambda_{t_0}-y)\leq1-\psi(y),\quad y\in[0,y_0].
\end{equation}
In case (b), we consider a version of $w(t_0-,\Lambda_{t_0-}-\cdot)$ that is right-continuous. Since $w(t_0,\Lambda_{t_0}-y)=w(t_0-,\Lambda_{t_0}-y)$, $y\in(0,\Lambda_{t_0})$ (see the paragraph following \eqref{eq.sec3.sergey.growth.cond.t0}), we deduce the right-continuity and local monotonicity of $w(t_0,\Lambda_{t_0}-\cdot)$ on $[0,\Lambda_{t_0})$. Further, the physicality condition \eqref{intro: phys} implies $w(t_0,\Lambda_{t_0}-)\leq1$. If $w(t_0,\Lambda_{t_0}-)<1$, then we are led back to case (a). If $w(t_0,\Lambda_{t_0}-)=1$, then $y\mapsto w(t_0,\Lambda_{t_0}-y)$ is non-increasing and strictly below $1$ in a right neighborhood of $0$, due to the physicality condition \eqref{intro: phys} and the local monotonicity of $w(t_0,\Lambda_{t_0}-\cdot)$. Hence, we can pick a $y_0\in(0,\Lambda_{t_0}/2)$ and a non-decreasing function $\psi\!:[0,y_0]\to[0,\infty)$, with $\psi(0)=0$ and $\psi>0$ on $(0,y_0]$, such that \eqref{1-psi} holds. This allows us to proceed independently of the case we are in. 

\medskip

By the right-continuity of $\Lambda$, we can find an $\upsilon_0>0$ for which $\Lambda_{t_0}-\Lambda_{t_0+\upsilon_0}\leq y_0/8$.~Taking limits in \eqref{eq.sec3.sergey.w.def} (cf.~\cite[proof of Lemma 3.2]{NaShsurface}) and repeating the proof of the first statement in Proposition \ref{prop:MarkovSystem} we obtain
for $s\in(t_0,t_0+\upsilon_0]$ and $y\in(0,\Lambda_{s-}-r_0)$, 
\begin{equation} \label{eq:Prop2.1split1}
w(s-,\Lambda_{s-}-y)
=\EE^{\Lambda_{s-}-y}\big[w(t_0,R_{s-t_0})\,\bone_{\{\tau^{-}_{\Lambda_{s-\cdot}}\wedge\tau_{r_0}>s-t_0\}}\big]
\leq\EE^{\Lambda_{s-}-y}\big[w(t_0,R_{s-t_0})\,\bone_{\{R_{s-t_0}<\Lambda_{t_0}\}}\big].
\end{equation}
We split the last term into three parts:
\begin{equation} \label{eq:Prop2.1split2}
\begin{split}
\EE^{\Lambda_{s-}-y}\big[w(t_0,R_{s-t_0})\,\bone_{\{R_{s-t_0}<\Lambda_{t_0}\}}\big] 
=&\;\EE^{\Lambda_{s-}-y}\big[w(t_0,R_{s-t_0})\,\bone_{\{R_{s-t_0}\in[\Lambda_{t_0}-y_0/4,\Lambda_{t_0})\}}\big] \\
&+\EE^{\Lambda_{s-}-y}\big[w(t_0,R_{s-t_0})\,\bone_{\{R_{s-t_0}\in[\Lambda_{t_0}-y_0,\Lambda_{t_0}-y_0/4)\}}\big] \\
&+\EE^{\Lambda_{s-}-y}\big[w(t_0,R_{s-t_0})\,\bone_{\{R_{s-t_0}\in(0,\Lambda_{t_0}-y_0)\}}\big].
\end{split}
\end{equation}
Since $\psi$ is non-decreasing, the first two summands can be estimated respectively by
\begin{eqnarray}
&& (1-\psi(0))\,\PP^{\Lambda_{s-}-y}\big(R_{s-t_0}\in[\Lambda_{t_0}-y_0/4,\Lambda_{t_0})\big), \label{eq:Prop2.1term1} \\
&&
(1-\psi(y_0/4))\,\PP^{\Lambda_{s-}-y}\big(R_{s-t_0}\in[\Lambda_{t_0}-y_0,\Lambda_{t_0}-y_0/4)\big).
\label{eq:Prop2.1term2}
\end{eqnarray}
Due to $w(t_0,y)\!\le\!\|w(0-,\cdot)\|_\infty$ the third summand is at most $\|w(0-,\cdot)\|_\infty\,\PP^{\Lambda_{s-}-y}(R_{s-t_0}\!\in\!(0,\Lambda_{t_0}\!-y_0))$.

\medskip

The aim now is to ``absorb'' the latter quantity into the ``gap'' provided by $\psi(y_0/4)$ in \eqref{eq:Prop2.1term2}. To this end, we note that for all $s\in(t_0,t_0+\upsilon_0]$ and all $y\in[0,y_0/8]$ (recall $\Lambda_{t_0}-\Lambda_{t_0+\upsilon_0}\le y_0/8$),
\begin{IEEEeqnarray*}{rCl}
\Lambda_{t_0}\ge\Lambda_{s-}-y=\Lambda_{t_0}-(\Lambda_{t_0}-\Lambda_{s-})-y
\ge\Lambda_{t_0}-(\Lambda_{t_0}-\Lambda_{t_0+\upsilon_0})-y
\ge\Lambda_{t_0}-y_0/4,
\end{IEEEeqnarray*}
i.e., $\Lambda_{s-}-y\in[\Lambda_{t_0}-y_0/4,\Lambda_{t_0}]$. Thinking of $R$ as the norm of a three-dimensional Brownian motion and letting $g_{s-t_0}(z):=(2\pi (s-t_0))^{-\frac32}e^{-\frac{|z|^2}{2(s-t_0)}}$ we then find 
\begin{align}
\PP^{\Lambda_{s-}-y}\big(R_{s-t_0}\in[\Lambda_{t_0}-y_0,\Lambda_{t_0}-y_0/4)\big) &=\int_{B(0,\Lambda_{t_0}-y_0/4)\setminus B(0,\Lambda_{t_0}-y_0)} g_{s-t_0}\big((\Lambda_{s-}-y,0,0)-z\big)\,\mathrm{d}z \notag \\
&\ge\int_{B((\Lambda_{t_0}-3y_0/8,0,0),y_0/8)}
g_{s-t_0}\big((\Lambda_{s-}-y,0,0)-z\big)\,\mathrm{d}z \label{eq:Prop2.1Gap1} \\
&\ge\frac{\pi y_0^3}{384}(2\pi(s-t_0))^{-\frac32}e^{-\frac{y_0^2}{8(s-t_0)}}, \notag
\end{align}
as it holds $|(\Lambda_{s-}-y,0,0)-z)|\le y_0/2$ for $y\in[0,y_0/8]$ and $z\in B((\Lambda_{t_0}-3y_0/8,0,0),y_0/8)$. Also,
\begin{equation}
\begin{split}
\PP^{\Lambda_{s-}-y}\big(R_{s-t_0}\in(0,\Lambda_{t_0}-y_0)\big)
&=\int_{B(0,\Lambda_{t_0}-y_0)} g_{s-t_0}\big((\Lambda_{s-}-y,0,0)-z\big)\,\mathrm{d}z \\
&\leq\frac{4\pi}{3}(\Lambda_{t_0}-y_0)^3(2\pi(s-t_0))^{-\frac32}e^{-\frac{9y_0^2}{32(s-t_0)}},
\end{split}
\end{equation}
since $|(\Lambda_{s-}-y,0,0)-z)|\ge 3y_0/4$ for $y\in[0,y_0/8]$ and $z\in B(0,\Lambda_{t_0}-y_0)$. Thanks to the fast decay of the exponential function, by shrinking $\upsilon_0>0$ if necessary, we can ensure for $s\in(t_0,t_0+\upsilon_0]$:
\begin{equation}\label{eq:Prop2.1term3}
\|w(0-,\cdot)\|_\infty\frac{4\pi}{3}(\Lambda_{t_0}-y_0)^3(2\pi(s-t_0))^{-\frac32}e^{-\frac{9y_0^2}{32(s-t_0)}}\leq\frac12\psi\Big(\frac{y_0}{4}\Big)\frac{\pi y_0^3}{384}(2\pi(s-t_0))^{-\frac32}e^{-\frac{y_0^2}{8(s-t_0)}}.
\end{equation}

\smallskip

Furthermore, 
\begin{align}
\PP^{\Lambda_{s-}-y}\big(R_{s-t_0}\in[\Lambda_{t_0}-y_0/4,\Lambda_{t_0})\big)
= \int_{B(0,\Lambda_{t_0})\setminus B(0,\Lambda_{t_0}-y_0/4)} g_{s-t_0}\big((\Lambda_{s-}-y,0,0)-z)\,\mathrm{d}z \notag \\
\ge \int_{(B(0,\Lambda_{t_0})\setminus B(0,\Lambda_{t_0}-y_0/4))\cap B((\Lambda_{s-}-y,0,0),y_0/4)} g_{s-t_0}\big((\Lambda_{s-}-y,0,0)-z\big)\,\mathrm{d}z \notag \\
= \int_0^{y_0/4} \int_{(B(0,\Lambda_{t_0})\setminus B(0,\Lambda_{t_0}-y_0/4))
\cap\partial B((\Lambda_{s-}-y,0,0),r)} g_{s-t_0}\big((\Lambda_{s-}-y,0,0)-z\big)\,\sigma(\mathrm{d}z)\,\mathrm{d}r \label{eq:Prop2.1Gap3} \\
= \int_0^{y_0/4} \int_{(B(0,\Lambda_{t_0})\setminus B(0,\Lambda_{t_0}-y_0/4))
	\cap\partial B((\Lambda_{s-}-y,0,0),r)} g_{s-t_0}((r,0,0))\,\sigma(\mathrm{d}z)\,\mathrm{d}r 
\notag \\
\ge c_{\Lambda_{t_0},y_0}\int_0^{y_0/4} \int_{\partial B((\Lambda_{s-}-y,0,0),r)}
g_{s-t_0}((r,0,0))\,\sigma(\mathrm{d}z)\,\mathrm{d}r
=c_{\Lambda_{t_0},y_0}\int_{B(0,y_0/4)} g_{s-t_0}(z)\,\mathrm{d}z
\ge\frac{c_{\Lambda_{t_0},y_0}}{2}, \notag
\end{align}
where $c_{\Lambda_{t_0},y_0}\in(0,1)$ is a constant guaranteeing that
\begin{IEEEeqnarray*}{rCl}
\sigma\big((B(0,\Lambda_{t_0})\setminus B(0,\Lambda_{t_0}-y_0/4))\cap \partial B((x,0,0),r)\big) \ge c_{\Lambda_{t_0},y_0}\sigma(\partial B((x,0,0),r)),
\end{IEEEeqnarray*}
for all $x\!\in\![\Lambda_{t_0}\!-y_0/4,\Lambda_{t_0}]$, $r\!\in\![0,y_0/4]$,
and we decrease $\upsilon_0>0$ if necessary, so that for $s\in(t_0,t_0\!+\upsilon_0]$:
\begin{IEEEeqnarray*}{rCl}
\int_{B(0,y_0/4)} g_{s-t_0}(z)\,\mathrm{d}z\ge\frac12.
\end{IEEEeqnarray*}

\smallskip

Plugging the bounds \eqref{eq:Prop2.1term1}--\eqref{eq:Prop2.1Gap3} into  \eqref{eq:Prop2.1split2} we get for all $y\in[0,y_0/8]$ and $s\in(t_0,t_0+\upsilon_0]$,
\begin{equation*}
\begin{split}
&\;\EE^{\Lambda_{s-}-y}\big[w(t_0,R_{s-t_0})\,\bone_{\{R_{s-t_0}<\Lambda_{t_0}\}}\big] \\
&\leq(1\!-\!\psi(0))\,\PP^{\Lambda_{s-}-y}\big(R_{s-t_0}\!\in\![\Lambda_{t_0}\!-\!y_0/4,\Lambda_{t_0})\big)
+(1\!-\!\psi(y_0/4))\,\PP^{\Lambda_{s-}-y}\big(R_{s-t_0}\!\in\![\Lambda_{t_0}-y_0,\Lambda_{t_0}-y_0/4)\big) \\
&\quad+\frac12\psi\Big(\frac{y_0}{4}\Big)\,\PP^{\Lambda_{s-}-y}\big(R_{s-t_0}\in[\Lambda_{t_0}\!-\! y_0,\Lambda_{t_0}\!-\! y_0/4)\big) \\
&\leq 1-\psi(0)\,\PP^{\Lambda_{s-}-y}\big(R_{s-t_0}\in[\Lambda_{t_0}-y_0/4,\Lambda_{t_0})\big)
-\frac12\psi(y_0/4)\,\PP^{\Lambda_{s-}-y}\big(R_{s-t_0}\in[\Lambda_{t_0}-y_0,\Lambda_{t_0}-y_0/4)\big) \\
&\leq1-\psi(0)\frac{c_{\Lambda_{t_0},y_0}}{2}
-\psi(y_0/4)\frac{\pi y_0^3}{768}(2\pi(s-t_0))^{-\frac32}e^{-\frac{y_0^2}{8(s-t_0)}}.
\end{split}
\end{equation*}
We conclude the proof upon letting 
\begin{IEEEeqnarray*}{rCl}
c(\upsilon):=\sup_{\upsilon\le s'\le\upsilon_0}
\bigg(1-\psi(0)\frac{c_{\Lambda_{t_0},y_0}}{2}
-\psi(y_0/4)\frac{\pi y_0^3}{768}(2\pi s')^{-\frac32}e^{-\frac{y_0^2}{8s'}}\bigg)\in(0,1)
\end{IEEEeqnarray*}
for $\upsilon\in(0,\upsilon_0)$, and if $\psi(0)>0$, extend $c$ to $\upsilon=0$ according to
\begin{IEEEeqnarray*}{rCl}
\qquad\qquad\qquad\qquad\quad\;\;\;\qquad\qquad c(0):=1-\psi(0)\frac{c_{\Lambda_{t_0},y_0}}{2}\in(0,1). \qquad\qquad\qquad\qquad\qquad\qquad\quad\;\;\; \qed
\end{IEEEeqnarray*}

\begin{proposition}
Fix an arbitrary $t_0\in[0,\zeta)$ and assume that there are $y_0>0$, $\upsilon_0\in(0,\zeta-t_0)$ and $c\in[0,1)$ for which
\begin{equation}\label{below1}
\int_0^y w(t-,\Lambda_{t-}-y')(\Lambda_{t-}-y')^{2}\,\mathrm{d}y'
\le c\int_0^y (\Lambda_{t-}-y')^2\,\mathrm{d}y',\quad y\in[0,y_0],\quad t\in[t_0,t_0+\upsilon_0].
\end{equation}
Then, there exist $\upsilon\in(0,\upsilon_0)$ and $L_{t_0,t_0+\upsilon_0}<\infty$ such that
\begin{IEEEeqnarray*}{rCl}
0\leq\Lambda_t-\Lambda_s\leq L_{t_0,t_0+\upsilon_0}\sqrt{s-t},\quad t_0\le t\le s\le(t+\upsilon)\wedge(t_0+\upsilon_0).   
\end{IEEEeqnarray*}
That is, $\Lambda$ is $1/2$-Hölder continuous on $[t_0,t_0+\upsilon_0]$.
\end{proposition}

\noindent\textbf{Proof.} By the assumption~\eqref{below1} and the physicality condition~\eqref{intro: phys}, 
$\Lambda$ is continuous on $[t_0,t_0+\upsilon_0]$ and we have $w(t-,x)=w(t,x)$, $x\neq\Lambda_t$ for any $t\in[t_0,t_0+\upsilon_0]$ (see the paragraph following \eqref{eq.sec3.sergey.growth.cond.t0}). Moreover, we can find an $\upsilon\in(0,\upsilon_0\wedge1)$ satisfying
\begin{IEEEeqnarray*}{rCl}
\sup_{t_0\le t\le s\le(t+\upsilon)\wedge(t_0+\upsilon_0)}\,(\Lambda_t-\Lambda_s)\leq y_0.    
\end{IEEEeqnarray*}
To estimate $\Lambda_t-\Lambda_s$ more accurately for $t_0\le t\le s\le(t+\upsilon)\wedge(t_0+\upsilon_0)$, we rely on \eqref{eq.sec3.sergey.growth.cond.t0}:
\begin{IEEEeqnarray*}{rCl}
\frac{\Lambda_t^3}{3}-\frac{\Lambda_s^3}{3}&\leq&\int_0^{\infty}w(t,x)\,\PP^x(\tau_{\Lambda_{t+\cdot}}\leq s-t)\,\nu(\mathrm{d}x) \\
&\leq&\int_{\Lambda_s}^{\Lambda_t} w(t,x)\,\nu(\mathrm{d}x)
+\int_{(0,\Lambda_s)\cup(\Lambda_t,\infty)} \|w(t,\cdot)\|_\infty\,\PP^x(\tau_{\Lambda_{t+\cdot}}\leq s-t)\,\nu(\mathrm{d}x). 
\end{IEEEeqnarray*}
Thanks to the assumption \eqref{below1},
\begin{IEEEeqnarray*}{rCl}
\int_{\Lambda_{s}}^{\Lambda_{t}} w(t,x)\,\nu(\mathrm{d}x)
=\int_0^{\Lambda_t-\Lambda_s} w(t,\Lambda_t-y')(\Lambda_t-y')^2\,\mathrm{d}y
\leq c\Big(\frac{\Lambda_t^3}{3}-\frac{\Lambda_s^3}{3}\Big).
\end{IEEEeqnarray*}

\smallskip

For $x\in(0,\Lambda_s)\cup(\Lambda_t,\infty)$ we apply the tail bounds of Lemma \ref{lm:BesselRunning}:
\begin{IEEEeqnarray*}{rCl}
&&\PP^x(\tau_{\Lambda_{t+\cdot}}\le s-t)
\le\PP^x\big(\max_{0\le s'\le s-t} R_{s'}\ge\Lambda_{s}\big)
\leq12\Phi\bigg(\!-\frac{\Lambda_{s}-x}{\sqrt{3(s-t)}}\bigg),\quad x\in(0,\Lambda_{s}), \\
&&\PP^x(\tau_{\Lambda_{t+\cdot}}\leq s-t)
\leq\PP^x\big(\min_{0\le s'\le s-t} R_{s'}\le\Lambda_t\big)
\leq2\Phi\bigg(\!-\frac{x-\Lambda_{t}}{\sqrt{s-t}}\bigg),\quad x\in(\Lambda_t,\infty).
\end{IEEEeqnarray*}
As a result, 
\begin{IEEEeqnarray*}{rCl}
&&\int_{(0,\Lambda_s)\cup(\Lambda_t,\infty)}\PP^x(\tau_{\Lambda_{t+\cdot}}\leq s-t)\,\nu(\mathrm{d}x) \\
&& \le\int_0^{\Lambda_{s}}12\Phi\bigg(\!-\frac{\Lambda_{s}-x}{\sqrt{3(s-t)}}\bigg)\,\nu(\mathrm{d}x)
+\int_{\Lambda_t}^\infty 2\Phi\bigg(\!-\frac{x-\Lambda_{t}}{\sqrt{s-t}}\bigg)\,\nu(\mathrm{d}x) \\
&& =\sqrt{s-t}\,\bigg(\int_0^{\frac{\Lambda_{s}}{\sqrt{3(s-t)}}}12\sqrt{3}\Phi(-x)
\big(\Lambda_{s}-\sqrt{3(s-t)}\,x\big)^2\,\mathrm{d}x 
+\int_0^\infty2\Phi(-x)\big(\Lambda_t+\sqrt{s-t}\,x\big)^2\,\mathrm{d}x\bigg) \\
&&\le C_{\Lambda_t}\sqrt{s-t},
\end{IEEEeqnarray*}
where
\begin{IEEEeqnarray*}{rCl}
C_{\Lambda_t}:=\big(12\sqrt{3}+4\big)\Lambda_{t}^2
\int_0^\infty\Phi(-x)\,\mathrm{d}x+4\int_0^\infty\Phi(-x)x^2\,\mathrm{d}x<\infty.
\end{IEEEeqnarray*}

\smallskip

Since $\|w(t,\cdot)\|_\infty\leq\|w(0-,\cdot)\|_\infty$ by \eqref{eq.sec3.sergey.w.def}, we have shown that
\begin{IEEEeqnarray*}{rCl}
\frac{\Lambda_t^3}{3}-\frac{\Lambda_s^3}{3}
\le c\Big(\frac{\Lambda_t^3}{3}-\frac{\Lambda_s^3}{3}\Big)
+C_{\Lambda_t}\|w(0-,\cdot)\|_\infty\sqrt{s-t}.
\end{IEEEeqnarray*}
Consequently,
\begin{IEEEeqnarray*}{rCl}
\frac{\Lambda_t^3}{3}-\frac{\Lambda_s^3}{3}
\le \frac{1}{1-c}C_{\Lambda_t}\|w(0-,\cdot)\|_\infty\sqrt{s-t}.    
\end{IEEEeqnarray*}
We note that $\frac{\Lambda_t^3}{3}-\frac{\Lambda_s^3}{3}\ge\Lambda_{t_0+\upsilon_0}^2(\Lambda_t-\Lambda_s)$ and that $C_{\Lambda_t}\le C_{\Lambda_{t_0}}$. The proof is finished by letting
\begin{IEEEeqnarray*}{rCl}
\qquad\qquad\qquad\qquad\qquad\quad\;
L_{t_0,t_0+\upsilon_0}:=\frac{1}{(1-c)\Lambda_{t_0+\upsilon_0}^2} C_{\Lambda_{t_0}}\|w(0-,\cdot)\|_\infty. 
\qquad\qquad\qquad\qquad\qquad\quad\; \qed  
\end{IEEEeqnarray*}

\begin{proposition}\label{KrySaf} 
Fix an arbitrary $t_0\in[0,\zeta)$ and assume that $\Lambda$ is $1/2$-Hölder continuous on $[t_0,t_0+\upsilon_0]$ for some $\upsilon_0\in(0,\zeta-t_0)$. Then, for all $\upsilon\in(0,\upsilon_0)$, there exist $C_{t_0,\upsilon_0,\upsilon}<\infty$, $y_0>0$ and $\varrho\in(0,1)$ such that
\begin{IEEEeqnarray*}{rCl}
w(s,\Lambda_s-y)\leq C_{t_0,\upsilon_0,\upsilon}\,y^\varrho,\quad y\in(0,y_0],\quad s\in[t_0+\upsilon,t_0+\upsilon_0].
\end{IEEEeqnarray*}
\end{proposition}

\begin{remark}
Since $w\le\|w(0-,\cdot)\|_\infty$ by \eqref{eq.sec3.sergey.w.def}, we can enlarge the constant $C_{t_0,\upsilon_0,\upsilon}<\infty$ to ensure
\begin{IEEEeqnarray*}{rCl}
w(s,\Lambda_s-y)\leq C_{t_0,\upsilon_0,\upsilon}\,y^\varrho,\quad y\in(0,\Lambda_s),\quad s\in[t_0+\upsilon,t_0+\upsilon_0].
\end{IEEEeqnarray*}
\end{remark}

\noindent\textbf{Proof of Proposition \ref{KrySaf}.} Please see Appendix \ref{app:KrySaf}. \qed

\medskip


Thanks to the non-increasing nature of $\Lambda$, it is much easier to estimate $w$ in the liquid phase. 

\begin{lemma}\label{lm:wLinearGrowthaboveLambda}
For any $t_0\in(0,\zeta)$, it holds
\begin{equation}\label{w liquid}
\frac{w(t_0,\Lambda_{t_0}+y)}{y}
\le 2\|w(0-,\cdot)\|_\infty\!\left(\frac{1}{\sqrt{2\pi t_0}}+\frac{1}{\Lambda_{t_0}}\right),\quad y>0. 
\end{equation}
\end{lemma} 

\begin{remark}
Since $w$ is bounded and by \eqref{w liquid}, any $t_0\in[0,\zeta)$ admits $C<\infty$, $\varrho\in(0,1)$ such that
\begin{IEEEeqnarray*}{rCl}
w(t_0,\Lambda_{t_0}+y)\leq Cy^\varrho,\quad y>0.    
\end{IEEEeqnarray*}
This estimate of $w$ for the liquid phase is in line with that of Proposition \ref{prop:LipschitzLambda1} for the solid phase.
\end{remark}

\noindent\textbf{Proof of Lemma \ref{lm:wLinearGrowthaboveLambda}.} We recall from \eqref{eq.sec3.sergey.w.def} that
\begin{equation}\label{eq:wabovesplit1}
w(t,\Lambda_{t_0}+y)
=\EE^{\Lambda_{t_0}+y}[w(0-,R_{t_0})\,\bone_{\{\tau_{\Lambda_{t_0-\cdot}}\wedge\tau_{r_0}>t_0}\}]
\le\|w(0-,\cdot)\|_\infty\,\PP^{\Lambda_{t_0}+y}(\tau_{\Lambda_{t_0-\cdot}}>t_0).
\end{equation}
By comparing the Bessel process $R$ to a Brownian motion with constant drift and subsequently applying the Girsanov theorem we find that
\begin{IEEEeqnarray*}{rCl}
\PP^{\Lambda_{t_0}+y}(\tau_{\Lambda_{t_0-\cdot}}>t_0)
&\le&\PP\big(\min_{0\le s\le t_0} (B_s+s/\Lambda_{t_0})>-y\big) \\
&=&\int_0^y \int_{-\infty}^m 
e^{-\frac{b}{\Lambda_{t_0}}-\frac{t_0}{2\Lambda_{t_0}^2}}\,
\frac{2(2m-b)}{t_0\sqrt{2\pi t_0}}\,e^{-\frac{(2m-b)^2}{2t_0}}\,\mathrm{d}b\,\mathrm{d}m 
\le 2\left(\frac{1}{\sqrt{2\pi t_0}}+\frac{1}{\Lambda_{t_0}}\right)y,
\end{IEEEeqnarray*}
where the latter inequality is shown in Lemma \ref{lem:BMestimates}(c) (see Appendix \ref{BMestimate}). \qed




\subsection{Interior Smoothness and Derivatives Estimates} \label{subsec:InteriorRegularity}

Let us consider $w$ on the interior set
\begin{IEEEeqnarray*}{rCl}
\stackrel{\circ}{D}\;:=\{(t,x)\in(0,\zeta)\times(r_0,\infty):\,x>\Lambda_{t-}\;\text{or}\;x<\Lambda_t\}.
\end{IEEEeqnarray*}

\begin{proposition}
\label{prop:InteriorSmoothness}
On $\stackrel{\circ}{D}$, the function $w$ is of class $C^\infty$  and satisfies pointwise 
\begin{IEEEeqnarray*}{rCl}
\partial_tw(t,x)=\frac12\partial_{xx}w(t,x)+\frac1x\partial_xw(t,x).
\end{IEEEeqnarray*}
\end{proposition}

\noindent\textbf{Proof.} Please see Appendix \ref{app:Cinf}. \qed

\medskip

The following derivatives estimates are used below to establish the analyticity of $w$ in $x$. 

\begin{proposition}\label{prop:GradientEstimateForw1}
There exists a universal constant $C<\infty$ such that 
\begin{IEEEeqnarray*}{rCl}
|\partial_xw(t,\Lambda_t-y)|\le\frac{C}{y}\,\|w(0-,\cdot)\|_\infty,\quad y\in\Big(0,\frac{\Lambda_t-r_0}{2}\wedge\sqrt{t}\Big),\quad t\in(0,\zeta).   
\end{IEEEeqnarray*}
More generally, for $(t,x)\in\stackrel{\circ}{D}$, define the parabolic distance to the boundary
\begin{IEEEeqnarray*}{rCl}
d_{\text{par}}(t,x)=\sup\big\{r>0:\,[t-r^2,t]\times[x-r,x+r]\subset\,\stackrel{\circ}{D}\!\big\}.
\end{IEEEeqnarray*}
Then, with the same $C<\infty$ as above and for all $m,n\in\NN_0$,
\begin{IEEEeqnarray*}{rCl}
|\partial_t^m\partial_x^n w(t,x)|\leq\frac{C^{2m+n}(2m+n)^{2m+n}}{d_{\text{par}}(t,x)^{2m+n}}\,\|w(0-,\cdot)\|_\infty,\quad (t,x)\in\stackrel{\circ}{D}.
\end{IEEEeqnarray*}
\end{proposition}

\noindent\textbf{Proof.} We prove the general statement. Consider the function $u(t,z):=w(t,|z|)$ on
\begin{IEEEeqnarray*}{rCl}
\stackrel{\circ}{D_3}\,:= \{(t,z)\in(0,\zeta)\times\RR^3:|z|\in(r_0,\Lambda_t)\;\text{or}\;|z|>\Lambda_{t-}\},
\end{IEEEeqnarray*}
which is of class $C^\infty$ and satisfies pointwise $\partial_t u(t,z)=\frac12\Delta_z u(t,z)$. Given $(t,x)\in\stackrel{\circ}{D}$, we pick $r>0$ so that $[t-r^2,t]\times[x-r,x+r]\subset\,\stackrel{\circ}{D}$ and apply the derivatives estimate for the heat equation of \cite[Theorem 8.4.4]{KrylovHolder} on the parabolic cylinder $(t-r^2,t)\times B_r((x,0,0))\subset\,\stackrel{\circ}{D_3}$ to get
\begin{IEEEeqnarray*}{rCl}
|\partial_t^m\partial_{z_1}^n u(t,(x,0,0))|
\leq\frac{C^{2m+n}(2m+n)^{2m+n}}{r^{2m+n}}\sup_{(t-r^2,t)\times B_r((x,0,0))} |u|,    
\end{IEEEeqnarray*}
where $C<\infty$ is a constant that only depends on the dimension $d=3$. (It differs from the $N(d)$ in \cite{KrylovHolder} by a constant factor since we use a different normalization of the heat equation.) Lastly, we observe $\partial_x^n w(t,x)=\partial_{z_1}^n u(t,(x,0,0))$ and  $0\le u(t,z)\le\|w(0-,\cdot)\|_\infty$, and take $r\uparrow d_{\text{par}}(t,x)$. \qed 

\begin{proposition}
\label{prop:InteriorAnalyticity}
For any $(t,x)\in\,\stackrel{\circ}{D}$, there is a neighborhood of $(t,x)$ in $\RR\times\CC$ to which $w$ and its time derivatives of all orders can be extended so that they are analytic in the space variable and jointly continuous in both variables. 
\end{proposition}

\noindent\textbf{Proof.} Let $u$ and $\stackrel{\circ}{D_3}$ be as in the proof of Proposition \ref{prop:GradientEstimateForw1}. The real analyticity of $u$ in $z$ is classical, see, e.g., \cite[Exercise 8.4.7]{KrylovHolder}.~The analyticity of $w$ follows directly from $w(t,x)=u(t,(x,0,0))$. More specifically, the analytic extension of $w$ reads
\begin{IEEEeqnarray*}{rCl}
w(t,\cdot):=\sum_{n=0}^\infty\frac{\partial_x^n w(t,x)}{n!}(\,\cdot-x)^n
=\sum_{n=0}^\infty\frac{\partial_{z_1}^n u(t,(x,0,0))}{n!}(\,\cdot-x)^n,
\end{IEEEeqnarray*}
where the radius of convergence is bounded away from zero locally uniformly due to the derivatives estimate of Proposition \ref{prop:GradientEstimateForw1}. Moreover, in each domain of convergence, the above power series can be differentiated (in both variables) term-by-term, which yields an extension of $\partial_t^m w$ by
\begin{IEEEeqnarray*}{rCl}
\partial_t^mw(t,\cdot):=\sum_{n=0}^\infty \frac{\partial_t^m\partial_x^n w(t,x)}{n!}(\,\cdot-x)^n.
\end{IEEEeqnarray*}
Indeed, by Proposition \ref{prop:GradientEstimateForw1} the latter power series in $\,\cdot-x$ has a radius of convergence of at least
\begin{IEEEeqnarray*}{rCl}
\qquad\qquad\qquad\qquad\qquad\quad\,
\bigg(\limsup_{n\to\infty}\,\bigg|\frac{\partial_t^m\partial_x^n w(t,x)}{n!}\bigg|^{1/n}\bigg)^{-1}
\ge\frac{d_{\text{par}}(t,x)}{Ce}. \qquad\qquad\qquad\qquad\qquad\quad\, \qed
\end{IEEEeqnarray*}

\subsection{Stefan Growth Condition} \label{subsec:StefanGrowth}

 In this subsection, we connect $\Lambda'_t$ to $\partial_x w(t,\Lambda_t\pm)$.
 
 \medskip

\noindent\textit{A few words on notation.} For given initial data $\Lambda_{0-}\!>\! r_0$ and $w(0-,\cdot)$, time horizon $T\!>\!0$, and non-increasing right-continuous function $\lambda\!:[-1,T]\to\mathbb{R}$ with $\lambda_t\equiv\Lambda_{0-}$ on $[-1,0)$, we denote by $w^\lambda$ the corresponding Feynman-Kac solution:
\begin{equation}\label{eq.sec3.sergey.wLambda.def}
w^\lambda(t,x):=\EE^x[w(0-,R_t)\,\bone_{\{\tau_{\lambda_{t-\cdot}}\wedge\tau_{r_0}>t\}}],
\quad x>0, \quad t\in[0,\zeta^\lambda\wedge T),
\end{equation}
where $\tau_{\lambda_{t-\cdot}}$ is defined per \eqref{eq.sec3.sergey.tauLambda.tminus.def} and $\zeta^\lambda:=\inf\{s\ge0\!:\lambda_s\le r_0\}$. We also let $\Gamma[\lambda]$ be the boundary ``generated'' by $\lambda$:
\begin{equation}\label{eq.section3.sergey.Gamma.def}
\begin{split}
\frac{\Lambda_{0-}^3}{3}-\frac{\Gamma[\lambda]_t^3}{3}
&:=\Gamma^+[\lambda]_t+\Gamma^-[\lambda]_t\\
&:=\int_{r_0}^{\Lambda_{0-}} \! w(0-,x)\,
\PP^x(\tau^-_{\lambda_{0+\cdot}}\!\le t\!\wedge\!\tau_{r_0})\,\nu(\mathrm{d}x)
+\int_{\Lambda_{0-}}^\infty \! w(0-,x)\,\PP^x(\tau^-_{\lambda_{0+\cdot}}\!\le t)\,\nu(\mathrm{d}x),
\end{split}
\end{equation}
with $\tau^-_{\lambda_{0+\cdot}}$ defined per \eqref{eq.sec3.sergey.tauLambda.minus.tplus.def}. Probabilistic solutions $\Lambda$ to \eqref{intro: heat rs}--\eqref{intro: Stefan rs} are then characterized by the fixed point equation $\Gamma[\Lambda]=\Lambda$. Moreover, arguing as for \eqref{eq.sec3.sergey.growth.cond.t0} we see that
\begin{IEEEeqnarray*}{rCl}
\Gamma^+[\lambda]_t&=&\Gamma^+[\lambda]_0 + \int_{r_0}^{\lambda_{0}}w^\lambda(0,x)\,\PP^x(\tau_{\lambda_{0+\cdot}}\leq t\wedge\tau_{r_0})\,\nu(\mathrm{d}x), \\
\Gamma^-[\lambda]_t&=&\Gamma^-[\lambda]_0 + \int_{\lambda_{0}}^\infty w^\lambda(0,x)\,\PP^x(\tau_{\lambda_{0+\cdot}}\le t)\,\nu(\mathrm{d}x).
\end{IEEEeqnarray*}



\begin{proposition}\label{prop:section3.sergey.prop.Lambda.via.w}
Suppose that $\zeta^\lambda>0$ and that $\lambda_0-\lambda_t=o(\sqrt{t})$ as $t\downarrow0$. Then,
\begin{eqnarray*}
\frac{\lambda_0^2}{2}\,\liminf_{y\downarrow0}\,\frac{w^\lambda(0,\lambda_0-y)}{y}
\le\liminf_{t\downarrow0}\,\frac{\Gamma^+[\lambda]_t-\Gamma^+[\lambda]_0}{t}
\le\limsup_{t\downarrow0}\,\frac{\Gamma^+[\lambda]_t-\Gamma^+[\lambda]_0}{t} \;\; \\
\le\frac{\lambda_0^2}{2}\,\limsup_{y\downarrow0}\,\frac{w^\lambda(0,\lambda_0-y)}{y}, \\
\frac{\lambda_0^{2}}{2}\,\liminf_{y\downarrow0}\,\frac{w^\lambda(0,\lambda_0+y)}{y}
\le\liminf_{t\downarrow0}\,\frac{\Gamma^-[\lambda]_t-\Gamma^-[\lambda]_0}{t}
\le\limsup_{t\downarrow0}\,\frac{\Gamma^-[\lambda]_t-\Gamma^-[\lambda]_0}{t} \;\, \\
\le\frac{\lambda_0^2}{2}\,\limsup_{y\downarrow0}\,\frac{w^\lambda(0,\lambda_0+y)}{y}.
\end{eqnarray*}
\end{proposition}

\noindent\textbf{Proof. Step 1.} To obtain a tractable expression for $\Gamma^+[\lambda]_t$ we use a coupling construction. Let $\{R^x\}_{x\ge0}$ be the respective strong solutions of $\mathrm{d}R^x_t=1/R^x_t\,\mathrm{d}t+\mathrm{d}B_t$, $R^x_0=x$ and note that
\begin{IEEEeqnarray*}{rCl}
\forall\,x>x'\ge0:\quad R_t^x>R_t^{x'},\;\; t\ge0\quad\text{almost surely.}
\end{IEEEeqnarray*}
Next, for any fixed $t\in(0,\zeta^\lambda\wedge T)$, we define the random variable
\begin{IEEEeqnarray*}{rCl}
\overline{Y}_t=\sup\Big\{y\in(0,\lambda_0-r_0):\,\inf_{0\le s\le t\wedge\tau^{\lambda_0-y}_{r_0}} (\lambda_s-R_s^{\lambda_0-y})\le0\Big\}\in(0,\lambda_0-r_0]\quad\text{almost surely,}
\end{IEEEeqnarray*}
where (for this proof only) we have introduced the stopping times
\begin{IEEEeqnarray*}{rCl}
\tau^x_{r_0}:=\inf\{s\geq0:\,R^x_s\le r_0\},
\quad \tau^x_\lambda:=\inf\{s\ge0:\,(R_s^x-\lambda_s)(x-\lambda_0)\le 0\}.    
\end{IEEEeqnarray*}
Then, for all $y\in(0,\lambda_0-r_0)$:
\begin{IEEEeqnarray*}{rCl}
\tau^{\lambda_0-y}_\lambda\leq t\wedge\tau^{\lambda_0-y}_{r_0}\quad\Longleftrightarrow\quad y\leq\overline{Y}_t.
\end{IEEEeqnarray*}
Hence, Fubini's theorem yields
\begin{equation*}
\begin{split}
\Gamma^+[\lambda]_t-\Gamma^+[\lambda]_0
&=\int_{r_0}^{\lambda_0} w^\lambda(0,x)\,
\PP(\tau^x_\lambda\le t\wedge\tau^x_{r_0})\,x^2\,\mathrm{d}x \\
&=\EE\bigg[\int_0^{\lambda_0-r_0} w^\lambda(0,\lambda_0-y)\,\bone_{\{\tau^{\lambda_0-y}_\lambda
\le t\wedge\tau^{\lambda_0-y}_{r_0}\}}(\lambda_0-y)^2\,\mathrm{d}y\bigg] \\
&=\EE\bigg[\int_0^{\overline{Y}_t} w^\lambda(0,\lambda_0-y)(\lambda_0-y)^2\,\mathrm{d}y\bigg].
\end{split}
\end{equation*}

\smallskip

To simplify notation, we set
\begin{IEEEeqnarray*}{rCl}
\overline{L}=\limsup_{y\downarrow0}\,\frac{w^\lambda(0,\lambda_0-y)}{y},
\quad\underline{L}=\liminf_{y\downarrow0}\,\frac{w^\lambda(0,\lambda_0-y)}{y}  
\end{IEEEeqnarray*}
and observe that $\overline{L}\ge\underline{L}\ge0$. Assume $\overline{L}<\infty$ and note that otherwise the last inequality in the first statement of the proposition becomes trivial. If $\overline{L}<\infty$, then for any given $\varepsilon>0$, there exists some $y_0\in(0,(\lambda_0-r_0)/4)$ such that
\begin{IEEEeqnarray*}{rCl}
(\underline{L}-\varepsilon)y\le w^\lambda(0,\lambda_0-y)\le(\overline{L}+\varepsilon)y, 
\quad y\in(0,y_0],
\end{IEEEeqnarray*}
and also
\begin{equation}\label{eq:estintw1}
\Big(\frac{\lambda_0^2}{2}\,\underline{L}-\varepsilon\Big)y^2
\le\int_0^y w^\lambda(0,\lambda_0-y')(\lambda_0-y')^2\,\mathrm{d}y'
\le\Big(\frac{\lambda_0^2}{2}\,\overline{L}+\varepsilon\Big)y^2,\quad y\in(0,y_0].
\end{equation}
With this $y_0>0$, we can take $t_0>0$ small enough to ensure
\begin{IEEEeqnarray*}{rCl}
\lambda_0-\lambda_{t_0}\le\frac{y_0}{4}\quad\text{and}\quad
\frac{t_0}{\lambda_0-2y_0}\le\frac{y_0}{4}.    
\end{IEEEeqnarray*}
Since $\lambda_0-\lambda_t=o(\sqrt{t})$ as $t\downarrow0$ by assumption, we can shrink $t_0>0$ if necessary, to guarantee 
\begin{IEEEeqnarray*}{rCl}
\lambda_0-\lambda_t\le\varepsilon\sqrt{t},\quad t\in(0,t_0].
\end{IEEEeqnarray*}

\smallskip

Consider now the running maximum processes
\begin{IEEEeqnarray*}{rCl}
S_t:=\max_{0\le s\le t} B_s,\quad B_t^*:=\max_{0\le s\le t} |B_s|.
\end{IEEEeqnarray*}
Then, for $t\in(0,t_0]$, on the event $\{B_t^*\le y_0/4\}$, it holds
\begin{equation} \label{eq.sec3.sergey.prop.3.13.eq1}
R_s^{\lambda_0-y_0}\ge R_0^{\lambda_0-y_0}+B_s
\ge\lambda_0-2y_0>r_0,\quad s\in[0,t].
\end{equation}
Therefore, $\tau^{\lambda_0-y}_{r_0}>t$ for $y\le y_0$ on the event $\{B_t^*\le y_0/4\}$, and
\begin{equation*}
\lambda_s-R_s^{\lambda_0-y_0}
=\lambda_s-R_0^{\lambda_0-y_0}
-\int_0^s\frac{1}{R_{s'}^{\lambda_0-y_0}}\,\mathrm{d}s'-B_s
\ge(\lambda_s-\lambda_0)+y_0-\frac{s}{\lambda_0-2y_0}-S_s
\ge\frac{y_0}{4}>0,\quad s\in[0,t].
\end{equation*}
Thus, on the event $\{B_t^*\le y_0/4\}$,
\begin{equation} \label{eq.sec3.sergey.prop.3.13.eq2}
\overline{Y}_t\le y_0
\quad\text{and}\quad\tau^{\lambda_0-\overline{Y}_t}_{r_0}>t.
\end{equation}

\smallskip

Next, we aim to compare $\overline{Y}_t$ to $S_t$. On the one hand, we have the trivial bound
\begin{IEEEeqnarray*}{rCl}
\inf_{0\le s\le t\wedge\tau^{\lambda_0-y}_{r_0}} (\lambda_s-R_s^{\lambda_0-y})
&=&\inf_{0\le s\le t\wedge\tau^{\lambda_0-y}_{r_0}} \bigg(\lambda_s-\lambda_0+y-\int_0^s\frac{1}{R_{s'}^{\lambda_0-y}}\,\mathrm{d}s'-B_s\bigg) \\
&\leq&y-\sup_{0\le s\le t\wedge\tau^{\lambda_0-y}_{r_0}} B_s
=y-S_{t\wedge\tau^{\lambda_0-y}_{r_0}},
\end{IEEEeqnarray*}
which yields $\overline{Y}_t\ge S_t$ on the event $\{B_t^*\le y_0/4\}$. On the other hand, using \eqref{eq.sec3.sergey.prop.3.13.eq1}, \eqref{eq.sec3.sergey.prop.3.13.eq2} we obtain for $t\in(0,t_0]$, on the event $\{B_t^*\le y_0/4\}$,
\begin{IEEEeqnarray*}{rCl}
0\ge\lambda_t-R_t^{\lambda_0-\overline{Y}_t}
=\lambda_t-\lambda_0+\overline{Y}_t
-\int_0^t\frac{1}{R_s^{\lambda_0-\overline{Y}_t}}\,\mathrm{d}s-B_t
\ge\lambda_t-\lambda_0+\overline{Y}_t-\frac{t}{\lambda_0-2y_0}-S_t.
\end{IEEEeqnarray*}
Hence, for $t\in(0,t_0]$, on the event $\{B_t^*\le y_0/4\}$,
\begin{IEEEeqnarray*}{rCl}
\overline{Y}_t\le S_t+\frac{t}{\lambda_0-2y_0}+\lambda_0-\lambda_t
\le S_t+\frac{t}{\lambda_0-2y_0}+\varepsilon\sqrt{t}.
\end{IEEEeqnarray*}

\smallskip

\noindent\textbf{Step 2.} We can now prove the desired result. Rewrite
\begin{IEEEeqnarray*}{rCl}
\Gamma^+[\lambda]_t-\Gamma^+[\lambda]_0&=&\EE\bigg[\int_0^{\overline{Y}_t} w^\lambda(0,\lambda_0-y)(\Lambda_0-y)^2\,\mathrm{d}y\,\bone_{\{B_t^*\le y_0/4\}}\bigg] \\
&&+\EE\bigg[\int_0^{\overline{Y}_t} w^\lambda(0,\lambda_0-y)(\lambda_0-y)^2\,\mathrm{d}y\,\bone_{\{B_t^*>y_0/4\}}\bigg].
\end{IEEEeqnarray*}
For the lower bound, we combine the non-negativity of $w^\lambda$ and \eqref{eq:estintw1} to derive
\begin{IEEEeqnarray*}{rCl}
\Gamma^+[\lambda]_t-\Gamma^+[\lambda]_0
\ge\EE\bigg[\int_0^{\overline{Y}_t} w^\lambda(0,\lambda_0-y)(\lambda_0-y)^2\,\mathrm{d}y
\,\bone_{\{B_t^*\le y_0/4\}}\bigg] 
\ge\EE\bigg[\Big(\frac{\lambda_0^2}{2}\underline{L}-\varepsilon\Big)_+\,\overline{Y}_t^2
\,\bone_{\{B_t^*\le y_0/4\}}\bigg] \\
\ge\Big(\frac{\lambda_0^2}{2}\underline{L}-\varepsilon\Big)_+\,
\EE[S_t^2\,\bone_{\{B_t^*\le y_0/4\}}].
\end{IEEEeqnarray*}
We note further that 
\begin{IEEEeqnarray*}{rCl}
\EE[S_t^2\,\bone_{\{B_t^*\le y_0/4\}}]
=\EE[S_t^2]-\EE[S_t^2\,\bone_{\{B_t^*>y_0/4\}}]
\ge t-\EE[S_t^4]^{1/2}\,\PP(B_t^*>y_0/4)^{1/2} 
= t-o(t)\quad\text{as}\quad t\downarrow0.
\end{IEEEeqnarray*}

\smallskip

For the upper bound, we use $w^\lambda(0,x)\leq\|w(0-,\cdot)\|_\infty$ and \eqref{eq:estintw1} to find
\begin{IEEEeqnarray*}{rCl}
\Gamma^+[\lambda]_t-\Gamma^+[\lambda]_0
&\leq&\EE\bigg[\Big(\frac{\lambda_0^2}{2}\overline{L}+\varepsilon\Big)\overline{Y}_t^2
\,\bone_{\{B_t^*\le y_0/4\}}\bigg]
+\int_0^{\lambda_0}\|w(0-,\cdot)\|_\infty\, x^2\,\mathrm{d}x
\,\PP(B_t^*>y_0/4) \\
&\le&\Big(\frac{\lambda_0^2}{2}\overline{L}+\varepsilon\Big)\,
\EE\bigg[\Big(S_t+\frac{t}{\lambda_0-2y_0}+\varepsilon\sqrt{t}\Big)^2\bigg]
+\frac{\|w(0-,\cdot)\|_\infty\,\lambda_0^3}{3}\,\PP(B_t^*>y_0/4).
\end{IEEEeqnarray*}
In view of $\EE[S_t^2]=t$, we have $\EE[S_t]\le\sqrt{t}$, from which it follows that, as $t\downarrow0$,
\begin{IEEEeqnarray*}{rCl}
\EE\bigg[\Big(S_t+\frac{t}{\lambda_0-2y_0}+\varepsilon\sqrt{t}\Big)^2\bigg]
=\EE\bigg[S_t^2+2\Big(\frac{t}{\lambda_0-2y_0}+\varepsilon\sqrt{t}\Big)S_t
+\Big(\frac{t}{\lambda_0-2y_0}+\varepsilon\sqrt{t}\Big)^2\bigg] \\
\le t+2\Big(\frac{t}{\lambda_0-2y_0}+\varepsilon\sqrt{t}\Big)\sqrt{t}
+\Big(\frac{t}{\lambda_0-2y_0}+\varepsilon\sqrt{t}\Big)^2
=(1+2\varepsilon+\varepsilon^2)t+O(t^{3/2}).
\end{IEEEeqnarray*}
Since $\PP(B_t^*>y_0/4)=o(t)$, one can take $t\downarrow0$ and then $\varepsilon\downarrow0$ to arrive at 
\begin{IEEEeqnarray*}{rCl}
\frac{\lambda_0^2}{2}\underline{L}
\le\liminf_{t\downarrow0}\,\frac{\Gamma^+[\lambda]_t-\Gamma^+[\lambda]_0}{t}
\le\limsup_{t\downarrow0}\,\frac{\Gamma^+[\lambda]_t-\Gamma^+[\lambda]_0}{t}
\le\frac{\lambda_0^2}{2}\overline{L}.  
\end{IEEEeqnarray*}

\smallskip

For the liquid phase, we use the same argument. To wit, for any fixed $t\in(0,\zeta^\lambda\wedge T)$, define
\begin{IEEEeqnarray*}{rCl}
\overline{Y}_t=\sup\Big\{y\in(0,\infty):\,\inf_{0\le s\le t} (R_s^{\lambda_0+y}-\lambda_s)\leq0\Big\}\ge0\quad\text{almost surely}
\end{IEEEeqnarray*}  
and observe that for all $y\in(0,\infty)$,
\begin{IEEEeqnarray*}{rCl}
\tau^{\lambda_0+y}_\lambda\leq t\quad\Longleftrightarrow
\quad y\leq\overline{Y}_t.
\end{IEEEeqnarray*}
Hence, by Fubini's theorem,
\begin{equation*}
\begin{split}
\Gamma^-[\lambda]_t-\Gamma^-[\lambda]_0=\int_{\lambda_0}^\infty w^\lambda(0,x)\,\PP(\tau^x_\lambda\le t)\,x^2\,\mathrm{d}x
&=\EE\bigg[\int_0^\infty w^\lambda(0,\lambda_0+y)\,
\bone_{\{\tau^{\lambda_0+y}_\lambda\le t\}}(\lambda_0+y)^2\,\mathrm{d}y\bigg] \\
&=\EE\bigg[\int_0^{\overline{Y}_t}w^\lambda(0,\lambda_0+y)(\lambda_0+y)^2\,\mathrm{d}y\bigg].
\end{split}
\end{equation*}
With $S^{-B}_t:=\max_{0\le s\le t}(-B_s)$, it is easy to see that
\begin{IEEEeqnarray*}{rCl}
S^{-B}_t-(\lambda_0-\lambda_t)-\frac{t}{\lambda_0-y_0}\le\overline{Y}_t\le S^{-B}_t
\end{IEEEeqnarray*}
as long as $\lambda_t\ge\lambda_0-y_0$, which allows to bound the latter expectation as for the solid phase above. \qed

\medskip



It is easy to see that the Markov property of probabilistic solutions, stated in \eqref{eq.sec3.sergey.growth.cond.t0.minus}, \eqref{eq.section3.sergey.w.back.Markov.tminus},
also holds for the pair $(\lambda,w^\lambda)$: For $0\leq t_0\leq t<\zeta^\lambda$,
\begin{IEEEeqnarray*}{rCl}
&& \frac{\Gamma[\lambda]_{t_0-}^3}{3}-\frac{\Gamma[\lambda]_t^3}{3}
=\int_{r_0}^{\infty}w^\lambda(t_0-,x)\,
\PP^x(\tau^-_{\lambda_{t_0+\cdot}}\le (t-t_0)\wedge\tau_{r_0})\,\nu(\mathrm{d}x), \\
&& w^\lambda(t,x)
=\EE^x[w^\lambda(t_0-,R_{t-t_0})\,\bone_{\{\tau_{\lambda_{t-\cdot}}\wedge\tau_{r_0}>t-t_0\}}],\quad x>0. 
\end{IEEEeqnarray*}
Thus, we get the following version of Proposition \ref{prop:section3.sergey.prop.Lambda.via.w} for arbitrary $t_0\in[0,\zeta^\lambda)$ in place of $t_0=0$.

\begin{corollary} \label{cor:BoundaryCondition2}
Let $t_0\in[0,\zeta^\lambda)$ be such that $\lambda_{t_0}-\lambda_{t_0+t}=o(\sqrt{t})$ as $t\downarrow0$ and that $\partial_x w^\lambda(t_0,x\pm)|_{x=\lambda_{t_0}}$ exist. Then, $t\mapsto\Gamma[\lambda]_t$ is right-differentiable at $t=t_0$ and
\begin{IEEEeqnarray*}{rCl}
\Gamma[\lambda]_{t_0}^2\,\frac{\mathrm{d}}{\mathrm{d}t}\Gamma[\lambda]_{t_0}
=\lambda_{t_0}^2
\Big(\frac12\partial_xw^\lambda(t_0,x-)|_{x=\lambda_{t_0}}
-\frac12\partial_xw^\lambda(t_0,x+)|_{x=\lambda_{t_0}}\Big).
\end{IEEEeqnarray*}
\end{corollary}

\subsection{Continuous Differentiability} \label{subsec:BoundaryRegularity}
The main result of this subsection is the next proposition.


\begin{proposition}\label{prop:LipschitzLambda1}
Fix an arbitrary $t_0\in[0,\zeta)$ and assume that for some $C<\infty$ and $\varrho\in(0,1)$, 
\begin{IEEEeqnarray*}{rCl}
&& 0\le w(t_0,\Lambda_{t_0}-y)\le Cy^\varrho,\quad y\in(0,\Lambda_{t_0}-r_0), \\
&& 0\le w(t_0,\Lambda_{t_0}+y)\le C(y^\varrho\wedge y),\quad y>0.  
\end{IEEEeqnarray*}
Then, there exists $\upsilon_0>0$ such that $\Lambda$ is continuously differentiable on $(t_0,t_0+\upsilon_0)$, and
both $w$ and $\partial_xw$ are continuous on
\begin{IEEEeqnarray*}{rCl}
\{(t,x)\in(t_0,t_0+\upsilon_0)\times(r_0,\infty):\,x\leq\Lambda_t\}\quad\text{and}\quad
\{(t,x)\in(t_0,t_0+\upsilon_0)\times(r_0,\infty):\,x\geq\Lambda_t\}.
\end{IEEEeqnarray*}
Moreover, the Stefan growth condition holds: 
\begin{IEEEeqnarray*}{rCl}
\Lambda'_t=\frac12\partial_xw(t,x-)|_{x=\Lambda_t}-\frac12\partial_xw(t,x+)|_{x=\Lambda_t},
\quad t\in(t_0,t_0+\upsilon_0).    
\end{IEEEeqnarray*}
\end{proposition}

\begin{remark}
By Proposition \ref{prop:LipschitzLambda1}, each $t_0\in[0,\zeta)$ admits an $\upsilon_0>0$ so that for all $\upsilon\in(0,\upsilon_0)$, there exist $C_{t_0,\upsilon_0,\upsilon}<\infty$ and $y_{t_0,\upsilon_0,\upsilon}>0$, with which on $[t_0+\upsilon,t_0+\upsilon_0]$,
\begin{IEEEeqnarray*}{rCl}
0\ge \Lambda'_t\ge -C_{t_0,\upsilon_0,\upsilon}
\quad\text{and}\quad 
0\le w(t,\Lambda_t\pm y)\le C_{t_0,\upsilon_0,\upsilon} y,\;\; y\in[0,y_{t_0,\upsilon_0,\upsilon}].
\end{IEEEeqnarray*}
\end{remark}

\smallskip

We turn to the proof of Proposition \ref{prop:LipschitzLambda1}. To simplify notation and without loss of generality, we take $t_0=0<\zeta$. Indeed, for any $0\le t_0<\zeta$, one can treat any physical solution $(\Lambda,w)$ from time $t_0$ on as if it starts at time $0$ with the initial condition $(\Lambda_{t_0-},w(t_0-,\cdot))$, thanks to \eqref{eq.sec3.sergey.growth.cond.t0.minus}, \eqref{eq.section3.sergey.w.back.Markov.tminus}. 
We also recall that $\Lambda_{0-}>r_0$ and $w(0-,\cdot)$ are fixed, that $\zeta>0$, and that $w^\lambda$ is constructed for given non-increasing right-continuous $\lambda\!:[-1,T]\to\RR$ with $\lambda_t\equiv\Lambda_{0-}$ on $[-1,0)$ by \eqref{eq.sec3.sergey.wLambda.def}. The next lemma bounds the Lipschitz constant of $w^{\lambda}$ near $\lambda$ through the weighted Lipschitz norm of $\lambda$.
\begin{lemma}
\label{lm:LambdaDerivativeInduction1}
Consider an $\upsilon_0>0$ and an absolutely continuous function $\lambda$ on $[0,\upsilon_0]$, and assume that there exist $C_1,C_2<\infty$, $\widehat{r}_0>r_0$ and $\varrho\in(0,1]$ such that $\lambda\ge\widehat{r}_0$ on $[0,\upsilon_0]$,
\begin{IEEEeqnarray*}{rCl}
&& 0\leq w^\lambda(0,\lambda_0-y)\le C_1y^\varrho,\quad y\in[0,\lambda_0-r_0), \\
&& 0\leq w^\lambda(0,\lambda_0+y)\le C_1y^\varrho,\quad y\in[0,\infty), \\
&& C_2:=\esssupmath{0\le t\le \upsilon_0}\;|\sqrt{t}\,\lambda'_t|<\infty.
\end{IEEEeqnarray*}
Then, for all $\varepsilon>0$, there exist $\upsilon,y_0>0$ such that 
\begin{IEEEeqnarray*}{rCl}
0\le\sqrt{t}\,\frac{w^\lambda(t,\lambda_t\pm y)}{y}&\leq&\frac{\varepsilon}{2}
\Big(1+\esssupmath{0\le t\le\upsilon_0}\;|\sqrt{t}\,\lambda'_t|^2\Big),\quad y\in(0,y_0],\quad t\in(0,\upsilon\wedge\upsilon_0].
\end{IEEEeqnarray*}
Moreover, the choice of $\upsilon,y_0$ depends only on the values of $\varepsilon,C_1,C_2,\widehat{r}_0,\varrho$. 
\end{lemma}

\noindent\textbf{Proof.} We fix an arbitrary $\varepsilon>0$ and $0<t\le\upsilon$, with the specific value of $\upsilon$ chosen later.

\medskip

\noindent\textit{Estimates for the solid phase.} Let us consider an arbitrary $y\in[0,\lambda_t-r_0)$ and bound $w^\lambda(t,\lambda_t-y)$. We abbreviate $\tau_{\lambda_{t-\cdot}}$, defined per \eqref{eq.sec3.sergey.tauLambda.tminus.def}, by $\tau$. From the definition \eqref{eq.sec3.sergey.wLambda.def} it is easy to see that $w(0-,R_t)=w^\lambda(0,R_t)$ on $\{\tau\wedge\tau_{r_0}>t\}$, $\PP^{\lambda_t-y}$-almost surely. Thus, thanks to  \eqref{eq.sec3.sergey.wLambda.def} again, 
\begin{equation}\label{eq:wsplit1}
0\leq w^\lambda(t,\lambda_t-y)
=\EE^{\lambda_t-y}[w^\lambda(0,R_t)\,\bone_{\{\tau\wedge\tau_{r_0}>t\}}]
\leq\EE^{\lambda_t-y}[w^\lambda(0,R_t)\,\bone_{\{\tau>t\}}].
\end{equation}
Next, we define the auxiliary boundary $\widehat{\lambda}$ by $\widehat\lambda_s:=\lambda_t+2C_2(t-s)/\sqrt{t}$, $s\in[0,t]$, which satisfies $\widehat{\lambda}_t=\lambda_t$ and $\widehat{\lambda}_0=\lambda_t+2C_2\sqrt{t}$. Moreover, since $\lambda'_s\ge-C_2/\sqrt{s}$ by the definition of $C_2$, 
\begin{IEEEeqnarray*}{rCl}
&& \widehat{\lambda}_s=\lambda_t+2C_2(t-s)/\sqrt{t}
=\lambda_t+2C_2(\sqrt{t}-s/\sqrt{t}) \\
&& \qquad\qquad\qquad\qquad\qquad\;\;\;\, \ge\lambda_t+2C_2(\sqrt{t}-\sqrt{s})
=\lambda_t+\int_s^t \frac{C_2}{\sqrt{s'}}\,\mathrm{d}s'
\ge\lambda_t+(\lambda_s-\lambda_t)=\lambda_s,\;\;\, s\in[0,t], \\   
&& w^\lambda(0,x)=w^\lambda(0,\lambda_{0}-(\lambda_{0}-x))
\le C_1(\lambda_{0}-x)^\varrho
\le C_1(\widehat{\lambda}_0-x)^\varrho
= C_1(\lambda_t+2C_2\sqrt{t}-x)^\varrho,\;\;\, x\in(0,\lambda_{0}].
\end{IEEEeqnarray*}

\smallskip

Using $\mathrm{d}R_s=R_s^{-1}\,\mathrm{d}s+\mathrm{d}B_s$ and the Girsanov theorem we bound
\begin{equation} \label{eq:Section3Est1}
\begin{split}
\frac{w^\lambda(t,\lambda_t-y)}{y}
\le \frac1y\,\EE^{\lambda_t-y}[w^\lambda(0,R_t)\,\bone_{\{\tau>t\}}]
\le\frac1y\,\EE^{\lambda_t-y}[C_1(\lambda_t+2C_2\sqrt{t}-R_t)^\varrho\,\bone_{\{\tau>t\}}] \\
\le\frac1y\,\EE^{\lambda_t-y}\big[C_1(y+2C_2\sqrt{t}-B_t)^\varrho\,\bone_{\{\max_{0\le s\le t}
(R_s-\widehat\lambda_{t-s})\le 0\}}\big] \\
\le\frac1y\,\EE\big[C_1(y+2C_2\sqrt{t}-B_t)^\varrho\,
\bone_{\{\max_{0\le s\le t} (B_s-2C_2s/\sqrt{t})\le y\}}\big] \\
=\frac1y\,\EE\Big[C_1(y-B_t)^\varrho\, e^{-2C_2B_t/\sqrt{t}-2C_2^2}\,\bone_{\{\max_{0\le s\le t} B_s\le y\}}\Big] \\
=\frac{C_1}{y} \int_0^y \int_{-\infty}^m (y-b)^\varrho\, e^{-2C_2b/\sqrt{t}-2C_2^2}
\,\frac{2(2m-b)}{t\sqrt{2\pi t}}\,e^{-\frac{(2m-b)^2}{2t}}\,\mathrm{d}b\,\mathrm{d}m \\
\le\frac{2C_1 y^\varrho}{\sqrt{t}}\left(\frac{1}{\sqrt{2\pi}}+2C_2\right)
+4t^{\frac{\varrho-1}{2}}C_1(1+2C_2^2),
\end{split}
\end{equation}
where, in the last inequality, we have applied the estimate \eqref{eq:GaussianMWestimate4} with $\mu=2C_2/\sqrt{t}$.

\medskip

\noindent\textit{Estimates for the liquid phase.} For the liquid phase, we proceed similarly, but now employing
\begin{IEEEeqnarray*}{rCl}
\widehat\lambda_s:=\lambda_t-2C_2(t-s)/\sqrt{t},\quad s\in[0,t].
\end{IEEEeqnarray*}
Notice that $\widehat\lambda_t=\lambda_t$ and that $\widehat\lambda_s\leq\lambda_s$, $s\in[0,t]$.
In addition, for $x\geq\lambda_{0}$,
\begin{IEEEeqnarray*}{rCl}
w^\lambda(0,x)&=&w^\lambda(0,\lambda_{0}+(x-\lambda_{0}))
\le C_1(x-\lambda_{0})^\varrho
\le C_1(x-\widehat{\lambda}_{0})^\varrho
= C_1(x-\lambda_t+2C_2\sqrt{t})^\varrho.
\end{IEEEeqnarray*}
Thus, for $y\ge0$, we obtain
\begin{equation}\label{eq:wliquidsplit1}
0\leq w^\lambda(t,\lambda_t+y)
=\EE^{\lambda_t+y}[w^\lambda(0,R_t)\,\bone_{\{\tau>t\}}]
\le\EE^{\lambda_t+y}\big[C_1(R_t+2C_2\sqrt{t}-\lambda_t)^\varrho\,\bone_{\{\tau>t\}}\big].
\end{equation}
For the latter, we observe that on $\{\tau>t\}$, it holds $\min_{0\le s\le t} R_s \ge \min_{0\le s\le t} \lambda_s \ge \widehat{r}_0$, $\PP^{\lambda_t+y}$-almost surely, and thus also 
\begin{IEEEeqnarray*}{rCl}
R_s\leq\lambda_t+y+B_s+\frac{s}{\widehat{r}_0},\quad s\in[0,t].
\end{IEEEeqnarray*}
Consequently, the Girsanov theorem yields
\begin{equation}\label{eq:Section3EstLiquid1}
\begin{split}
\frac1y\, w^\lambda(t,\lambda_t+y)
&\le\frac1y\,\EE^{\lambda_t+y}\big[C_1(R_t+2C_2\sqrt{t}-\lambda_t)^\varrho\,
\bone_{\{\max_{0\le s\le t} (\widehat\lambda_{t-s}-R_s)\le0\}}\big] \\
&\le\frac1y\,\EE\big[C_1(y+2C_2\sqrt{t}+t/\widehat{r}_0+B_t)^\varrho\,
\bone_{\{\max_{0\le s\le t} (-B_s-2C_2s/\sqrt{t}-s/\widehat{r}_0)\le y\}}\big] \\
&=\frac1y\,\EE\big[C_1(y+B_t)^\varrho\, e^{(2C_2/\sqrt{t}+1/\widehat{r}_0) B_t
-\frac12(2C_2/\sqrt{t}+1/\widehat{r}_0)^2 t}\,\bone_{\{\max_{0\le s\le t} (-B_s)\le y\}}\big] \\
&=\frac{C_1}{y}\int_0^y\int_{-\infty}^m (y-b)^\varrho\, e^{-(2C_2/\sqrt{t}+1/\widehat{r}_0)b-\frac12(2C_2/\sqrt{t}+1/\widehat{r}_0)^2t}\,
\frac{2(2m-b)}{t\sqrt{2\pi t}}e^{-\frac{(2m-b)^2}{2t}}\,\mathrm{d}b\,\mathrm{d}m \\
&\le \frac{2C_1y^\varrho}{\sqrt{t}}\bigg(\frac{1}{\sqrt{2\pi}}+2C_2+\frac{\sqrt{t}}{\widehat{r}_0}\bigg)+4t^{\frac{\varrho-1}{2}}C_1\bigg(1+2\Big(C_2+\frac{\sqrt{t}}{2\widehat{r}_0}\Big)^2\bigg),
\end{split}
\end{equation}
where, in the last inequality, we have used the estimate \eqref{eq:GaussianMWestimate4} with $\mu=2C_2/\sqrt{t}+1/\widehat{r}_0$.

\medskip

\noindent\textit{Choice of the parameters and conclusion.}  All in all, we have shown that
\begin{IEEEeqnarray*}{rCl}
\sqrt{t}\,\frac{w^\lambda(t,\lambda_t\pm y)}{y}\le C(y,t,C_1,C_2,\widehat{r}_0,\varrho),\quad y\in(0,y_0],\quad t\in(0,\upsilon],
\end{IEEEeqnarray*}
where
\begin{IEEEeqnarray*}{rCl}
C(y,t,C_1,C_2,\widehat{r}_0,\varrho)
:=2C_1y^\varrho\bigg(\frac{1}{\sqrt{2\pi}}+2C_2+\frac{\sqrt{t}}{\widehat{r}_0}\bigg)
+4t^{\frac{\varrho}{2}}C_1\bigg(1+2\Big(C_2+\frac{\sqrt{t}}{2\widehat{r}_0}\Big)^2\bigg).
\end{IEEEeqnarray*}
Since $C_2\le(1+C_2^2)/2$, we can pick small enough $\upsilon,y_0>0$, depending on $\varepsilon,C_1,C_2,\widehat{r}_0,\varrho$, to ensure
\begin{IEEEeqnarray*}{rCl}
\qquad\qquad\quad\;  
C(y_0,\upsilon,C_1,C_2,\widehat{r}_0,\varrho):=
\sup_{y\in(0,y_0],\,t\in(0,\upsilon]} C(y,t,C_1,C_2,\widehat{r}_0,\varrho)
\le\frac\varepsilon2(1+C_2^2). \qquad\qquad\quad\;   \qed
\end{IEEEeqnarray*}

\smallskip

\noindent\textbf{Proof of Proposition \ref{prop:LipschitzLambda1}.} Recall that we assume $t_0=0<\zeta$ without loss of generality. Consider a sequence $(\Lambda^n)_{n\in\NN_0}$ of non-increasing continuous functions $[0,\infty)\to\RR$ defined recursively per
\begin{IEEEeqnarray*}{rCl}
\Lambda^0\equiv\Lambda_0,\quad\Lambda^{n+1}=\Gamma[\Lambda^n],
\end{IEEEeqnarray*}
where $\Gamma$ is given by \eqref{eq.section3.sergey.Gamma.def}. Every $\Lambda^n$ is non-increasing and continuous since the mapping $\Gamma$ preserves these properties, by its definition. Concurrently, we keep track of the corresponding functions
\begin{IEEEeqnarray*}{rCl}
w_n(t,x):=w^{\Lambda^n}(t,x)=\EE^x[w(0-,R_t)\,\bone_{\{\tau_{\Lambda_{t-\cdot}^n}\wedge\tau_{r_0}>t\}}],\quad x>0,\quad t\in[0,\zeta^{\Lambda^n}),
\end{IEEEeqnarray*}
where we set $\Lambda^n_s=\Lambda_{0-}$ for $s\in[-1,0)$ in the definition of $\tau_{\Lambda_{t-\cdot}^n}$.

\smallskip

\noindent\textbf{Step 1.} Let $\widehat{r}_0=(\Lambda_0+r_0)/2$, fix an arbitrary $\varepsilon\in(0,\widehat{r}_0^2/(2(\Lambda_0)^2))$, and define $C_\varepsilon\in(0,\frac12)$ by
\begin{IEEEeqnarray*}{rCl}
\frac{\widehat\varepsilon}{2}(1+C_\varepsilon^2)=C_\varepsilon,\quad\text{where}\quad \widehat{\varepsilon}:=\frac{(\Lambda_0)^2}{\widehat{r}_0^2}\varepsilon\in\Big(0,\frac{1}{2}\Big).
\end{IEEEeqnarray*}
Specifically, $C_\varepsilon:=\widehat\varepsilon/(1+\sqrt{1-\widehat\varepsilon^2})$. Recall also $\upsilon=\upsilon(\varepsilon,C_1,C_2:=C_\varepsilon,\widehat{r}_0,\varrho)$ from Lemma \ref{lm:LambdaDerivativeInduction1}, where $C_1$, $\varrho$ are determined by the assumption in Proposition \ref{prop:LipschitzLambda1}. Finally, take
\begin{align*}
&\upsilon_0=\upsilon\wedge\frac{((\Lambda_0)^3-\widehat{r}_0^3)^2}{36C_\varepsilon^2\Lambda_0^4}.
\end{align*}

\smallskip

We are going to show by induction that, for all $n\in\NN_0$,
\begin{equation} \label{eq:InductionGrowthEst1}
\Lambda^n\in C([0,\upsilon_0))\cap C^1((0,\upsilon_0)),\quad
\Lambda^n_t\ge\widehat{r}_0,\;\; t\in[0,\upsilon_0]\quad\text{and}\quad
\sup_{0<t<\upsilon_0} |\sqrt{t}\,(\Lambda_t^n)'|\le C_\varepsilon.
\end{equation}
The claim is trivial for $n=0$, so we now assume that it holds for some $n\ge0$ and aim to prove it for $n+1$. Arguing as in the proof of Proposition \ref{prop:InteriorSmoothness} we see that $w_n$ is of class $C^\infty$ and satisfies 
\begin{IEEEeqnarray*}{rCl}
\partial_tw_n(t,x)=\frac12\partial_{xx}w_n(t,x)+\frac1x\partial_xw_n(t,x)    
\end{IEEEeqnarray*}
pointwise on $\stackrel{\circ}{D_n}\,:=\{(t,x)\in(0,\zeta^{\Lambda^n})\times(r_0,\infty):\,x>\Lambda^n_{t-}\;\text{or}\;x<\Lambda^n_t\}$. Thanks to the physicality of~$\Lambda$, we obtain by induction that $\Lambda^n_0=\Lambda_0$, which in turn implies $w_n(0,\cdot)=w^{\Lambda^n}(0,\cdot)=w(0,\cdot)$. Thus, the assumptions of Lemma \ref{lm:LambdaDerivativeInduction1} are satisfied, with $\lambda=\Lambda^n$, $C_2=C_{\varepsilon}$, and with $C_1$, $\varrho$ determined by the assumption in Proposition \ref{prop:LipschitzLambda1}. Let $\upsilon,y_0>0$ be as in the statement of Lemma \ref{lm:LambdaDerivativeInduction1}. Then,
\begin{equation}\label{eq:InductionLipEst1}
0\le\sqrt{t}\,\frac{w_n(t,\Lambda_t^n\pm y)}{y}
\le\frac{\varepsilon}{2}\Big(1+\sup_{0<t\le\upsilon_0}
\;|\sqrt{t}\,(\Lambda_t^n)'|^2\Big)
\le\frac\varepsilon2(1+C_\varepsilon^2)\le C_\varepsilon,\; y\in(0,y_0],\; t\in(0,\upsilon\wedge\upsilon_0=\upsilon_0).
\end{equation}

\smallskip

\noindent\textbf{Step 2.} The goal of this step is to estimate $\partial_xw_n$ near the boundary $\Lambda^n$. For any $t\in(0,\upsilon_0)$ and $y\in(0,\sqrt{t/2}\wedge (y_0/(2C_\varepsilon+2))\wedge((\widehat{r}_0-r_0)/2))$, we argue as in the proof of Proposition \ref{prop:GradientEstimateForw1} to get
\begin{IEEEeqnarray*}{rCl}
|\partial_xw_n(t,\Lambda_t^n-y)|\le
\frac{C}{y}\sup_{(s,x)\in(t-y^2,t)\times(\Lambda_t^n-2y,\Lambda_t^n)} |w_n(s,x)|.    
\end{IEEEeqnarray*}
For $(s,x)\in(t-y^2,t)\times(\Lambda_t^n-2y,\Lambda_t^n)\subset\,\stackrel{\circ}{D_n}$, it holds $s>t-y^2\ge \frac{t}{2}$ and 
\begin{IEEEeqnarray*}{rCl}
0<\Lambda_s^n-x=(\Lambda_s^n-\Lambda_t^n)+(\Lambda_t^n-x)
\le 2C_\varepsilon(\sqrt{t}-\sqrt{s})+2y
\le 2C_\varepsilon\sqrt{t-s}+2y
\le (2C_\varepsilon+2)y\le y_0.
\end{IEEEeqnarray*}
Thanks to the inequality \eqref{eq:InductionLipEst1} this implies
\begin{IEEEeqnarray*}{rCl}
0\leq w_n(s,x)=w_n(s,\Lambda_s^n-(\Lambda_s^n-x))
\le \frac{C_\varepsilon}{\sqrt{s}}(\Lambda_s^n-x)
\le \frac{4C_\varepsilon}{\sqrt{t}}(C_\varepsilon+1)y.   
\end{IEEEeqnarray*}
As a result, for $t\in(0,\upsilon_0)$ and $y\in(0,\sqrt{t/2}\wedge (y_0/(2C_\varepsilon+2))\wedge((\widehat{r}_0-r_0)/2))$, we have
\begin{IEEEeqnarray*}{rCl}
|\partial_xw_n(t,\Lambda_t^n-y)|\leq \frac{4CC_\varepsilon}{\sqrt{t}}(C_\varepsilon+1).    
\end{IEEEeqnarray*}
This, together with the continuity of $\partial_xw_n$ on $\stackrel{\circ}{D_n}$, proves that $\partial_xw_n$ is bounded on
\begin{IEEEeqnarray*}{rCl}
\{(t,x)\in[\upsilon_1,\upsilon_0)\times[r_1,\infty):\,x<\Lambda_t^n\},
\end{IEEEeqnarray*}
for any $\upsilon_1\in(0,\upsilon_0)$ and $r_1\in(r_0,\widehat{r}_0)$.

\medskip

Consider $\partial_xw_n$ now above the boundary $\Lambda^n$. For $t\in(0,\upsilon_0)$ and $y\in(0,(C_\varepsilon\sqrt{t})\wedge((2y_0)/(1+C^{-1}_\varepsilon)))$,
\begin{IEEEeqnarray*}{rCl}
\Lambda_t^n+\frac{y}{2}-\Lambda_{t-\frac{y^2}{16C_\varepsilon^2}}^n
\ge\frac{y}{2}-2C_\varepsilon\Bigg(\sqrt{t}-\sqrt{t-\frac{y^2}{16C_\varepsilon^2}}\Bigg)
\ge\frac{y}{2}-\frac{y}{2}=0.
\end{IEEEeqnarray*}
Hence,  $(t-y^2/(16C^2_\varepsilon),t)\times(\Lambda_t^n+y/2,\Lambda_t^n+y(1+C^{-1}_\varepsilon)/2)
\subset\,\stackrel{\circ}{D_n}$. The proof of Proposition \ref{prop:GradientEstimateForw1} yields
\begin{IEEEeqnarray*}{rCl}
|\partial_xw_n(t,\Lambda_t^n+y)|\le\frac{2C}{y}\sup_{(s,x)\in(t-\frac{y^2}{16C^2_\varepsilon},t)\times(\Lambda_t^n+\frac{y}{2},\Lambda_t^n+\frac{y}{2}(1+C^{-1}_\varepsilon))}|w_n(s,x)|.    
\end{IEEEeqnarray*}
For $(s,x)\in(t-y^2/(16C^2_\varepsilon),t)\times(\Lambda_t^n+y/2,\Lambda_t^n+y(1+C^{-1}_\varepsilon)/2)$, we have $s\ge t/2$ and   
\begin{IEEEeqnarray*}{rCl}
0<x-\Lambda_s^n\leq\Lambda_t^n+\frac{y}{2}(1+C^{-1}_\varepsilon)-\Lambda_t^n
=\frac{y}{2}(1+C^{-1}_\varepsilon)\le y_0.    
\end{IEEEeqnarray*}
This and the inequality \eqref{eq:InductionLipEst1} result in
\begin{IEEEeqnarray*}{rCl}
0\leq w_n(s,x)=w_n(s,\Lambda_s^n+(x-\Lambda_s^n))
\le\frac{1+C^{-1}_\varepsilon}{2\sqrt{s}}C_\varepsilon y
\le \frac{C_\varepsilon+1}{\sqrt{t}} y.  
\end{IEEEeqnarray*}
Consequently, for $t\in(0,\upsilon_0)$ and $y\in(0,(C_\varepsilon\sqrt{t})\wedge ((2y_0)/(1+C^{-1}_\varepsilon)))$, it holds that
\begin{IEEEeqnarray*}{rCl}
|\partial_xw_n(t,\Lambda_t^n+y)|\le \frac{2C}{\sqrt{t}}(C_\varepsilon+1).
\end{IEEEeqnarray*}
This and the continuity of $\partial_xw_n$ on $\stackrel{\circ}{D_n}$ show that $\partial_xw_n$ is bounded on
\begin{IEEEeqnarray*}{rCl}
\{(t,x)\in[\upsilon_1,\upsilon_0)\times(0,r_2]:\,x>\Lambda_t^n\},
\end{IEEEeqnarray*}
for any $\upsilon_1\in(0,\upsilon_0)$ and $r_2\in(\Lambda_0,\infty)$.

\medskip

\noindent\textbf{Step 3.} Our next goal is to establish the continuity of $\partial_xw_n$ up to the boundary $\Lambda^n$. To this end, let $L\in(0,(\widehat{r}_0-r_0)/2)$ and $\upsilon_1\in(0,\upsilon_0)$, and consider the auxiliary function
\begin{IEEEeqnarray*}{rCl}
\widecheck{w}_n(s,y):=w_n(\upsilon_1+s,\Lambda_{\upsilon_1+s}^n-y)
-\frac{y}{L}w_n(\upsilon_1+s,\Lambda_{\upsilon_1+s}^n-L),\quad (s,y)\in[0,\upsilon_0-\upsilon_1)\times[0,L].
\end{IEEEeqnarray*}
It satisfies
\begin{equation} \label{eq:AuxiliaryPDEtildew2}
\begin{cases}
\partial_s\widecheck{w}_n(s,y)=\frac12\partial_{yy}\widecheck{w}_n(s,y)+f(s,y)-\frac{y}{L}g'(s), &\quad(s,y)\in(0,\upsilon_0-\upsilon_1)\times(0,L), \\
\widecheck{w}_n(s,0)=\widecheck{w}_n(s,L)=0,&\quad s\in[0,\upsilon_0-\upsilon_1), \\
\widecheck{w}_n(0,y)=h(y)-\frac{y}{L}g(0),&\quad y\in(0,L),
\end{cases}
\end{equation}
where
\begin{equation*}
\begin{split}
& f(s,y):=\Big(\frac{1}{\Lambda_{\upsilon_1+s}^n-y}+(\Lambda_{\upsilon_1+s}^n)'\Big)
\,\partial_xw_n(\upsilon_1+s,\Lambda_{\upsilon_1+s}^n-y), \\ 
& g(s):=w_n(\upsilon_1+s,\Lambda_{\upsilon_1+s}^n-L), \\
& h(y):=w_n(\upsilon_1,\Lambda_{\upsilon_1}^n-y).
\end{split}
\end{equation*}

\smallskip

It follows from the assumption $|(\Lambda_t^n)'|\le C_\varepsilon/\sqrt{t}$, $t\in(0,\upsilon_0)$, from the result of Step 2 and from the interior smoothness of $w_n$ that $f$ is bounded on $(0,\upsilon_0-\upsilon_1)\times(0,L)$ and that $g$ is bounded and absolutely continuous on $[0,\upsilon_0-\upsilon_1)$, with its weak derivative
\begin{IEEEeqnarray*}{rCl}
g'(s)=\partial_tw_n(\upsilon_1+s,\Lambda_{\upsilon_1+s}^n-L)+\partial_xw_n(\upsilon_1+s,\Lambda_{\upsilon_1+s}^n-L)\,(\Lambda_{\upsilon_1+s}^n)'
\end{IEEEeqnarray*}
being bounded by a constant for almost every $s\in[0,\upsilon_0-\upsilon_1)$. Consequently, the auxiliary function $\widecheck w_n$ solves the problem \eqref{eq:AuxiliaryPDEtildew2} weakly in the sense of \cite[Section 7.1]{Ev}. On the other hand, the function $y\mapsto h(y)-\frac{y}{L}g(0)$ is in $\mathring H^1(0,L)$ by the result of Step 2 and the interior smoothness of $w_n$. Thus, the problem \eqref{eq:AuxiliaryPDEtildew2} admits a unique weak solution in the function class
\begin{IEEEeqnarray*}{rCl}
\big\{\widecheck{w}\in L^2(0,\upsilon_0-\upsilon_1;\mathring H^1(0,L)):\,\partial_s\widecheck{w}\in L^2(0,\upsilon_0-\upsilon_1;H^{-1}(0,L))\big\}  
\end{IEEEeqnarray*}
by \cite[Theorem 7.1.4]{Ev}. Due to the boundedness of $\widecheck{w}_n$ and $\partial_x\widecheck{w}_n$, the function $\widecheck{w}_n$ does belong to the above function class.~Hence, we identify $\widecheck{w}_n$ as the unique weak solution of the problem~\eqref{eq:AuxiliaryPDEtildew2}. 

\medskip

Next, we apply the energy estimates for the heat equation with the Dirichlet boundary condition (see, e.g., \cite[Theorem 7.1.5]{Ev}) to find that $\widecheck{w}_n\in L^2(0,\upsilon_0-\upsilon_1;H^2(0,L))$. In particular, there is a sequence $t_m\downarrow0$ such that each $\widecheck{w}_n(t_m,\cdot)$ is in $H^2(0,L)$. Since $H^2(0,L)\subset W_q^{2-\frac2q}(0,L)$ for every $q\in(3,6)$ (see \cite[Section II.2]{Lady} for the definition of the Sobolev spaces $W_q^l$ for non-integral~$l$), we can view $\widecheck{w}_n(t_m,\cdot)$ as the initial condition in the problem for $\widecheck{w}_n$, so that $\widecheck{w}_n$, restricted to $[t_m,\upsilon_0-\upsilon_1)\times[0,L]$, is in $W_q^{1,2}([t_m,\upsilon_0-\upsilon_1)\times[0,L])$ by \cite[Theorem IV.9.1]{Lady}. The Sobolev embedding theorem (see \cite[Lemma II.3.3]{Lady}) then implies that $\partial_y\widecheck{w}_n$ is $\frac\lambda2$-Hölder continuous on $[t_m,\upsilon_0-\upsilon_1)\times[0,L]$ for every $\lambda\in(0,1-3/q)$. Thus, letting $m\to\infty$ we deduce that $\partial_y\widecheck{w}_n$ has a version that is continuous on $(0,\upsilon_0-\upsilon_1)\times[0,L]$. Recalling the relation 
\begin{IEEEeqnarray*}{rCl}
\partial_xw_n(t,x)=-\partial_y\widecheck{w}_n(t-\upsilon_1,\Lambda_t^n-x)-\frac{1}{L}w_n(t,\Lambda_t^n-L)     
\end{IEEEeqnarray*}
and that the choice of $\upsilon_1\in(0,\upsilon_0)$ is arbitrary, we get the continuity of $\partial_xw_n$ on 
\begin{IEEEeqnarray*}{rCl}
\{(t,x)\in(0,\upsilon_0)\times(r_0,\infty):\,x\le\Lambda_t^n\}.
\end{IEEEeqnarray*}
The continuity of $\partial_xw_n$ on 
\begin{IEEEeqnarray*}{rCl}
\{(t,x)\in(0,\upsilon_0)\times(0,\infty):\,x\ge\Lambda_t^n\}
\end{IEEEeqnarray*}
is proved in the same way by considering $\widecheck{w}_n(s,y):=w_n(\upsilon_1+s,\Lambda_{\upsilon_1+s}^n+y)-\frac{y}{L}w_n(\upsilon_1+s,\Lambda^n_{\upsilon_1+s}+L)$.

\medskip

The Stefan growth condition stated in Corollary \ref{cor:BoundaryCondition2}, together with \eqref{eq:InductionLipEst1}, now gives
\begin{equation*}
\frac{(\Lambda^{n+1}_t)^3}{3}
\ge\frac{(\Lambda^{n+1}_0)^3}{3}-C_\varepsilon\int_0^t\frac{(\Lambda^n_s)^2}{\sqrt{s}}\,\mathrm{d}s 
\ge\frac{(\Lambda_0)^3}{3}-2C_\varepsilon(\Lambda_0)^2\sqrt{t},\quad t\in(0,\upsilon_0),
\end{equation*}
so $\Lambda^{n+1}_t\ge\widehat{r}_0$, $t\in(0,\upsilon_0)$ thanks to the choice of $\upsilon_0$. This, Corollary \ref{cor:BoundaryCondition2} and the continuity of $\partial_x w_n$ show the continuity of $(\Lambda^{n+1})'$ on $(0,\upsilon_0)$. Moreover, combining Corollary \ref{cor:BoundaryCondition2} and \eqref{eq:InductionLipEst1} we find
\begin{equation*}
\sup_{0<t<\upsilon_0} |\sqrt{t}\,(\Lambda_t^{n+1})'|
\le\frac{(\Lambda_0)^2}{\widehat{r}_0^2}\frac{\varepsilon}{2}(1+C_\varepsilon^2)=C_\varepsilon.	
\end{equation*}
The proof of \eqref{eq:InductionGrowthEst1} is complete.

\medskip

\noindent\textbf{Step 4.} 
In view of \eqref{eq:InductionGrowthEst1} and Remark \ref{rmk dep}, we can apply the estimates of Propositions \ref{prop multi}, \ref{prop 2nd phase} to deduce the existence of a $c<1$, independent of $n\in\NN_0$, such that
\begin{IEEEeqnarray*}{rCl}
\max_{0\le t\le\upsilon_0} \big|(\Lambda^{n+1}_t)^3-(\Lambda^n_t)^3\big|
=\max_{0\le t\le\upsilon_0} \big|\Gamma[\Lambda^n]_t^3-\Gamma[\Lambda^{n-1}]_t^3\big|
\le c\max_{0\le t\le\upsilon_0} \big|(\Lambda^n_t)^3-(\Lambda^{n-1}_t)^3\big|,
\end{IEEEeqnarray*}
upon shrinking $\upsilon_0>0$ if necessary. This and \eqref{eq:InductionGrowthEst1} show that $(\Lambda^n)_{n\in\NN_0}$ tends uniformly on $[0,\upsilon_0]$ to a continuous limit $\widetilde\Lambda$, which is a fixed point of $\Gamma$. (The proof of the latter is standard and follows from the crossing property of the Bessel process.) As $\Lambda$ is also a fixed point of $\Gamma$, we conclude that $\widetilde\Lambda=\Lambda$, by again applying Remark \ref{rmk dep}, upon shrinking $\upsilon_0>0$ further if necessary.

\medskip

Next, we notice that \eqref{eq:InductionGrowthEst1} reveals the non-decreasing nature of every $t\mapsto\Lambda_t^n+2C_\varepsilon\sqrt{t}$ on $[0,\upsilon_0]$.
Passing to the limit $n\to\infty$ we deduce that $t\mapsto\Lambda_t+2C_\varepsilon\sqrt{t}$ is non-decreasing on $[0,\upsilon_0]$, i.e.:
\begin{IEEEeqnarray*}{rCl}
0\ge\Lambda_{t_2}-\Lambda_{t_1}\ge-2C_\varepsilon(\sqrt{t_2}-\sqrt{t_1}),\quad 
0\le t_1<t_2\le\upsilon_0.
\end{IEEEeqnarray*}
This proves the Lipschitz continuity of $\Lambda$ on $(0,\upsilon_0]$, as well as
\begin{IEEEeqnarray*}{rCl}
\esssupmath{0\le t\le\upsilon_0}\,|\sqrt{t}\,\Lambda'_t|\le C_\varepsilon<\infty.
\end{IEEEeqnarray*}
Having obtained the absolute continuity of $\Lambda$, we can rely on Lemma \ref{lm:LambdaDerivativeInduction1} to infer 
\begin{equation} \label{eq:InductionLipEst3}
0\leq\sqrt{t}\,\frac{w(t,\Lambda_t\pm y)}{y}
\le\frac{\varepsilon}{2}\Big(1+\esssupmath{0\le t\le\upsilon_0}\,|\sqrt{t}\,\Lambda'_t|^2\Big)
\le\frac\varepsilon2(1+C_\varepsilon^2)
\le C_\varepsilon,\quad y\in(0,y_0],\quad t\in(0,\upsilon\wedge\upsilon_0].
\end{equation}
Recalling that $\Gamma[\Lambda]=\Lambda$ we repeat the arguments in Steps 2 and 3 (observe that the latter do not require the continuity of $\Lambda'$) to deduce the continuity of $\partial_xw$ on
\begin{IEEEeqnarray*}{rCl}
\{(t,x)\in(0,\upsilon\wedge\upsilon_0)\times(r_0,\infty):\,x\le\Lambda_t\}\quad\text{and}\quad
\{(t,x)\in(0,\upsilon\wedge\upsilon_0)\times(0,\infty):\,x\ge\Lambda_t\},
\end{IEEEeqnarray*}
as well as the continuity of $\Lambda'$ on $(0,\upsilon\wedge\upsilon_0)$ and the identity
\begin{IEEEeqnarray*}{rCl}
\qquad\qquad\qquad\quad\;\;
\Lambda'_t=\frac12\partial_xw(t,x-)|_{x=\Lambda_t}-\frac12\partial_xw(t,x+)|_{x=\Lambda_t},\quad t\in(0,\upsilon\wedge\upsilon_0).
\qquad\qquad\qquad\quad\;\; \qed    
\end{IEEEeqnarray*}

\smallskip

By combining Proposition \ref{prop:LipschitzLambda1} with the results of Subsection \ref{sec:Holder} we obtain the next proposition. 

\begin{proposition}
\label{prop:BoundaryRegularityC1}
Fix an arbitrary $t_0\in[0,\zeta)$ and assume that $w(t_0-,\Lambda_{t_0-}-\cdot)$ is monotone in a right neighborhood of any point in $[0,\Lambda_{t_0-}-r_0)$. Then, there exists an $\upsilon_0>0$ such that $\Lambda\in C^1((t_0,t_0+\upsilon_0))$
and $w$, $\partial_xw$ are continuous on
\begin{IEEEeqnarray*}{rCl}
\{(t,x)\in(t_0,t_0+\upsilon_0)\times(r_0,\infty):\,x\leq\Lambda_t\}\quad\text{and}\quad
\{(t,x)\in(t_0,t_0+\upsilon_0)\times(0,\infty):\,x\geq\Lambda_t\}.
\end{IEEEeqnarray*}
Moreover, the Stefan growth condition holds:
\begin{IEEEeqnarray*}{rCl}
\Lambda'_t=\frac12\partial_xw(t,x-)|_{x=\Lambda_t}-\frac12\partial_xw(t,x+)|_{x=\Lambda_t},
\quad t\in(t_0,t_0+\upsilon_0).
\end{IEEEeqnarray*}
\end{proposition}

\subsection{Auxiliary Results for the Uniqueness Argument}

In the statements of this subsection, we denote by $(\Lambda,w)$ an arbitrary (fixed) physical solution with the initial data $(\Lambda_{0-},w(0-,\cdot))$.

\begin{lemma}\label{lm:wStrictlyPositive}
$w(t,x)>0$, $x\in(r_0,\Lambda_t)$, $t\in(0,\zeta)$ unless $w(0-,x)=0$ for almost every $x\in(r_0,\Lambda_0)$.    
\end{lemma}

\noindent\textbf{Proof.}~By its monotonicity-changing property $w(0-,\cdot)$ is continuous almost everywhere in $(r_0,\Lambda_{0-})$. Unless $w(0-,x)=0$ for almost every $x\in(r_0,\Lambda_0)$, there exist $\varepsilon>0$ and $(a_\varepsilon,b_\varepsilon)\subset(r_0,\Lambda_0)$ such that $w(0,x)=w(0-,x)\ge\varepsilon$, $x\in(a_\varepsilon,b_\varepsilon)$. As a result, for any $t\in(0,\zeta)$ and $x\in(r_0,\Lambda_t)$,
\begin{IEEEeqnarray*}{rCl}
\qquad\quad\;\;\; w(t,x)=\EE^x[w(0,R_t)\,\bone_{\{\tau_{\Lambda_{t-\cdot}}\wedge\tau_{r_0}>t\}}]\ge\varepsilon\PP^x\big(R_t\in (a_\varepsilon,b_\varepsilon),\,\tau_{\Lambda_{t-\cdot}}\wedge\tau_{r_0}>t\big)>0. \qquad\quad\;\;\; \qed  
\end{IEEEeqnarray*}

\smallskip

The following result shows that the boundary $\Lambda$ decreases at a positive rate on its intervals of $C^1$-smoothness, due to the solid phase. This is used in Appendix \ref{se:4}.
\begin{proposition}
\label{prop:partialxwnegativeRegular}
Fix an arbitrary $t_0\in[0,\zeta)$ and assume that $w(t_0-,\Lambda_{t_0-}-\cdot)$ is monotone in a right neighborhood of any point in $[0,\Lambda_{t_0-}-r_0)$. Let $\upsilon_0$ be the constant of Proposition \ref{prop:BoundaryRegularityC1}. Then, 
\begin{IEEEeqnarray*}{rCl}
\Lambda'_t\le\frac12\partial_xw(t,x-)|_{x=\Lambda_t}<0,\quad t\in(t_0,t_0+\upsilon_0)
\end{IEEEeqnarray*}
unless $w(0-,x)=0$ for almost every $x\in(r_0,\Lambda_0)$.
\end{proposition}

\noindent\textbf{Proof.} The first inequality follows from the Stefan growth condition
\begin{IEEEeqnarray*}{rCl}
\Lambda'_t=\frac12\partial_xw(t,x-)|_{x=\Lambda_t}-\frac12\partial_xw(t,x+)|_{x=\Lambda_t},\quad t\in(t_0,t_0+\upsilon_0)   
\end{IEEEeqnarray*}
and from the fact that $\partial_xw(t,x+)|_{x=\Lambda_t}\ge0$ (by $w(t,x+)|_{x=\Lambda_t}=0$ and $w(t,\Lambda_t+y)\ge0$, $y>0$). To get the second inequality, we fix arbitrary $\upsilon_1\in(0,\upsilon_0)$, $L\in(0,(\Lambda_{t_0+\upsilon_0}-r_0)/2)$ and apply the classical Hopf lemma for parabolic equations (see, e.g., \cite[Theorem 2]{Friedman1958RemarksOT}) to 
\begin{IEEEeqnarray*}{rCl}
\widecheck{w}(s,y):=w(t_0+\upsilon_1+s,\Lambda_{t_0+\upsilon_1+s}-y),    
\end{IEEEeqnarray*}
which is a non-negative classical ($C^{1,2}$) solution of the PDE
\begin{IEEEeqnarray*}{rCl}
\partial_s\widecheck{w}(s,y)=\frac12\partial_{yy}\widecheck{w}(s,y)-\Big(\frac{1}{\Lambda_{t_0+\upsilon_1+s}-y}+\Lambda'_{t_0+\upsilon_1+s}\Big)\partial_y\widecheck{w}(s,y),\quad (s,y)\in(0,\upsilon_0-\upsilon_1)\times(0,L).    
\end{IEEEeqnarray*}
The coefficients of the PDE are continuous on $(0,\upsilon_0-\upsilon_1)\times[0,L]$. The domain is rectangular, so that the interior ball condition trivially holds for any lateral boundary point $(s,0)\in(0,\upsilon_0-\upsilon_1)\times\{0\}$. Moreover, $\widecheck{w}(s,0+)=0$ and $\widecheck{w}$ is strictly positive in the interior $(0,\upsilon_0-\upsilon_1)\times(0,L)$ according to Lemma \ref{lm:wStrictlyPositive}. As a result, all the conditions in \cite[Theorem 2]{Friedman1958RemarksOT} apply to the function $\widecheck{w}$ and its domain $(0,\upsilon_0-\upsilon_1)\times(0,L)$, which then yields $\partial_y\widecheck{w}(s,0+)>0$, and thus $\partial_xw(t,x-)|_{x=\Lambda_t}=-\partial_y\widecheck{w}(t-t_0-\upsilon_1,0+)<0$, $t\in(t_0+\upsilon_1,t_0+\upsilon_0)$. It remains to recall that $\upsilon_1\in(0,\upsilon_0)$ is arbitrary. \qed

\begin{lemma}\label{lm:wSolidDecay}
There exists a constant $\overline C<\infty$ such that $w(t,x)\le\overline{C}/x$, $x>\Lambda_t$, $t\in[0,\zeta)$.   
\end{lemma}

\noindent\textbf{Proof.} By Assumption \ref{intro: main ass}, we can find a constant $\overline C<\infty$ such that $w(0-,x)\le\overline{C}/x$, $x>0$. Given $x>\Lambda_t$, $t\in[0,\zeta)$, under $\PP^x$, $1/R$ is a non-negative local martingale, hence a supermartingale. Consequently,
\begin{IEEEeqnarray*}{rCl}
\qquad\qquad\qquad\qquad\qquad\quad\;\;
w(t,x)\leq\EE^x[w(0-,R_t)]\le\EE^x\bigg[\frac{\overline C}{R_t}\bigg]
\le\frac{\overline C}{x}.\qquad\qquad\qquad\qquad\qquad\quad\;\; \qed
\end{IEEEeqnarray*}

\bigskip\bigskip

\bibliographystyle{amsalpha}
\bibliography{Main}

\providecommand{\bysame}{\leavevmode\hbox to3em{\hrulefill}\thinspace}
\providecommand{\MR}{\relax\ifhmode\unskip\space\fi MR }
\providecommand{\MRhref}[2]{%
  \href{http://www.ams.org/mathscinet-getitem?mr=#1}{#2}
}
\providecommand{\href}[2]{#2}
\begin{thebibliography}{DIRT15}

\bibitem[AF88]{AF}
S.~B. Angenent and B.~Fiedler, \emph{The dynamics of rotating waves in scalar
  reaction diffusion equations}, Trans. Amer. Math. Soc. \textbf{307} (1988),
  no.~2, 545--568. \MR{940217}

\bibitem[BS22]{BaSh}
G.~Baker and M.~Shkolnikov, \emph{Zero kinetic undercooling limit in the
  supercooled {S}tefan problem}, Ann. Inst. Henri Poincar\'{e} Probab. Stat.
  \textbf{58} (2022), no.~2, 861--871. \MR{4421610}

\bibitem[CH82]{ChowHale}
S.~N. Chow and J.~K. Hale, \emph{Methods of bifurcation theory}, Grundlehren
  der Mathematischen Wissenschaften [Fundamental Principles of Mathematical
  Science], vol. 251, Springer-Verlag, New York-Berlin, 1982.

\bibitem[DIRT15]{DIRT2}
F.~Delarue, J.~Inglis, S.~Rubenthaler, and E.~Tanr\'{e}, \emph{Particle systems
  with a singular mean-field self-excitation. {A}pplication to neuronal
  networks}, Stochastic Process. Appl. \textbf{125} (2015), no.~6, 2451--2492.
  \MR{3322871}

\bibitem[DNS22]{dns}
F.~Delarue, S.~Nadtochiy, and M.~Shkolnikov, \emph{Global solution to
  super-cooled {S}tefan problem with blow-ups: Regularity and uniqueness},
  Probab. Math. Phys. \textbf{3} (2022), no.~1, 171--213.

\bibitem[EPS03]{eps}
J.~Escher, J.~Pr\"{u}ss, and G.~Simonett, \emph{Analytic solutions for a
  {S}tefan problem with {G}ibbs-{T}homson correction}, J. Reine Angew. Math.
  \textbf{563} (2003), 1--52. \MR{2009238}

\bibitem[Eva10]{Ev}
L.~C. Evans, \emph{Partial differential equations}, second ed., Graduate
  Studies in Mathematics, vol.~19, American Mathematical Society, Providence,
  RI, 2010. \MR{2597943}

\bibitem[FR91]{FrRe}
A.~Friedman and F.~Reitich, \emph{The {S}tefan problem with small surface
  tension}, Trans. Am. Math. Soc. \textbf{328} (1991), no.~2, 465--515.

\bibitem[Fri58]{Friedman1958RemarksOT}
A.~Friedman, \emph{Remarks on the maximum principle for parabolic equations and
  its applications.}, Pac. J. Math. \textbf{8} (1958), 201--211.

\bibitem[Gli10]{Gli}
M.~E. Glicksman, \emph{Principles of solidification: an introduction to modern
  casting and crystal growth concepts}, Springer Science \& Business Media,
  2010.

\bibitem[Had12]{ha}
M.~Had{\v{z}}i{\'c}, \emph{Orthogonality conditions and asymptotic stability in
  the {S}tefan problem with surface tension}, Arch. Ration. Mech. Anal.
  \textbf{203} (2012), no.~3, 719--745.

\bibitem[HG10]{HaGu}
M.~Had{\v{z}}i{\'c} and Y.~Guo, \emph{Stability in the {S}tefan problem with
  surface tension (i)}, Commun. Partial. Differ. Equ. \textbf{35} (2010),
  no.~2, 201--244.

\bibitem[Kry96]{KrylovHolder}
N.~V. Krylov, \emph{Lectures on elliptic and parabolic equations in
  {H}\"{o}lder spaces}, Graduate Studies in Mathematics, vol.~12, American
  Mathematical Society, Providence, RI, 1996.

\bibitem[KS91]{KaSh}
I.~Karatzas and S.~E. Shreve, \emph{Brownian motion and stochastic calculus},
  second ed., Graduate Texts in Mathematics, vol. 113, Springer-Verlag, New
  York, 1991. \MR{1121940}

\bibitem[Law18]{Law}
G.~F. Lawler, \emph{Notes on the {B}essel process}, Lecture notes. Available at
  http://www.math.uchicago.edu/$\sim$lawler/bessel18new.pdf (2018).

\bibitem[LC31]{LC}
G.~Lam\'{e} and B.~P. Clapeyron, \emph{M\'{e}moire sur la solidification par
  r\'{e}froidissement d'un globe liquide}, Ann. Chimie Physique \textbf{47}
  (1831), 250--256.

\bibitem[LSU68]{Lady}
O.~A. Lady\v{z}enskaja, V.~A. Solonnikov, and N.~N. Ural'ceva, \emph{Linear and
  quasilinear equations of parabolic type}, Translated from the Russian by S.
  Smith. Translations of Mathematical Monographs, Vol. 23, American
  Mathematical Society, Providence, R.I., 1968. \MR{0241822}

\bibitem[Luc90]{Lu}
S.~Luckhaus, \emph{Solutions for the two-phase {S}tefan problem with the
  {G}ibbs-{T}homson law for the melting temperature}, European J. Appl. Math.
  \textbf{1} (1990), no.~2, 101--111. \MR{1117346}

\bibitem[McK05]{McK}
H.~P. McKean, \emph{Stochastic integrals}, AMS Chelsea Publishing, Providence,
  RI, 2005, Reprint of the 1969 edition, with errata. \MR{2169626}

\bibitem[Mei94]{meir}
A.~M. Meirmanov, \emph{The {S}tefan problem with surface tension in the three
  dimensional case with spherical symmetry: Nonexistence of the classical
  solution}, European J. Appl. Math. \textbf{5} (1994), no.~1, 1--19.

\bibitem[NS23]{NaShsurface}
S.~Nadtochiy and M.~Shkolnikov, \emph{Stefan problem with surface tension:
  global existence of physical solutions under radial symmetry}, to appear in
  Probab. Theory Related Fields (2023).

\bibitem[NSZ23]{nsz}
S.~Nadtochiy, M.~Shkolnikov, and X.~Zhang, \emph{Scaling limits of external
  multi-particle {DLA} on the plane and the supercooled {Stefan} problem}, to
  appear in Ann. Inst. Henri Poincar\'{e} Probab. Stat. (2023).

\bibitem[RS06]{RoSa}
R.~Rossi and G.~Savar\'{e}, \emph{Gradient flows of non convex functionals in
  {H}ilbert spaces and applications}, ESAIM Control Optim. Calc. Var.
  \textbf{12} (2006), no.~3, 564--614. \MR{2224826}

\bibitem[Ste89]{Stefan1}
J.~Stefan, \emph{\"{U}ber einige {P}robleme der {T}heorie der
  {W}\"{a}rmeleitung}, Sitzungber., Wien, Akad. Mat. Natur. \textbf{98} (1889),
  473--484.

\bibitem[Ste90a]{Stefan2}
\bysame, \emph{\"{U}ber die {T}heorie der {E}isbildung}, Monatsh. Math. Phys.
  \textbf{1} (1890), no.~1, 1--6. \MR{1546138}

\bibitem[Ste90b]{Stefan3}
\bysame, \emph{{\"U}ber die {V}erdampfung und die {A}ufl\"{o}sung als
  {V}org\"{a}nge der {D}iffusion}, Ann. Physik \textbf{277} (1890), 725--747.

\bibitem[Ste91]{Stefan4}
\bysame, \emph{\"{U}ber die {T}heorie der {E}isbildung, insbesondere \"{u}ber
  die {E}isbildung im {P}olarmeere}, Ann. Physik Chemie \textbf{42} (1891),
  269--286.

\bibitem[Vis87]{vis}
A.~Visintin, \emph{Stefan problem with a kinetic condition at the free
  boundary}, Ann. Mat. Pura Appl. (4) \textbf{146} (1987), 97--122. \MR{916689}

\end{thebibliography}

\bigskip\bigskip

\appendix\label{se:app}

\section{Tail bounds for the Bessel process}

We frequently employ the following tail bounds for the Bessel process $R$. 

\begin{lemma}\label{lm:BesselRunning}
Let $R$ be a three-dimensional Bessel process started from $x>0$ under $\PP^x$, and let $B$ be a standard one-dimensional Brownian motion.~Then,
\begin{eqnarray*}
&& \PP^x\big(\min_{0\le s\le t} R_s\leq x-a\big)
\leq\PP\big(\min_{0\le s\le t} B_s\leq-a\big)=2\Phi(-a/\sqrt{t}),\quad a>0, \\
&& \PP^x\big(\max_{0\le s\le t}\,|R_s-x|\ge a\big) 
\leq 3\PP\big(\max_{0\le s\le t}\,|B_s|\ge a/\sqrt{3}\big)
\leq 12\Phi(-a/\sqrt{3t}),\quad a>0.
\end{eqnarray*}
\end{lemma}

\noindent\textbf{Proof.} Since $\mathrm{d}R_t=R_t^{-1}\,\mathrm{d}t+\mathrm{d}B_t$, $R_0=x$, we have $\min_{0\le s\le t} (R_s-x)\ge\min_{0\le s\le t} B_s$, which yields the first estimate. For the second estimate, we use the fact that $R$ can be realized as the Euclidean norm of a three-dimensional Brownian motion $(x+B^1,\,B^2,\,B^3)$. Hence,
\begin{equation*}
\PP^x\big(\max_{0\le s\le t}\, |R_s-x|\ge a\big) 
\leq\PP\bigg(\max_{0\le s\le t}\,\sum_{i=1}^3\, (B^i_s)^2 \ge a^2\bigg)
\leq\PP\bigg(\bigcup_{i=1}^3\big\{\max_{0\le s\le t}\,|B_s^i|\ge a/\sqrt{3}\big\}\!\bigg).
\end{equation*}
The claimed estimate now follows from a union bound. \qed

\section{Proof of Proposition \ref{KrySaf}}\label{app:KrySaf}

Fix an $\upsilon\in(0,\upsilon_0)$ and denote by $\kappa\ge0$ the $1/2$-Hölder norm of $\Lambda$ on $[t_0+\upsilon,t_0+\upsilon_0]$.~Let $L:=2(\kappa+1)$, pick $y_0>0$ so that
\begin{IEEEeqnarray*}{rCl}
	L^2y_0\leq\sqrt{\upsilon}\wedge\frac{\Lambda_{t_0+\upsilon_0}-r_0}{2},
\end{IEEEeqnarray*}
and consider an arbitrary $s\in[t_0+\upsilon,t_0+\upsilon_0]$. For $0<y\leq\delta$, with $\delta\leq\sqrt{\upsilon}\wedge (\Lambda_{t_0+\upsilon_0}-r_0)/(2L)$, consider the stopping time
\begin{IEEEeqnarray*}{rCl}
	\rho:=\inf\left\{s'>0:\,R_{s'}\le\Lambda_{s-s'}-L\delta\right\}\wedge\delta^2\wedge\tau_{\Lambda_{s-\cdot}}.
\end{IEEEeqnarray*}
By arguing as in the proof of Proposition \ref{prop:MarkovSystem}, we see that 
\begin{equation}\label{eq:Prop2.5split1}
	w(s,\Lambda_s-y)
	=\EE^{\Lambda_s-y}\big[w(s-\rho,R_\rho)\,\bone_{\{\tau_{\Lambda_{s-\cdot}}>\rho\}}\big].
\end{equation}

\smallskip

Since $\mathrm{d}R_t=\frac{1}{R_t}\,\mathrm{d}t+\mathrm{d}B_t$, $R_0=\Lambda_s-y$ under $\PP^{\Lambda_s-y}$, we find that for any $s'\in[0,\rho]\subset[0,\delta^2]\subset[0,\upsilon]$,
\begin{IEEEeqnarray*}{rCl}
	\Lambda_{s-s'}-R_{s'}\le \Lambda_{s-s'}-\Lambda_s+y-B_{s'}
	\le\kappa\delta+\delta-B_{s'},
\end{IEEEeqnarray*}
where we have used the $1/2$-Hölder continuity of $\Lambda$. Recalling that $L=2(\kappa+1)$ we obtain
\begin{IEEEeqnarray*}{rCl}
	\{\tau_{\Lambda_{s-\cdot}}=\rho\}\supset\big\{\max_{0\le s'\le \delta^2} B_{s'}\ge(\kappa+1)\delta,\,\min_{0\le s'\le\delta^2} B_{s'}>-(\kappa+1)\delta\big\}.
\end{IEEEeqnarray*}
By the scaling property of Brownian motion, there exists a $c\in(0,1)$, depending only on $\kappa$, so that 
\begin{IEEEeqnarray*}{rCl}
	L(1-c)>1 \quad\text{and}\quad \PP^{\Lambda_s-y}(\tau_{\Lambda_{s-\cdot}}=\rho)\ge c.
\end{IEEEeqnarray*}
Consequently,
\begin{IEEEeqnarray*}{rCl}
	w(s,\Lambda_s-y)
	=\EE^{\Lambda_s-y}\big[w(s-\rho,R_\rho)\,\bone_{\{\tau_{\Lambda_{s-\cdot}}>\rho\}}\big]
	\le(1-c)\sup_{s-\delta^2\le s'\le s,\,0<y'\le L\delta} w(s',\Lambda_{s'}-y').
\end{IEEEeqnarray*}
Iterating this argument, replacing $\delta$ by $L\delta$ in each iteration, we arrive at
\begin{equation}\label{KrSa1}
	w(s,\Lambda_s-y)
	\le(1-c)^{k+1}\sup_{s-(1+L^2+\cdots+L^{2k})\delta^2\le s'\le s,\,0<y'\le L^{k+1}\delta} w(s',\Lambda_{s'}-y'),
\end{equation}
as long as $(1+L^2+\cdots+L^{2k})\delta^2\leq\upsilon$ and 
$L^{k+1}\delta\le(\Lambda_{t_0+\upsilon_0}-r_0)/2$.

\medskip

Notice that $L=2(\kappa+1)\ge2$. Thus, as long as $L^{k+1}\delta\le L^2 y_0$, we have 
\begin{IEEEeqnarray*}{rCl}
	(1+L^2+\cdots+L^{2k})\delta^2=\frac{L^{2k+2}-1}{L^2-1}\delta^2
	\leq L^4y_0^2\leq\upsilon\quad\text{and}\quad
	L^{k+1}\delta\leq L^2 y_0\leq\frac{\Lambda_{t_0+\upsilon_0}-r_0}{2}
\end{IEEEeqnarray*}
by the definition of $y_0$. Finally, for fixed $y\in(0,y_0]$, we choose $\delta=y$ and 
\begin{IEEEeqnarray*}{rCl}
	k=\left\lfloor\frac{\log(y_0/y)}{\log(L)}\right\rfloor+1\ge1.
\end{IEEEeqnarray*}
Then, the constraint $L^{k+1}\delta \le L^2y_0$ is satisfied, and for $\varrho:=-\log(1-c)/\log(L)\in(0,1)$,
\begin{equation}\label{KrSa2}
	(1-c)^{k+1}\leq(1-c)^{\frac{\log(y_0/y)}{\log(L)}}=(y/y_0)^\varrho.
\end{equation}
Combining \eqref{KrSa1}, \eqref{KrSa2} and 
\begin{IEEEeqnarray*}{rCl}
	\sup_{s-(1+L^2+\cdots+L^{2k})\delta^2\le s'\le s,\,0<y'\le L^{k+1}\delta} w(s',\Lambda_{s'}-y')
	\le\|w(0-,\cdot)\|_\infty
\end{IEEEeqnarray*}
(see \eqref{eq.sec3.sergey.w.def}) we obtain the desired 
\begin{IEEEeqnarray*}{rCl}
	\qquad\qquad\quad\;\, w(s,\Lambda_s-y)\le\|w(0-,\cdot)\|_\infty\,y_0^{-\varrho}\,y^\varrho,\quad y\in(0,y_0],\quad s\in[t_0+\upsilon,t_0+\upsilon_0]. \qquad\qquad\quad\;\,\qed
\end{IEEEeqnarray*}

\section{Estimates on some Brownian densities} \label{BMestimate}

We group some frequently used estimates on Brownian densities into the following lemma.
\begin{lemma}\label{lem:BMestimates}
\begin{enumerate}[(a)]
\item For all $\mu\in\RR$,
\begin{equation} \label{eq:GaussianMWestimate1}
\int_{-\infty}^\infty e^{\mu b-\frac{\mu^2}{2}}\,\frac{1}{\sqrt{2\pi}}\,e^{-\frac{b^2}{2}}\,\mathrm{d}b=1.  
\end{equation}
\item For all $\mu\ge0$ and $s>0$,
\begin{equation} \label{eq:GaussianMWestimate2}
\int_0^\infty e^{\mu\sqrt{s}b-\frac{\mu^2s}{2}}\frac{b}{\sqrt{2\pi s}}\,e^{-\frac{b^2}{2}}\,\mathrm{d}b\le\frac{1}{\sqrt{2\pi s}}+\mu.
\end{equation}
\item For all $\mu\ge0$, $y>0$ and $s>0$,
\begin{equation} \label{eq:GaussianMWestimate3}
\frac1y\int_0^y\int_{-\infty}^m e^{-\mu b-\frac{\mu^2s}{2}}\,
\frac{2(2m-b)}{s\sqrt{2\pi s}}\,e^{-\frac{(2m-b)^2}{2s}}\,\mathrm{d}b\,\mathrm{d}m
\le 2\left(\frac{1}{\sqrt{2\pi s}}+\mu\right). 
\end{equation}
\item For all $\varrho\in(0,1]$, $\mu\ge0$, $y>0$ and $s>0$,
\begin{equation}\label{eq:GaussianMWestimate4}
\frac1y\int_0^y \int_{-\infty}^m (y-b)^\varrho\,
e^{-\mu b-\frac{\mu^2s}{2}}\,
\frac{2(2m-b)}{s\sqrt{2\pi s}}\,
e^{-\frac{(2m-b)^2}{2s}}\,\mathrm{d}b\,\mathrm{d}m 
\le 2y^\varrho \left(\frac{1}{\sqrt{2\pi s}}+\mu\right)+2s^{\frac{\varrho-1}{2}}(2+\mu^2s).   
\end{equation}
\end{enumerate}
\end{lemma}

\noindent\textbf{Proof. (a)} This is simply a restatement of the formula for the moment generating function of the standard normal distribution. 

\medskip
	
\noindent\textbf{(b)} We integrate by parts to bound the integral according to
\begin{IEEEeqnarray*}{rCl}
\int_0^\infty e^{\mu\sqrt{s}b-\frac{\mu^2s}{2}}\,
\frac{b}{\sqrt{2\pi s}}\,e^{-\frac{b^2}{2}}\,\mathrm{d}b
&=&\frac{1}{\sqrt{2\pi s}}\,e^{-\frac{\mu^2s}{2}}+\mu\int_0^\infty 
e^{\mu \sqrt{s}b-\frac{\mu^2s}{2}}\,\frac{1}{\sqrt{2\pi}}\,e^{-\frac{b^2}{2}}\,\mathrm{d}b \\
&\leq&\frac{1}{\sqrt{2\pi s}}+\mu\int_{-\infty}^\infty 
e^{\mu\sqrt{s}b-\frac{\mu^2 s}{2}}\,\frac{1}{\sqrt{2\pi}}\,e^{-\frac{b^2}{2}}\,\mathrm{d}b
=\frac{1}{\sqrt{2\pi s}}+\mu,
\end{IEEEeqnarray*}
where the latter equality is due to the identity \eqref{eq:GaussianMWestimate1}.
	
\medskip
	
\noindent\textbf{(c)} Changing variables from $b$ to $2m-b$ and then from $b$ to $\sqrt{s}b$ in the inner integral we find 
\begin{IEEEeqnarray*}{rCl}
\int_{-\infty}^m e^{-\mu b-\frac{\mu^2 s}{2}}
\,\frac{2(2m-b)}{s\sqrt{2\pi s}}\,e^{-\frac{(2m-b)^2}{2s}}\mathrm{d}b
&=&\int_m^\infty e^{\mu(b-2m)-\frac{\mu^2s}{2}}\,\frac{2b}{s\sqrt{2\pi s}}\,e^{-\frac{b^2}{2s}}\,\mathrm{d}b \\
&\leq&\int_0^\infty e^{\mu b-\frac{\mu^2s}{2}}\,\frac{2b}{s\sqrt{2\pi s}}\,e^{-\frac{b^2}{2s}}\,\mathrm{d}b \\
&=&\int_0^\infty e^{\mu\sqrt{s}b-\frac{\mu^2s}{2}}\,\frac{2b}{\sqrt{2\pi s}}\, e^{-\frac{b^2}{2}}\,\mathrm{d}b
\le 2\left(\frac{1}{\sqrt{2\pi s}}+\mu\right)\!,
\end{IEEEeqnarray*}
where the last inequality is a consequence of \eqref{eq:GaussianMWestimate2}. The estimate \eqref{eq:GaussianMWestimate3} readily follows. 
	
\medskip
	
\noindent\textbf{(d)} We slightly simplify the inner integral by the change of variables from $2m-b$ to $b$: 
\begin{equation} \label{eq:GaussianMWestimate4a}
\begin{split}
\int_{-\infty}^m (y-b)^\varrho\,
e^{-\mu b-\frac{\mu^2s}{2}}\,
\frac{2(2m-b)}{s\sqrt{2\pi s}}\,
e^{-\frac{(2m-b)^2}{2s}}\,\mathrm{d}b
&=\int_m^\infty (y+b-2m)^\varrho\,e^{\mu(b-2m)-\frac{\mu^2s}{2}}\,
\frac{2b}{s\sqrt{2\pi s}}\,e^{-\frac{b^2}{2s}}\,\mathrm{d}b \\
&\le\int_0^\infty (y+b)^\varrho\,e^{\mu b-\frac{\mu^2s}{2}}\,\frac{2b}{s\sqrt{2\pi s}}\,e^{-\frac{b^2}{2s}}\,\mathrm{d}b.
\end{split}
\end{equation}
Note that $m$ does not appear in the latter expression. Thanks to the elementary  $(y+b)^\varrho\leq y^\varrho+b^\varrho$, $y,b>0$ (recall that $\varrho\in(0,1]$), we obtain the further bounds
\begin{equation} \label{eq:GaussianMWestimate4b}
\begin{split}
&\,\int_0^\infty (y+b)^\varrho\,e^{\mu b-\frac{\mu^2s}{2}}\,
\frac{2b}{s\sqrt{2\pi s}}\,e^{-\frac{b^2}{2s}}\,\mathrm{d}b \\
&\le y^\varrho\int_0^\infty e^{\mu b-\frac{\mu^2s}{2}}\,\frac{2b}{s\sqrt{2\pi s}}\,e^{-\frac{b^2}{2s}}\,\mathrm{d}b+\int_0^\infty b^\varrho\,e^{\mu b-\frac{\mu^2s}{2}}\,\frac{2b}{s\sqrt{2\pi s}}\,e^{-\frac{b^2}{2s}}\,\mathrm{d}b \\
&= y^\varrho\int_0^\infty e^{\mu\sqrt{s}b-\frac{\mu^2s}{2}}\,\frac{2b}{\sqrt{2\pi s}}\,
e^{-\frac{b^2}{2}}\,\mathrm{d}b
+s^{\frac{\varrho-1}{2}}\int_0^\infty b^{\varrho+1}\,e^{\mu \sqrt{s}b-\frac{\mu^2s}{2}}\,\frac{2}{\sqrt{2\pi}}\,e^{-\frac{b^2}{2}}\mathrm{d}b.
\end{split}
\end{equation}
Finally, observing the elementary $b^{\varrho+1}\leq1+b^2$, $b>0$ and integrating by parts twice we get 
\begin{equation} \label{eq:GaussianMWestimate4c}
\begin{split}	
\int_0^\infty b^{\varrho+1}\,e^{\mu \sqrt{s}b-\frac{\mu^2s}{2}}\,
\frac{1}{\sqrt{2\pi}}\,e^{-\frac{b^2}{2}}\mathrm{d}b
&\le\int_{-\infty}^\infty(1+b^2)\,e^{\mu \sqrt{s}b-\frac{\mu^2s}{2}}\,\frac{1}{\sqrt{2\pi}}\,e^{-\frac{b^2}{2}}\,\mathrm{d}b \\
&=1+\int_{-\infty}^\infty (1+\mu\sqrt{s}b)\,
e^{\mu\sqrt{s}b-\frac{\mu^2s}{2}}\,\frac{1}{\sqrt{2\pi}}\,e^{-\frac{b^2}{2}}\,\mathrm{d}b \\
&=2+\mu^2s\int_{-\infty}^\infty e^{\mu \sqrt{s}b-\frac{\mu^2s}{2}} \,\frac{1}{\sqrt{2\pi}}\,e^{-\frac{b^2}{2}}\,\mathrm{d}b=2+\mu^2s.
\end{split}
\end{equation}
The desired result now follows by combining \eqref{eq:GaussianMWestimate2}, \eqref{eq:GaussianMWestimate4a}, \eqref{eq:GaussianMWestimate4b} and \eqref{eq:GaussianMWestimate4c}. \qed

\section{Proof of Proposition \ref{prop:InteriorSmoothness}} \label{app:Cinf}

For any $\varphi\in C_c^\infty(\RR^2)$ supported in $\stackrel{\circ}{D}$, we may apply It\^o's formula to get
\begin{IEEEeqnarray*}{rCl}
	\varphi(\tau^-_{\Lambda_{0+\cdot}}\!\wedge\!\tau_{r_0}\!\wedge\!\zeta,R_{\tau^-_{\Lambda_{0+\cdot}}\wedge\tau_{r_0}\wedge\zeta})
	&=&\varphi(0,R_0)+\int_0^{\tau^-_{\Lambda_{0+\cdot}}\!\wedge\!\tau_{r_0}\!\wedge\!\zeta} \partial_t\varphi(t,R_t)\!+\!\frac{1}{R_t}\partial_x\varphi(t,R_t)\!+\!
	\frac12\partial_{xx}\varphi(t,R_t)\,\mathrm{d}t \\
	&&+\int_0^{\tau^-_{\Lambda_{0+\cdot}}\!\wedge\!\tau_{r_0}\!\wedge\!\zeta}\partial_x\varphi(t,R_t)\,\mathrm{d}B_t.    
\end{IEEEeqnarray*}
Since $\varphi$ vanishes outside of $\stackrel{\circ}{D}$,  $\varphi(\tau^-_{\Lambda_{0+\cdot}}\!\wedge\!\tau_{r_0}\!\wedge\!\zeta,R_{\tau^-_{\Lambda_{0+\cdot}}\wedge\tau_{r_0}\wedge\zeta})=\varphi(0,R_0)=0$ almost surely. Also,
\begin{IEEEeqnarray*}{rCl}
	\EE\bigg[\int_0^{\tau^-_{\Lambda_{0+\cdot}}\!\wedge\!\tau_{r_0}\!\wedge\!\zeta}\partial_x\varphi(t,R_t)\,\mathrm{d}B_t\bigg]=0,
\end{IEEEeqnarray*}
as $\partial_x\varphi$ is bounded. Consequently,
\begin{IEEEeqnarray*}{rCl}
	0&=&\EE\bigg[\int_0^{\tau^-_{\Lambda_{0+\cdot}}\!\wedge\!\tau_{r_0}\!\wedge\!\zeta}\partial_t\varphi(t,R_t)+\frac{1}{R_t}\partial_x\varphi(t,R_t)+\frac12\partial_{xx}\varphi(t,R_t)\,\mathrm{d}t\bigg] \\
	&=&\int_0^\zeta\EE\Big[\!\Big(\partial_t\varphi(t,R_t)+\frac{1}{R_t}\partial_x\varphi(t,R_t)+\frac12\partial_{xx}\varphi(t,R_t)\!\Big)\,\bone_{\{\tau^-_{\Lambda_{0+\cdot}}\!\wedge\!\tau_{r_0}>t\}}\Big]\,\mathrm{d}t.
\end{IEEEeqnarray*}
Recalling from Section \ref{se:1.5} that $x^2w(t,x)$ is the density of $R_t$ on $\{\tau^-_{\Lambda_{0+\cdot}}\!\wedge\!\tau_{r_0}>t\}$ we arrive at
\begin{IEEEeqnarray*}{rCl}
	0=\int_0^\zeta \int_{r_0}^\infty \Big(\partial_t\varphi(t,x)+\frac{1}{x}\partial_x\varphi(t,x)+\frac12\partial_{xx}\varphi(t,x)\!\Big)\,x^2w(t,x)\,\mathrm{d}x\,\mathrm{d}t.
\end{IEEEeqnarray*}
Since $\frac1x\partial_x+\frac12\partial_{xx}$ is self-adjoint on  $L^2((r_0,\infty),x^2\,\mathrm{d}x)$, we infer that $w$ is a weak solution of
\begin{IEEEeqnarray*}{rCl}
	\partial_tw(t,x)=\frac12\partial_{xx}w(t,x)+\frac1x\partial_xw(t,x)    
\end{IEEEeqnarray*}
on $\stackrel{\circ}{D}$. The proposition now follows from Weyl's lemma (see, e.g., \cite[p.85]{McK}). \qed

\section{Changes of Monotonicity} \label{se:4}

In this appendix, we focus on (the interior of) the solid phase 
\begin{IEEEeqnarray*}{rCl}
	\stackrel{\circ}{D}_{sol}\,:=\{(t,x)\in(0,\zeta)\times(r_0,\infty):\,x<\Lambda_t\},
\end{IEEEeqnarray*}
and notice in passing that, on $\stackrel{\circ}{D}_{sol}$, the function $\partial_xw$ is a classical solution of the PDE
\begin{equation} \label{eq:PDEpartialxwd=3}
	\partial_t\partial_xw(t,x)=\frac12\partial_{xx}\partial_xw(t,x)+\frac1x\partial_x\partial_xw(t,x)-\frac{1}{x^2}\partial_xw(t,x). 
\end{equation}

The next proposition is the main result of this appendix.

\begin{proposition}
	\label{le:wMonotonicityChanging}
	For any $t\in[0,\zeta)$, the function $w(t-,\Lambda_{t-}-\cdot)$ changes monotonicity finitely often on compact sub-intervals of $[0,\Lambda_{t-}-r_0)$.
\end{proposition}


\noindent\textbf{Proof.} We follow \cite[proof of Lemma 4.1]{dns} and only explain the differences in detail. 

\medskip

\noindent\textit{First Step.} This step is the same. The task is to obtain, for each fixed $R\in(0,\Lambda_{t_*-}-r_0)$, a uniform upper bound $M<\infty$ on the number of monotonicity changes of $w(s,\Lambda_s-\cdot)$ on $[0,R]$ for all $s\in[t_*/2,t_*)$, where we suppose 
\begin{IEEEeqnarray*}{rCl}
t_*:=\inf\cT:=\inf\big\{t\in(0,\zeta):\,w(t-,\Lambda_{t-}-\cdot)\text{ violates the monotonicity-changing property}\big\}<\infty.
\end{IEEEeqnarray*}
Propositions \ref{prop:partialxwnegativeRegular} and \ref{prop:InteriorSmoothness} imply that $t_*\in\cT$. In particular, $t_*>0$. If the described upper bound $M<\infty$ does exist, we arrive at a contradiction. Indeed, we can then find a sequence $t_n\in[t_*/2,t_*)$, $t_n\uparrow t_*$ and $0\le x_1^{(n)}\le x_2^{(n)}\le\cdots\le x_{M+1}^{(n)}=R$, converging as $n\to\infty$ to $0\le x_1\leq x_2\le\cdots\le x_{M+1}=R$, respectively, such that each $w(t_n,\Lambda_{t_n}-\cdot)$ is non-decreasing on $[0,x_1^{(n)}],\,[x_2^{(n)},x_3^{(n)}],\,\ldots$ and non-increasing on $[x_1^{(n)},x_2^{(n)}],\, [x_3^{(n)},x_4^{(n)}],\,\ldots\,$. But then $w(t_*-,\Lambda_{t_*-}-\cdot)$ is non-decreasing on $[0,x_1],\,[x_2,x_3],\,\ldots$ and non-increasing on $[x_1,x_2],\,[x_3,x_4],\,\ldots\,$, a contradiction to $t_*\in\cT$.

\medskip

\noindent\textit{Second Step.} This step is the same, apart from replacing \cite[Lemma 4.2]{dns} by Lemma \ref{lm:ZeroCurvesLocal} below. The sawtooth-shaped curve $\theta$ now lives in $[t_*/2,t_*]\times[\Lambda_{t_*-}-R-\delta,\Lambda_{t_*-}-R]$. We need to study the local behavior of the zero curves of $\partial_xw$ to bound the number of their intersections with~$\theta$. This is done in Lemma \ref{lm:ZeroCurvesLocal}, where the statement on the curve $\xi$ is weaker than in \cite[Lemma 4.2]{dns}, but strong enough to ensure that $\xi$ intersects every straight segment of $\theta$ at most once. 

\medskip

\noindent\textit{Third Step.} This step is the same up to an obvious change of direction: We now track the highest zero curves. More specifically, the assertion therein now goes as follows.

\smallskip

\noindent\textbf{Assertion.} For all $x\in(\theta_s,\Lambda_s)$ with $\partial_xw(s,x)=0$, there exist a $\underline{t}\in[t_*/2,s)$ and a continuous function $\xi\!:[\underline{t},s]\to(r_0,\infty)$ such that
\begin{enumerate}[(a)]
\item $\xi_s=x$ and $\xi_t\in(\theta_t,\Lambda_t)$, $\partial_xw(t,\xi_t)=0$ when $t\in(\underline{t},s]$;
\item for all $t\in(\underline{t},s)$, there exists a neighborhood
\begin{IEEEeqnarray*}{rCl}
(t-\delta_1,t+\delta_1)\times(\xi_t-\delta_2,\xi_t+\delta_2)\subset\{(t,x)\in(t_*/2,s)\times(r_0,\infty):\,\theta_t<x<\Lambda_t\}=:\Gamma_{t_*/2,s}    
\end{IEEEeqnarray*}
on which $\partial_xw(s',y)=0$ implies $y\le\xi_{s'}$;
\item one has
\begin{IEEEeqnarray*}{rCl}
(\underline{t},\xi_{\underline{t}})\in\theta_{[t_*/2,s]}\cup(\{t_*/2\}\times[\Lambda_{t_*-}-R,\Lambda_{t_*/2}])\cup\Lambda_{[t_*/2,s]}=:\partial_{\text{par}}\Gamma_{t_*/2,s}.   
\end{IEEEeqnarray*}
\end{enumerate}

\smallskip

\noindent\textit{Fourth Step.} We can use the same argument. Recall that $\partial_xw$ solves the PDE \eqref{eq:PDEpartialxwd=3} and note that $1/x$, $-1/x^2$ remain bounded for $x\ge r_0>0$. Thus, Feynman-Kac formulas (the standard one and the one shown in the proof of Lemma \ref{lemE3}) rule out the intersection of zero curves. 

\medskip

\noindent\textit{Fifth Step.} This step is the same up to an obvious change of direction. \qed

\begin{lemma}
\label{lm:ZeroCurvesLocal}
Let $(t,x)\in(0,\zeta)\times(r_0,\infty)$ satisfy $\partial_xw(t,x)=0$ and let $k\ge1$ be the smallest integer such that $\partial_x^k\partial_xw(t,x)\neq0$. Then, there is a neighborhood $(t-\delta_1,t+\delta_1)\times(x-\delta_2,x+\delta_2)$ on which the zero set of $\partial_xw$ is the union of $k$ curves. Of those, $2\lfloor k/2\rfloor$ are graphs of continuous functions on $(t-\delta_1,t]$, formed by solutions of $(x'-x)^2\sigma(x'-x)=t'-t$ for smooth $\sigma$ with $\sigma(0)<0$. If $k$ is odd, there is another curve given by the graph of a continuous function $\xi$ on $(t-\delta_1,t+\delta_1)$  which satisfies $\xi(t+s)=O(|s|)$, $s\to0$.
\end{lemma}

\noindent\textbf{Proof.} By Proposition \ref{prop:InteriorAnalyticity}, $\partial_xw(t+s,x+y)$ can be extended analytically in $y$ to a neighborhood $V_1$ of $0\in\CC$, and $V_1$ can be chosen independently of $s\in U_1$, for some neighborhood $U_1$ of $0\in\RR$. Then, for $(s,y)$ in a neighborhood $W_1\subset U_1\times V_1$ of $(0,0)\in\RR\times\CC$, it holds $\partial_x^{k+1}w(t+s,x+y)\neq0$. Upon shrinking $U_1$ and $V_1$ if necessary, we achieve $W_1=U_1\times V_1$ and $\inf_{W_1} |\partial_x^{k+1}w(t+\cdot,x+\cdot)|>0$.

\medskip

As discussed in the last paragraph of \cite[proof of Lemma 4.2]{dns}, the following intermediate result between the Weierstrass Preparation Theorem \cite[Theorem 6.2]{ChowHale} and the Malgrange Preparation Theorem \cite[Theorem 7.1]{ChowHale} holds: There exists a neighborhood $W_2$ of $(0,0)\in\RR\times\CC$ such that $\partial_xw(t+s,x+y)=q(s,y)\,\Gamma(s,y)$, $(s,y)\in W_2$, where (a) $q$ is a smooth nowhere-vanishing function; (b) $\Gamma(s,y)=y^k-\sum_{l=0}^{k-1}a_l(s)y^l$ is a degree-$k$ polynomial in $y$ with coefficients $a_l(s)$ being real-valued smooth functions of $s$. By differentiating with respect to $y$, one sees that $\Gamma(0,y)=y^k$. Since in $W_2$ the zero set of $\partial_xw(t+\cdot,x+\cdot)$ equals that of $\Gamma$, it follows from the continuous dependence of the roots upon the coefficients that all the roots of $\Gamma(s,\cdot)$ approach $0$ as $|s|\downarrow0$. More specifically, we can assume that $W_2=U_2\times V_2\subset U_1\times V_1$. Then we can find a neighborhood $U_3$ of $0\in\RR$ so that for $s\in U_3$, the set $V_2$ contains exactly $k$ zeros (counting multiplicity) of $y\mapsto\partial_xw(t+s,x+y)$. These zeros, with $s$ ranging over $U_3$, exhaust the zero set of $\partial_xw(t+\cdot,x+\cdot)$ in $U_3\times V_2$. In the sequel, we restrict ourselves to $U_3\times V_2$ and locate these zeros.

\medskip

The function $\partial_xw$ solves the PDE \eqref{eq:PDEpartialxwd=3}, which is of the form studied in \cite[Section 5]{AF}. The PDE \eqref{eq:PDEpartialxwd=3} and the definition of $k$ yield for $2m+n<k$, 
\begin{IEEEeqnarray*}{rCl}
\partial_t^m\partial_x^n\partial_xw(t,x)=\partial_x^n\partial_t^m\partial_xw(t,x)=0.    
\end{IEEEeqnarray*}
Thus, for $\sigma<0$, we obtain the expansion
\begin{IEEEeqnarray*}{rCl}
\partial_xw(t+y^2\sigma,x+y)
&=& y^k\sum_{2m+n=k}\frac{1}{m! n!} \partial_t^m\partial_x^n\partial_xw(t,x)\sigma^m+y^{k+1}H(\sigma,y) \\
&=& y^k\partial_x^{k+1}w(t,x)\sum_{2m+n=k}\frac{1}{m!n!2^m}\sigma^m+y^{k+1}H(\sigma,y),
\end{IEEEeqnarray*}
where the smooth $H$ varies from line to line. Repeating \cite[proof of Lemma 4.2, First Step]{dns} one finds $\lfloor k/2\rfloor$ ``parabolas'' of the form $s=y^2\sigma(y)$ for $s<0$, where $\sigma$ is some smooth and strictly negative function, giving $2\lfloor k/2\rfloor$ zero curves that are smooth in $s<0$ and behave like $\sqrt{-s}$ as $s\uparrow0$. The $2\lfloor k/2\rfloor$ roots are all simple, as follows from the argument in \cite[proof of Theorem 5.1]{AF}. We record a rectangle $(-\delta_{11},0]\times(-\delta_{21},\delta_{21})\subset U_3\times V_2$
containing the $\lfloor k/2\rfloor$ parabolas. 

\medskip


For $s>0$, the situation differs from the one in \cite[proof of Lemma 4.2, First Step]{dns} because here the complex roots of $\partial_xw(t+s,x+\cdot)$ may not be purely imaginary, making it harder to locate them. Rather than to capture all the complex-valued roots, we focus on the behavior of the real roots. For $s>0$ and even $k=2l$, we use the expansion
\begin{IEEEeqnarray*}{rCl}
\partial_xw(t+s,x+y)&=&\sum_{2m+n=2l}\frac{1}{m!n!}\partial_t^m\partial_x^n\partial_xw(t,x)s^my^n+o(s^l+y^{2l}) \\
&=&\sum_{m=0}^l\frac{1}{m!(2l-2m)!2^m}\partial_x^{2l+1}w(t,x)s^my^{2l-2m}+o(s^l+y^{2l}),
\end{IEEEeqnarray*}
which has the same sign as $\partial_x^{2l+1}w(t,x)\neq0$ for all $(s,y)\in(0,\delta_{12})\times(-\delta_{22},\delta_{22})$, where $\delta_{12},\delta_{22}>0$ are small enough that $[0,\delta_{12})\times(-\delta_{22},\delta_{22})\subset U_3\times V_2$. Thus, $\partial_xw(t+s,x+\cdot)$ has no real roots in $(-\delta_{22},\delta_{22})$ for $s\in(0,\delta_{12})$. With $\delta_1:=\delta_{11}\!\wedge\!\delta_{12}$, $\delta_2:=\delta_{21}\!\wedge\!\delta_{22}$ we get the desired picture for even $k$.

\medskip

For $s>0$ and odd $k=2l+1$, we expand $\partial_{xx}w(t+s,x+y)$ and find that it has the same sign as $\partial_x^{2l+2}w(t,x)\neq0$ for all $(s,y)\in(0,\delta_{12})\times(-\delta_{22},\delta_{22})\subset U_3\times V_2$. As a result, $\partial_xw(t+s,x+\cdot)$ is strictly monotone on $(-\delta_{22},\delta_{22})$, and hence has at most one real root thereon, for all $s\in(0,\delta_{12})$. Since $\Gamma(s,\cdot)$ is a degree-$(2l+1)$ polynomial with real coefficients, $\partial_xw(t+s,x+\cdot)$ has exactly one real root in $(-\delta_{22},\delta_{22})$ for $s\in(0,\delta_{12})$. Recall that for $s<0$, we are also missing one root if $k$ is odd. This root $\xi(s)$ is real and simple also for $s<0$. Moreover, $\xi(\cdot)$ is continuous on a neighborhood of zero by the continuous dependence of the roots on the coefficients, and by the implicit function theorem, $\xi(\cdot)$ is $C^\infty$ on $(-\delta_{13},0)\cup(0,\delta_{13})$ for some $\delta_{13}\in(0,\delta_{11}\wedge\delta_{12})$. It remains to analyze $\xi(\cdot)$ near zero. We distinguish between $r<\infty$ and $r=\infty$, where $r$ is the first positive integer such that
\begin{IEEEeqnarray*}{rCl}
	\partial_t^{l+r}\partial_xw(t,x)=\partial_x\partial_t^{l+r}w(t,x)=\partial_x\Big(\frac12\partial_{xx}+\frac1x\partial_x\Big)^{l+r}w(t,x)\neq0.    
\end{IEEEeqnarray*}

Let $r<\infty$. Then, we expand $\partial_xw(t+s,x+\eta s^r)$ for $s\to0$:
\begin{IEEEeqnarray*}{rCl}
	\partial_xw(t+s,x+\eta s^r)&=&\sum_{2m+n\ge k,\,m+rn\le l+r}
	\frac{1}{m!n!}\partial_t^m\partial_x^n\partial_xw(t,x) s^m(\eta s^r)^n+s^{l+r+1}L(s,\eta)\\
	&=&\frac{1}{(l+r)!}\partial_x\Big(\frac12\partial_{xx}\!+\!\frac1x\partial_x\Big)^{l+r}w(t,x)s^{l+r}+\frac{1}{l!2^l}\partial_x^{k+1}w(t,x)s^{l+r}\eta+s^{l+r+1}L(s,\eta),
\end{IEEEeqnarray*}
where $L$ is smooth, and in the second equality we have exploited that the only non-negative integer solutions of the system of inequalities $2m+n\ge k$, $m+rn\le l+r$ are
\begin{IEEEeqnarray*}{rCl}
\{(m,n):\;\; n=0,\;m=l+1,\,l+2,\,\ldots,\,l+r\quad\text{or}\quad n=1,\,m=l\},
\end{IEEEeqnarray*}
as well as the definition of $r$. By the implicit function theorem, we find for $s$ in a neighborhood of zero a root of the form $y=\eta(s)s^r$, where $\eta$ is a smooth function of $s$ satisfying
\begin{IEEEeqnarray*}{rCl}
	\eta(0)=-\frac{l!2^l}{(l+r)!}\frac{\partial_x\left(\frac12\partial_{xx}+\frac1x\partial_x\right)^{l+r}w(t,x)}{\partial_x^{k+1}w(t,x)}\neq0.
\end{IEEEeqnarray*}
This root behaves polynomially in $s$, and thus, $\eta(s)s^r=\xi(s)$, which proves $\xi(s)=O(|s|)$, $s\to0$. 

\medskip

Let now $r=\infty$. From the definition of $k=2l+1$ and the PDE \eqref{eq:PDEpartialxwd=3} for $\partial_xw$ we infer that
\begin{IEEEeqnarray*}{rCl}
\partial_xw(t,x)=\partial_t\partial_xw(t,x)=\cdots=\partial_t^l\partial_xw(t,x)=0.    
\end{IEEEeqnarray*}
On the other hand, the definition of $r$ together with $r=\infty$ imply
\begin{IEEEeqnarray*}{rCl}
0=\partial_t^{l+1}\partial_xw(t,x)=\partial_t^{l+2}\partial_xw(t,x)=\cdots,
\end{IEEEeqnarray*}
so that $\partial_t^m\partial_xw(t,x)=0$, $m\in\NN_0$. For $s<0$, we can argue as in \cite[proof of Lemma 4.2, Fourth Step]{dns}: Since $\partial_xw(t+s,x+\cdot)$ has $2\lfloor k/2\rfloor$ roots behaving like $\sqrt{-s}$ as $s\uparrow0$, the decomposition $\partial_xw(t+s,x)=q(s,0)c(s)s^{\lfloor k/2\rfloor}\xi(s)$ is valid, with a nowhere-vanishing smooth $c(\cdot)$. Differentiating this identity repeatedly with respect to $s$ we obtain $\xi^{(m)}(0-)=0$ for all $m\in\NN_0$.

\medskip

For $s>0$, we take an arbitrary $\rho\in\NN$ and expand $\partial_xw(t+s,x+\eta s^\rho)$ for $s\downarrow 0$: 
\begin{IEEEeqnarray*}{rCl}
\partial_xw(t+s,x+\eta s^\rho)&=&\sum_{2m+n\ge k,\,m+\rho n\le l+\rho} \frac{1}{m!n!}\partial_t^m\partial_x^n\partial_xw(t,x)s^m(\eta s^\rho)^n+s^{l+\rho+1}L(s,\eta) \\
&=&\frac{1}{l!2^l}\partial_x^{k+1}w(t,x)s^{l+\rho}\eta+s^{l+\rho+1}L(s,\eta),
\end{IEEEeqnarray*}
where $L$ is a smooth function on a neighborhood of $0\in[0,\infty)\times\RR$. Thus, for  any $\eta\neq0$ fixed, $\partial_xw(t+s,x+\eta s^\rho)$ has the same sign as $\partial_x^{k+1}w(t,x)\eta$ for all $s>0$ small enough.~Therefore, for $\eta>0$ fixed, $\partial_xw(t+s,x+\eta s^\rho)$ and $\partial_xw(t+s,x-\eta s^\rho)$ are of different sign for all $s>0$ small enough, meaning that $\partial_xw(t+s,x+\cdot)$ has a zero in $(-\eta s^\rho,\eta s^\rho)$. We deduce $-\eta s^\rho<\xi(s)<\eta_1 s^\rho$, for all $s>0$ small enough. In particular, $-\eta\le\liminf_{s\downarrow0} s^{-\rho}\xi(s)\le
\limsup_{s\downarrow0} s^{-\rho}\xi(s)\le\eta$, and since $\eta>0$ and $\rho\in\NN$ were arbitrary, $\xi(s)=o(s^\rho)$ as $s\downarrow0$ for all $\rho\in\NN$. This finishes the proof. \qed

\begin{lemma}
\label{lemE3}
Given $0<s_*<s<\zeta$, let $\xi^1,\xi^2:[s_*,s]\to(r_0,\infty)$ be two continuous paths satisfying
\begin{enumerate}[(a)]
\item $\xi_{s_*}^1=\xi_{s_*}^2=\Lambda_{s_*}$;
\item $\Lambda_t>\xi_t^1>\xi_t^2$, $t\in(s_*,s]$; and 
\item $\partial_xw(t,\xi_t^1)=\partial_xw(t,\xi_t^2)=0$, $t\in(s_*,s]$. 
\end{enumerate}
Then, $\partial_xw(t,x)=0$, $x\in(\xi_t^2,\xi_t^1)$, $t\in(s_*,s]$.
\end{lemma}

\noindent\textbf{Proof.}
Fix $t\in(s_*,s]$ and $x\in(\xi_t^2,\xi_t^1)$. Let $s_n:=s_*+1/n$, where $n\in\NN$ is sufficiently large to guarantee $s_n<t$. Consider the stopping times
\begin{IEEEeqnarray*}{rCl}
\tau_n=\inf\{t'>0:\,R_{t'}\notin(\xi_{t-t'}^2,\xi_{t-t'}^1)\}\wedge(t-s_n).    
\end{IEEEeqnarray*}
Since $\{(s',x)\!:s'\in[s_n,t],\,x\in[\xi_{s'}^2,\xi_{s'}^1]\}\subset\,\stackrel{\circ}{D}_{sol}$, we have the Feynman-Kac formula 
\begin{IEEEeqnarray*}{rCl}
\partial_xw(t,x)=\EE^x\Big[e^{-\int_0^{\tau_n}\frac{1}{R_{t'}^2}\,\mathrm{d}t'}\,\partial_xw(t-\tau_n,R_{\tau_n})\Big]
=\EE^x\Big[e^{-\int_0^{t-s_n}\frac{1}{R_{t'}^2}\,\mathrm{d}t'}\,\partial_xw(s_n,R_{t-s_n})\,\bone_{\{\tau_n=t-s_n\}}\Big].
\end{IEEEeqnarray*}
(Recall the assumption (c).) Girsanov's theorem and Hölder's inequality then yield
\begin{equation} \label{eq:partialxwestimate1}
\begin{split}
|\partial_xw(t,x)|&\le\EE^x\left[|\partial_xw(s_n,R_{t-s_n})|\,\bone_{\{\tau_n=t-s_n\}}\right] \\
&=\EE^x\Big[e^{\int_0^{t-s_n} \frac{1}{B_{t'}}\,\mathrm{d}B_{t'}-\frac12\int_0^{t-s_n} \frac{1}{B_{t'}^2}\,\mathrm{d}t'}\,|\partial_xw(s_n,B_{t-s_n})|\,\bone_{\{\tau_n^B=t-s_n\}}\Big] \\
&\le\EE^x\bigg[\Big(e^{\int_0^{t-s_n}\frac{1}{B_{t'}}\,\mathrm{d}B_{t'}
-\frac12\int_0^{t-s_n}\frac{1}{B_{t'}^2}\,\mathrm{d}t'}\Big)^3\,
\bone_{\{\tau_n^B=t-s_n\}}\bigg]^{\frac13}\,
\EE^x\big[|\partial_xw(s_n,B_{t-s_n})|^{3/2}\,\bone_{\{\tau_n^B=t-s_n\}}\big]^{\frac23},
\end{split}
\end{equation}
where $B$ is a standard Brownian motion with $B_0=x$ under $\PP^x$, and
\begin{IEEEeqnarray*}{rCl}
\tau_n^B:=\inf\{t'>0:\,B_{t'}\notin(\xi_{t-t'}^2,\xi_{t-t'}^1)\}\wedge(t-s_n).    
\end{IEEEeqnarray*}
The change of measure is licit thanks to $\min_{0\le t'\le t-s_n} R_{t'}\ge\min_{s_n\le s'\le t} \xi_{s'}^2>r_0$ on $\{\tau_n=t-s_n\}$.
	
\medskip
	
In view of $\mathrm{d}\log B_{t'}=\frac{1}{B_{t'}}\,\mathrm{d}B_{t'}-\frac{1}{2B_{t'}^2}\,\mathrm{d}t'$ on $\{\tau^B_n=t-s_n\}$,
\begin{IEEEeqnarray*}{rCl}
\EE^x\bigg[\Big(e^{\int_0^{t-s_n}\frac{1}{B_{t'}}\,\mathrm{d}B_{t'}
	-\frac12\int_0^{t-s_n}\frac{1}{B_{t'}^2}\,\mathrm{d}t'}\Big)^3\,
\bone_{\{\tau_n^B=t-s_n\}}\bigg] 
=\EE^x\big[e^{3\log B_{t-s_n}-3\log x}\,\bone_{\{\tau_n^B=t-s_n\}}\big] \\
=\EE^x\bigg[\frac{B_{t-s_n}^3}{x^3}\,\bone_{\{\tau_n^B=t-s_n\}}\bigg] 
\le\frac{\EE^x[|B_{t-s_n}|^3]}{x^3}\stackrel{n\to\infty}{\longrightarrow}
\frac{\EE^x[|B_{t-s_*}|^3]}{x^3}<\infty. 
\end{IEEEeqnarray*}
To bound the last expectation in \eqref{eq:partialxwestimate1}, we let $W:=B-x$, $M_{s'}:=\max_{0\le s''\le s'} W_{s''}$. Recall that
\begin{IEEEeqnarray*}{rCl}
\frac{2(2m-b)}{s'\sqrt{2\pi s'}}\,
e^{-\frac{(2m-b)^2}{2s'}}\,\bone_{\{m\ge0,\,b\leq m\}}
\end{IEEEeqnarray*}
is the density of $(M_s,W_s)$. Since $\xi_{s'}^1<\Lambda_{s'}\le\Lambda_{s_n}$, $s'\in[s_n,t]$, we have under $\PP^x$,
\begin{IEEEeqnarray*}{rCl}
\{\tau_n^B=t-s_n\}
\subset\{M_{t-s_n}\leq\Lambda_{s_n}-x,\,W_{t-s_n}\in[\xi_{s_n}^2-x,\xi_{s_n}^1-x]\}.    
\end{IEEEeqnarray*}
Given an $\varepsilon\in(0,((\Lambda_{s_*}-r_0)/2)\wedge\sqrt{s_*})$, it holds $\Lambda_{s_n}-\xi_{s_n}^2<\varepsilon$ and $\xi_{s_n}^2>x$ for all $n\in\NN$ large enough. Thus, by Proposition \ref{prop:GradientEstimateForw1}, there exists a constant $C<\infty$ such that
\begin{IEEEeqnarray*}{rCl}
|\partial_xw(s_n,x')|\le\frac{C}{\Lambda_{s_n}-x'},\quad x'\in[\xi_{s_n}^2,\xi_{s_n}^1].    
\end{IEEEeqnarray*}
Consequently, the last expectation in \eqref{eq:partialxwestimate1} is at most 
\begin{equation*}
\begin{split}
&\;C^{3/2}\,\EE^x\left[\frac{1}{(\Lambda_{s_n}-B_{t-s_n})^{3/2}}\,\bone_{\{\tau_n^B=t-s_n\}}\right] \\
&\le C^{3/2}\, \EE\left[\frac{1}{(\Lambda_{s_n}-x-W_{t-s_n})^{3/2}}\,\bone_{\{M_{t-s_n}\leq\Lambda_{s_n}-x,\,W_{t-s_n}\in[\xi_{s_n}^2-x,\xi_{s_n}^1-x]\}}\right] \\
&=C^{3/2} \int_{\xi_{s_n}^2-x}^{\xi_{s_n}^1-x}\int_b^{\Lambda_{s_n}-x} \frac{1}{(\Lambda_{s_n}-x-b)^{3/2}}\,\frac{2(2m-b)}{(t-s_n)\sqrt{2\pi(t-s_n)}}\,e^{-\frac{(2m-b)^2}{2(t-s_n)}}\,\mathrm{d}m\,\mathrm{d}b \\
&=C^{3/2}\int_{\xi_{s_n}^2-x}^{\xi_{s_n}^1-x}\int_b^{2(\Lambda_{s_n}-x)-b} \frac{1}{(\Lambda_{s_n}-x-b)^{3/2}}\,\frac{m}{(t-s_n)\sqrt{2\pi(t-s_n)}}\,e^{-\frac{m^2}{2(t-s_n)}}\,\mathrm{d}m\,\mathrm{d}b,
\end{split}
\end{equation*}
where in the last equality we made the change of variables $2m-b\to m$. 
	
\medskip	
	
By differentiating the function $b\mapsto be^{-\frac{b^2}{2}}$, it is easy to prove that
\begin{IEEEeqnarray*}{rCl}
\frac{m}{(t-s_n)\sqrt{2\pi(t-s_n)}}\,e^{-\frac{m^2}{2(t-s_n)}}
\le\frac{1}{(t-s_n)\sqrt{2\pi}}\,e^{-\frac12},\quad m>0.
\end{IEEEeqnarray*}
Putting everything together,
\begin{equation*}
\begin{split}
\EE^x\big[|\partial_xw(s_n,B_{t-s_n})|^{3/2}\,\bone_{\{\tau_n^B=t-s_n\}}\big]
&\le C^{3/2}\int_{\xi_{s_n}^2-x}^{\xi_{s_n}^1-x} \frac{2(\Lambda_{s_n}-x-b)}{(\Lambda_{s_n}-x-b)^{3/2}}\,\frac{1}{(t-s_n)\sqrt{2\pi}}\,e^{-\frac12}\,\mathrm{d}b \\
&=C^{3/2}\frac{4e^{-\frac12}}{(t-s_n)\sqrt{2\pi}}\big((\Lambda_{s_n}-\xi_{s_n}^2)^{1/2}
-(\Lambda_{s_n}-\xi_{s_n}^1)^{1/2}\big)\stackrel{n\to\infty}{\longrightarrow}0,
\end{split}
\end{equation*}
since $\lim_{n\to\infty}  \Lambda_{s_n}=\lim_{n\to\infty} \xi_{s_n}^2 = \lim_{n\to\infty} \xi_{s_n}^1=\Lambda_{s_*}$. It remains to take $n\to\infty$ in \eqref{eq:partialxwestimate1}. \qed

\bigskip\bigskip\bigskip

\end{document}